\documentclass[12pt]{article}

\usepackage{amsmath,amssymb,amsthm,amscd,a4wide,upref}
\usepackage[toc]{appendix}

\usepackage{hyperref}
\hypersetup{colorlinks,citecolor=blue,filecolor=black,linkcolor=blue,urlcolor=blue}
\usepackage{enumerate}
\usepackage{mathtools}

\begin{document}

\newcommand{\ad}{{\rm ad}}
\newcommand{\id}{{\rm id}}
\newcommand{\cri}{{\rm cri}}
\newcommand{\ext}{{\rm ext}}
\renewcommand{\exp}{{\rm exp}}
\newcommand{\End}{{\rm{End}\ts}}
\newcommand{\Rep}{{\rm{Rep}\ts}}
\newcommand{\Hom}{{\rm{Hom}}}
\newcommand{\Mat}{{\rm{Mat}}}
\newcommand{\ch}{{\rm{ch}\ts}}
\newcommand{\chara}{{\rm{char}\ts}}
\newcommand{\diag}{{\rm diag}}
\newcommand{\non}{\nonumber}
\newcommand{\wt}{\widetilde}
\newcommand{\wh}{\widehat}
\newcommand{\ot}{\otimes}
\newcommand{\la}{\lambda}
\newcommand{\La}{\Lambda}
\newcommand{\De}{\Delta}
\newcommand{\al}{\alpha}
\newcommand{\be}{\beta}
\newcommand{\ga}{\gamma}
\newcommand{\Ga}{\Gamma}
\newcommand{\ep}{\epsilon}
\newcommand{\ka}{\kappa}
\newcommand{\vk}{\varkappa}
\newcommand{\si}{\sigma}
\newcommand{\vs}{\varsigma}
\newcommand{\vp}{\varphi}
\newcommand{\de}{\delta}
\newcommand{\ze}{\zeta}
\newcommand{\om}{\omega}
\newcommand{\Om}{\Omega}
\newcommand{\ee}{\epsilon^{}}
\newcommand{\su}{s^{}}
\newcommand{\hra}{\hookrightarrow}
\newcommand{\ve}{\varepsilon}
\newcommand{\ts}{\,}
\newcommand{\pr}{^{\tss\prime}}
\newcommand{\vac}{\mathbf{1}}
\newcommand{\di}{\partial}
\newcommand{\qin}{q^{-1}}
\newcommand{\tss}{\hspace{1pt}}
\newcommand{\Sr}{ {\rm S}}
\newcommand{\U}{ {\rm U}}
\newcommand{\BL}{ {\overline L}}
\newcommand{\BE}{ {\overline E}}
\newcommand{\BP}{ {\overline P}}
\newcommand{\AAb}{\mathbb{A}\tss}
\newcommand{\CC}{\mathbb{C}\tss}
\newcommand{\KK}{\mathbb{K}\tss}
\newcommand{\QQ}{\mathbb{Q}\tss}
\newcommand{\SSb}{\mathbb{S}\tss}
\newcommand{\TT}{\mathbb{T}\tss}
\newcommand{\ZZ}{\mathbb{Z}\tss}
\newcommand{\DY}{ {\rm DY}}
\newcommand{\X}{ {\rm X}}
\newcommand{\Y}{ {\rm Y}}
\newcommand{\Z}{{\rm Z}}
\newcommand{\Ac}{\mathcal{A}}
\newcommand{\Lc}{\mathcal{L}}
\newcommand{\Mc}{\mathcal{M}}
\newcommand{\Pc}{\mathcal{P}}
\newcommand{\Qc}{\mathcal{Q}}
\newcommand{\Rc}{\mathcal{R}}
\newcommand{\Sc}{\mathcal{S}}
\newcommand{\Tc}{\mathcal{T}}
\newcommand{\Bc}{\mathcal{B}}
\newcommand{\Cc}{\mathcal{C}}
\newcommand{\Ec}{\mathcal{E}}
\newcommand{\Fc}{\mathcal{F}}
\newcommand{\Gc}{\mathcal{G}}
\newcommand{\Hc}{\mathcal{H}}
\newcommand{\Uc}{\mathcal{U}}
\newcommand{\Vc}{\mathcal{V}}
\newcommand{\Wc}{\mathcal{W}}
\newcommand{\Xc}{\mathcal{X}}
\newcommand{\Yc}{\mathcal{Y}}
\newcommand{\Ar}{{\rm A}}
\newcommand{\Br}{{\rm B}}
\newcommand{\Ir}{{\rm I}}
\newcommand{\Fr}{{\rm F}}
\newcommand{\Jr}{{\rm J}}
\newcommand{\Or}{{\rm O}}
\newcommand{\GL}{{\rm GL}}
\newcommand{\Spr}{{\rm Sp}}
\newcommand{\Rr}{{\rm R}}
\newcommand{\Zr}{{\rm Z}}
\newcommand{\ZX}{{\rm ZX}}
\newcommand{\gl}{\mathfrak{gl}}
\newcommand{\middd}{{\rm mid}}
\newcommand{\ev}{{\rm ev}}
\newcommand{\Pf}{{\rm Pf}}
\newcommand{\Norm}{{\rm Norm\tss}}
\newcommand{\oa}{\mathfrak{o}}
\newcommand{\spa}{\mathfrak{sp}}
\newcommand{\osp}{\mathfrak{osp}}
\newcommand{\e}{\mathfrak{e}}
\newcommand{\f}{\mathfrak{f}}
\newcommand{\g}{\mathfrak{g}}
\newcommand{\h}{\mathfrak h}
\newcommand{\n}{\mathfrak n}
\newcommand{\z}{\mathfrak{z}}
\newcommand{\Zgot}{\mathfrak{Z}}
\newcommand{\p}{\mathfrak{p}}
\newcommand{\sll}{\mathfrak{sl}}
\newcommand{\agot}{\mathfrak{a}}
\newcommand{\qdet}{ {\rm qdet}\ts}
\newcommand{\Ber}{ {\rm Ber}\ts}
\newcommand{\HC}{ {\mathcal HC}}
\newcommand{\cdet}{{\rm cdet}}
\newcommand{\rdet}{{\rm rdet}}
\newcommand{\tr}{ {\rm tr}}
\newcommand{\tra}{ {\rm t}}
\newcommand{\gr}{ {\rm gr}\ts}
\newcommand{\str}{ {\rm str}}
\newcommand{\loc}{{\rm loc}}
\newcommand{\Gr}{{\rm G}}
\newcommand{\sgn}{ {\rm sgn}\ts}
\newcommand{\sign}{{\rm sgn}}
\newcommand{\ba}{\bar{a}}
\newcommand{\bb}{\bar{b}}
\newcommand{\eb}{\bar{e}}
\newcommand{\bi}{\bar{\imath}}
\newcommand{\bj}{\bar{\jmath}}
\newcommand{\bk}{\bar{k}}
\newcommand{\bl}{\bar{l}}
\newcommand{\hb}{\mathbf{h}}
\newcommand{\Sym}{\mathfrak S}
\newcommand{\fand}{\quad\text{and}\quad}
\newcommand{\Fand}{\qquad\text{and}\qquad}
\newcommand{\For}{\qquad\text{or}\qquad}
\newcommand{\OR}{\qquad\text{or}\qquad}
\newcommand{\grpr}{{\rm gr}^{\tss\prime}\ts}
\newcommand{\degpr}{{\rm deg}^{\tss\prime}\tss}

\numberwithin{equation}{section}

\newtheorem{thm}{Theorem}[section]
\newtheorem{lem}[thm]{Lemma}
\newtheorem{prop}[thm]{Proposition}
\newtheorem{cor}[thm]{Corollary}
\newtheorem{conj}[thm]{Conjecture}
\newtheorem*{mthm}{Main Theorem}
\newtheorem*{mthma}{Theorem A}
\newtheorem*{mthmb}{Theorem B}
\newtheorem*{mthmc}{Theorem C}
\newtheorem*{mthmd}{Theorem D}

\theoremstyle{definition}
\newtheorem{defin}[thm]{Definition}

\theoremstyle{remark}
\newtheorem{remark}[thm]{Remark}
\newtheorem{example}[thm]{Example}

\newcommand{\bth}{\begin{thm}}
\renewcommand{\eth}{\end{thm}}
\newcommand{\bpr}{\begin{prop}}
\newcommand{\epr}{\end{prop}}
\newcommand{\ble}{\begin{lem}}
\newcommand{\ele}{\end{lem}}
\newcommand{\bco}{\begin{cor}}
\newcommand{\eco}{\end{cor}}
\newcommand{\bde}{\begin{defin}}
\newcommand{\ede}{\end{defin}}
\newcommand{\bex}{\begin{example}}
\newcommand{\eex}{\end{example}}
\newcommand{\bre}{\begin{remark}}
\newcommand{\ere}{\end{remark}}
\newcommand{\bcj}{\begin{conj}}
\newcommand{\ecj}{\end{conj}}

\newcommand{\bal}{\begin{aligned}}
\newcommand{\eal}{\end{aligned}}
\newcommand{\beq}{\begin{equation}}
\newcommand{\eeq}{\end{equation}}
\newcommand{\ben}{\begin{equation*}}
\newcommand{\een}{\end{equation*}}

\newcommand{\bpf}{\begin{proof}}
\newcommand{\epf}{\end{proof}}

\def\beql#1{\begin{equation}\label{#1}}

\title{\Large\bf Isomorphism between the $R$-matrix and Drinfeld presentations
of quantum affine algebra: type $C$}

\author{{Naihuan Jing,\quad Ming Liu\quad and\quad Alexander Molev}}

\date{} 
\maketitle


\begin{abstract}
An explicit isomorphism between the $R$-matrix and Drinfeld presentations
of the quantum affine algebra in type $A$ was given by Ding and I.~Frenkel (1993).
We show that this result can be extended to types $B$, $C$ and $D$
and give a detailed construction for type $C$ in this paper.
In all classical types the Gauss decomposition of the generator matrix
in the $R$-matrix presentation yields the Drinfeld generators.
To prove that the resulting map is an isomorphism we follow the work
of E.~Frenkel and Mukhin (2002) in type $A$ and employ the universal $R$-matrix
to construct the inverse map. A key role in our construction is played by a homomorphism theorem
which relates the quantum affine algebra of rank $n-1$ in the $R$-matrix
presentation with a subalgebra of the corresponding algebra of rank
$n$ of the same type.



\end{abstract}

\vspace{-0.7cm}


\tableofcontents

\section{Introduction}
\label{sec:int}

The quantum affine algebras $U_q(\wh\g)$
were introduced independently by Drinfeld~\cite{d:ha} and
Jimbo~\cite{j:qd} as deformations of the universal
enveloping algebras of the affine Lie algebras $\wh\g$ in the class of
Hopf algebras. Another presentation of these algebras
was given by Drinfeld~\cite{d:nr}, which is known as the {\em new realization} or
{\em Drinfeld presentation}. A detailed construction of the isomorphism between
the presentations was given by Beck~\cite{b:bg}. Yet another $R$-{\em matrix presentation}
of the quantum affine algebras was introduced by Reshetikhin and Semenov-Tian-Shansky
\cite{rs:ce} and further developed by I.~Frenkel and Reshetikhin~\cite{fri:qa}.

The algebras $U_q(\wh\g)$ possess a substantive algebraic structure
and rich representation theory. Their
finite-dimensional irreducible representations
were classified by Chari and Pressly in terms of the Drinfeld presentation; see~\cite[Chapter~12]{cp:gq}.
A theory of $q$-{\em characters} of these representations was developed in \cite{fm:cq} and \cite{fr:qc};
its connections with classical and quantum integrable systems were reviewed in
the expository paper \cite{kns:ts}.

The $R$-matrix presentation of the quantum affine algebras can also be used to describe
finite-dimensional irreducible representations by following the approach of Tarasov~\cite{t:im};
see also \cite{gm:rt}. The role of this presentation in the theory of
Knizhnik--Zamolodchikov equations is discussed in detail in the lectures~\cite{efk:lr}.

An explicit isomorphism between the Drinfeld and $R$-matrix presentations of the algebras $U_q(\wh\g)$
should provide a bridge between the two sides of the theory and widen the spectrum of methods
for their investigation. Such an isomorphism was already constructed in the case of
simple Lie algebras $\g$ of type $A$ by
Ding and I.~Frenkel~\cite{df:it}. We aim to extend
this result to the Lie algebras $\g$ of types $B$, $C$ and $D$. The present article is
concerned with type $C$, while types $B$ and $D$ will be dealt with in
a forthcoming paper.\footnote{The paper is now published [SIGMA {\bf 16} (2020), 043].
Note that the definition of the zero modes of the $L$-operators therein requires a correction:
the matrix $L^+[0]=\big[l_{ij}^+[0]\big]$
should be lower triangular, while $L^-[0]=\big[l_{ij}^-[0]\big]$ should be upper triangular.}

Our approach is similar to \cite{df:it}; it is based on the {\em Gauss decomposition} of the
generator matrices in the $R$-matrix presentation. The first part of the construction is
the verification that the generators arising from the Gauss decomposition do satisfy the required
relations of the Drinfeld presentation. The second part is the proof that the resulting
homomorphism is injective. We use an argument alternative to \cite{df:it} and follow
the work of E.~Frenkel and Mukhin~\cite{fm:ha} instead, where the map inverse to the
Ding--Frenkel isomorphism was constructed. This map relies on the formula for the universal
$R$-matrix corresponding to the algebra $U_q(\wh\g)$ due to
Khoroshkin and Tolstoy~\cite{kt:ur} and Damiani~\cite{d:rm}. It turns out to be possible
to use this formula in types $B$, $C$ and $D$ to construct a similar map in those cases.

Similar to the Yangian case in our previous work \cite{jlm:ib}, in this paper we
will mainly work with the {\em extended quantum affine algebra} in type $C$
defined by an $R$-matrix presentation. We prove
an embedding theorem which will allow us to regard
the extended algebra of rank $n-1$ as a subalgebra of the
corresponding algebra of rank $n$. We also
produce a Drinfeld-type presentation for the extended quantum affine algebra
and give explicit formulas for generators
of its center. Note that a different approach to the equivalences
between Yangian presentations and to the description of the centers of the extended Yangians
was developed in \cite{grw:eb} and \cite{w:rm} which should also be applicable
to quantum affine algebras.

To state our isomorphism theorem, choose simple roots for the
symplectic Lie algebra $\g=\spa_{2n}$ in the form
\ben
\al_i=\ep_i-\ep_{i+1}\quad\text{for}\quad i=1,\dots,n-1\Fand\al_n=2\epsilon_n,
\een
where
$\ep_1,\dots,\ep_n$ is an orthonormal basis of a Euclidian space with the bilinear form
$(\cdot\ts,\cdot)$. The Cartan matrix $[A_{ij}]$ is defined by
\beql{cartan}
A_{ij}=\frac{2(\al_i,\al_j)}{(\al_i,\al_i)}.
\eeq
For a variable $q$ we set $q_i=q^{r_i}$ for $i=1,\dots,n$, where $r_i=(\al_i,\al_i)/2$,
so that $q_i=q$ for $i<n$ and $q_n=q^2$.
We will use the standard notation
\beql{kq}
[k]_q=\frac{q^k-q^{-k}}{q-q^{-1}}
\eeq
for a nonnegative integer $k$, and
\ben
[k]_q!=\prod_{s=1}^{k}[s]_q,\qquad
{k\brack r}_{q}=\frac{[k]_q!}{[r]_q!\ts[k-r]_q!}.
\een

The {\em quantum affine algebra $U_q(\wh{\spa}_{2n})$ in its Drinfeld presentation}
is the associative algebra over $\CC(q)$
with generators
$x_{i,m}^{\pm}$, $a_{i,l}$, $k_{i}^{\pm}$ and $q^{\pm c/2}$ for $i=1,\dots,n$ and
$m,l\in\ZZ$ with $l\ne 0$, subject to the following defining relations:
the elements $q^{\pm c/2}$ are central,
\ben
\bal
  k_{i}k_i^{-1}&=k_i^{-1}k_{i}=1, \qquad q^{c/2}q^{-c/2}=q^{-c/2}q^{c/2}=1,\\
  k_ik_j&=k_jk_i, \qquad k_i\ts a_{j,k}=a_{j,k}\ts k_i,\qquad
  k_i\ts x_{j,m}^{\pm}\ts k_i^{-1}=q_i^{\pm A_{ij}}x_{j,m}^{\pm},
\eal
\een
\smallskip
\ben
\bal[]
  [a_{i,m},a_{j,l}]&=\delta_{m,-l}\ts
  \frac{[mA_{ij}]_{q_i}}{m}\ts\frac{q^{mc}-q^{-mc}}{q_j-q_j^{-1}},\\[0.4em]
  [a_{i,m}, x_{j,l}^{\pm}]&=\pm \frac{[mA_{ij}]_{q_i}}{m}\ts q^{\mp |m|c/2}\ts x^{\pm}_{j,m+l},
\eal
\een
\smallskip
\ben
\bal
  x^{\pm}_{i,m+1}x^{\pm}_{j,l}-q_i^{\pm A_{ij}}x^{\pm}_{j,l}x^{\pm}_{i,m+1}&=
  q_i^{\pm A_{ij}}x^{\pm}_{i,m}x^{\pm}_{j,l+1}-x^{\pm}_{j,l+1}x^{\pm}_{i,m},\\[0.5em]
  [x_{i,m}^{+},x_{j,l}^{-}]&=\delta_{ij}\frac{q^{(m-l)\ts c/2}\ts\psi_{i,m+l}
  -q^{-(m-l)\ts c/2}\ts\varphi_{i,m+l}}{q_i-q_i^{-1}},
\eal
\een
\ben
  \sum_{\pi\in \Sym_{r}}\sum_{l=0}^{r}(-1)^l{{r}\brack{l}}_{q_i}
  x^{\pm}_{i,s_{\pi(1)}}\dots x^{\pm}_{i,s_{\pi(l)}}
  x^{\pm}_{j,m}x^{\pm}_{i,s_{\pi(l+1)}}\dots x^{\pm}_{i,s_{\pi(r)}}=0, \quad i\neq j,
\een
where in the last relation we set $r=1-A_{ij}$. The elements
$\psi_{i,m}$ and $\varphi_{i,-m}$ with $m\in \ZZ_+$ are defined by
\begin{align}
\label{psiiu}
\psi_i(u)&:=\sum_{m=0}^{\infty}\psi_{i,m}u^{-m}=k_i\ts\exp\big((q_i-q_i^{-1})
\sum_{s=1}^{\infty}a_{i,s}u^{-s}\big),\\
\label{phiiu}
\varphi_{i}(u)&:=\sum_{m=0}^{\infty}\varphi_{i,-m}u^{m}=k_i^{-1}\exp\big({-}(q_i-q_i^{-1})
\sum_{s=1}^{\infty}a_{i,-s}u^{s}\big),
\end{align}
whereas $\psi_{i,m}=\varphi_{i,-m}=0$ for $m<0$.

To introduce the $R$-matrix presentation of the quantum affine algebra, consider the
following elements of the endomorphism algebra $\End(\CC^{2n}\ot\CC^{2n})
\cong \End\CC^{2n}\ot\End\CC^{2n}$:
\begin{align}
P&=\sum_{i,j=1}^{2n}e_{ij}\otimes e_{ji},\qquad
Q=\sum_{i,j=1}^{2n} q^{\bi-\bj}\ts\varepsilon_i\tss\varepsilon_j\ts e_{i'j'}\otimes e_{ij},
\non\\
\intertext{and}
R&=q\ts\sum_{i=1}^{2n}e_{ii}\otimes e_{ii}+\sum_{i\neq j,j'}e_{ii}\otimes e_{jj}
+q^{-1}\sum_{i\neq i'}e_{ii}\otimes e_{i'i'}
\non\\[0.3em]
{}&\qquad\qquad\qquad\qquad+(q-q^{-1})\sum_{i<j}e_{ij}\otimes
e_{ji}-(q-q^{-1})\sum_{i>j}q^{\bi-\bj}\ts\varepsilon_i\tss\varepsilon_j\ts e_{i'j'}\otimes e_{ij},
\non
\end{align}
where $e_{ij}\in\End\CC^{2n}$ are the matrix units and we used the notation
\ben
\ve_i=\begin{cases} \phantom{-}1\qquad\text{for}\quad i=1,\dots,n,\\
-1\qquad\text{for}\quad i=n+1,\dots,2n,
\end{cases}
\een
$i'=2n-i+1$ and
$
(\ts\overline{1},\overline{2},\dots,\overline{2n}\ts)=(n,n-1,\dots,1,-1,\dots,-n).
$
Furthermore, consider the formal power series
\ben
f(u)=1+\sum_{k=1}^{\infty}f_k\tss u^k
\een
whose coefficients $f_k$ are
rational functions in $q$ uniquely determined
by the relation
\beql{furel}
f(u)\tss f(u\xi)=\frac{1}{(1-u\tss q^{-2})(1-u\tss q^{2})(1-u\tss\xi)(1-u\tss\xi^{-1})},
\eeq
where $\xi=q^{-2n-2}$. Equivalently, $f(u)$ is given by
the infinite product formula
\beql{fu}
f(u)=\prod_{r=0}^{\infty}\frac{(1-u\ts\xi^{2r})(1-u\ts q^{-2}\ts\xi^{2r+1})
(1-u\ts q^2\ts\xi^{2r+1})(1-u\ts\xi^{2r+2})}
{(1-u\ts\xi^{2r-1})(1-u\ts\xi^{2r+1})
(1-u\ts q^2\xi^{2r})(1-u\ts q^{-2}\xi^{2r})}.
\eeq
Introduce the $R$-{\em matrix} $R(u)$ by
\beql{ru}
R(u)=f(u)\Big(q^{-1}(u-1)(u-\xi)R-(q^{-2}-1)(u-\xi)P+(q^{-2}-1)(u-1)\ts \xi\ts Q\Big).
\eeq
This formula goes back to Jimbo~\cite{j:qr}; for the significance of the scalar
function $f(u)$ see the paper by Frenkel and Reshetikhin~\cite{fri:qa}.
The $R$-matrix is a solution of the {\em Yang--Baxter equation}
\beql{ybe}
R_{12}(u)\ts R_{13}(u\tss v)\ts R_{23}(v)=R_{23}(v)\ts R_{13}(u\tss v)\ts R_{12}(u).
\eeq

The associative algebra $U^{R}_q(\wh{\spa}_{2n})$ over $\CC(q)$
is generated by an invertible central element $q^{c/2}$ and
elements  ${l}^{\pm}_{ij}[\mp m]$ with $1\leqslant i,j\leqslant 2n$ and $m\in \ZZ_{+}$
subject to the following defining relations. We have\footnote{The condition $i<j$
was erroneously replaced by the opposite inequality in the previous version.}
\ben
{l}_{ij}^{+}[0]={l}^{-}_{ji}[0]=0\quad\text{for}\quad i<j
\Fand {l}^{+}_{ii}[0]\ts {l}_{ii}^{-}[0]={l}_{ii}^{-}[0]\ts{l}^{+}_{ii}[0]=1,
\een
while the remaining relations will be written in terms of the formal power series
\beql{liju}
{l}^{\pm}_{ij}(u)=\sum_{m=0}^{\infty}{l}^{\pm}_{ij}[\mp m]\ts u^{\pm m}
\eeq
which we combine into the respective matrices
\ben
{L}^{\pm}(u)=\sum\limits_{i,j=1}^{2n}e_{ij}\ot {l}_{ij}^{\pm}(u)\in
\End\CC^{2n}\ot U^{R}_q(\wh{\spa}_{2n})[[u,u^{-1}]].
\een
Consider the tensor product algebra $\End\CC^{2n}\ot\End\CC^{2n}\ot U^{R}_q(\wh{\spa}_{2n})[[u,u^{-1}]]$
and introduce the series with coefficients in this algebra by
\beql{lonetwo}
{L}^{\pm}_1(u)=\sum\limits_{i,j=1}^{2n}e_{ij}\ot 1\ot {l}_{ij}^{\pm}(u)
\Fand
{L}^{\pm}_2(u)=\sum\limits_{i,j=1}^{2n}1\ot e_{ij}\ot {l}_{ij}^{\pm}(u).
\eeq
The defining relations then take the form
\begin{align}
\label{rllss}
{R}(u/v)L^{\pm}_{1}(u)L^{\pm}_2(v)&=L^{\pm}_2(v)L^{\pm}_{1}(u){R}(u/v),\\[0.4em]
{R}(uq^c/v)L^{+}_1(u)L^{-}_2(v)&=L^{-}_2(v)L^{+}_1(u){R}(uq^{-c}/v),
\label{rllmp}
\end{align}
together with the relations
\beql{unitary}
{L}^{\pm}(u)D{L}^{\pm}(u\ts\xi)^{\tra}D^{-1}=1,
\eeq
where $\tra$ denotes the matrix
transposition with $e_{ij}^{\tra}=\ve_i\ts\ve_j\ts e_{j',i'}$ and $D$ is the diagonal matrix
\beql{D}
D=\diag\ts [q^n,\dots,q,q^{-1},\dots,q^{-n}].
\eeq

Now apply the
{\em Gauss decomposition} to the matrices ${L}^{+}(u)$ and ${L}^{-}(u)$.
There exist unique matrices
of the form
\ben
F^{\pm}(u)=\begin{bmatrix}
1&0&\dots&0\ts\\
f^{\pm}_{21}(u)&1&\dots&0\\
\vdots&\vdots&\ddots&\vdots\\
f^{\pm}_{2n\ts1}(u)&f^{\pm}_{2n\ts2}(u)&\dots&1
\end{bmatrix},
\qquad
E^{\pm}(u)=\begin{bmatrix}
\ts1&e^{\pm}_{12}(u)&\dots&e^{\pm}_{1\ts 2n}(u)\ts\\
\ts0&1&\dots&e^{\pm}_{2\ts 2n}(u)\\
\vdots&\vdots&\ddots&\vdots\\
0&0&\dots&1
\end{bmatrix},
\een
and $H^{\pm}(u)=\diag\ts\big[h^{\pm}_1(u),\dots,h^{\pm}_{2n}(u)\big]$,
such that
\beql{gaussdec}
L^{\pm}(u)=F^{\pm}(u)H^{\pm}(u)E^{\pm}(u).
\eeq
For $i=1,\dots,n$ set
\ben
X^{+}_i(u)=e^{+}_{i,i+1}(u\tss q^{c/2})-e_{i,i+1}^{-}(u\tss q^{-c/2}),
\qquad X^{-}_i(u)=f^{+}_{i+1,i}(u\tss q^{-c/2})-f^{-}_{i+1,i}(u\tss q^{c/2}).
\een

To state our main result, combine the generators $x^{\pm}_{i,m}$ of the algebra
$U_q(\wh{\spa}_{2n})$ into the series
\beql{xpm}
x^{\pm}_{i}(u)=\sum_{m\in\ZZ}x^{\pm}_{i,m}\ts u^{-m}.
\eeq

\begin{mthm}
The maps $q^{c/2}\mapsto q^{c/2}$,
\ben
\bal
x^{\pm}_{i}(u)&\mapsto (q_i-q_i^{-1})^{-1}X^{\pm}_i(uq^i),\\[0.4em]
\psi_{i}(u)&\mapsto h^{-}_{i+1}(uq^i)\ts h^{-}_{i}(uq^i)^{-1},\\[0.4em]
\varphi_{i}(u)&\mapsto h^{+}_{i+1}(uq^i)\ts h^{+}_{i}(uq^i)^{-1},
\eal
\een
for $i=1,\dots,n-1$, and
\ben
\bal
x^{\pm}_{n}(u)&\mapsto (q_n-q_n^{-1})^{-1}X^{\pm}_n(uq^{n+1}),\\[0.4em]
\psi_{n}(u)&\mapsto h^{-}_{n+1}(uq^{n+1})\ts h^{-}_{n}(uq^{n+1})^{-1},\\[0.4em]
\varphi_{n}(u)&\mapsto h^{+}_{n+1}(uq^{n+1})\ts h^{+}_{n}(uq^{n+1})^{-1},
\eal
\een
define an isomorphism $U_q(\wh{\spa}_{2n})\to U^{R}_q(\wh{\spa}_{2n})$.
\end{mthm}

For the proof of the Main Theorem we embed $U_q(\wh{\spa}_{2n})$ into an extended
quantum affine algebra $U^{\ext}_q(\wh{\spa}_{2n})$ which is defined by
a Drinfeld-type presentation. The next step is to use the Gauss decomposition to
construct a homomorphism from
the extended quantum affine algebra to the
algebra $U(R)$ which is defined by the same
presentation as the algebra $U^{R}_q(\wh{\spa}_{2n})$, except that the relation
\eqref{unitary} is omitted. The expressions on the left hand side of \eqref{unitary},
considered in the algebra $U(R)$, turn out to be scalar matrices,
\ben
{L}^{\pm}(u)D{L}^{\pm}(u\xi)^{\tra}D^{-1}={z}^{\pm}(u)\ts 1,
\een
for certain formal series ${z}^{\pm}(u)$. Moreover, all coefficients of these series
are central in $U(R)$. We will give explicit formulas for ${z}^{\pm}(u)$,
regarded as series with coefficients in the algebra $U^{\ext}_q(\wh{\spa}_{2n})$,
in terms of its Drinfeld generators.
The quantum affine algebra $U_q(\wh{\spa}_{2n})$ can therefore be considered
as the quotient of $U^{\ext}_q(\wh{\spa}_{2n})$ by the relations ${z}^{\pm}(u)=1$.

As a final step, we construct the inverse map $U(R)\to U^{\ext}_q(\wh{\spa}_{2n})$
by using the universal $R$-matrix for the quantum affine algebra and producing
the associated $L$-operators corresponding to the vector representation of the
algebra $U_q(\wh{\spa}_{2n})$.

An immediate consequence of the Main Theorem is the Poincar\'e--Birkhoff--Witt
theorem for the $R$-matrix presentation $U^{R}_q(\wh{\spa}_{2n})$ of the
quantum affine algebra which is implied by the corresponding result
of Beck~\cite{b:cb} for $U_q(\wh{\spa}_{2n})$. As another application,
it is straightforward to transfer the classification theorem
for finite-dimensional irreducible representations of the algebra
$U_q(\wh{\spa}_{2n})$ to its $R$-matrix presentation $U^{R}_q(\wh{\spa}_{2n})$;
see \cite[Chapter~12]{cp:gq}.

\section{Quantum affine algebras}

We start by recalling the original definition of the quantum affine algebra
$U_q(\wh\g)$ as introduced by Drinfeld~\cite{d:ha} and Jimbo~\cite{j:qd}.
We suppose that $\g$ is a simple Lie algebra over $\CC$ of rank $n$
and $\wh\g$ is the corresponding (untwisted) affine Kac--Moody algebra with the affine
Cartan matrix $[A_{ij}]_{i,j=0}^n$. We let $\al_0,\al_1,\dots,\al_n$
denote the simple roots and use the notation as in \cite[Secs~9.1 and 12.2]{cp:gq}
so that $q_i=q^{r_i}$ for $r_i=(\al_i,\al_i)/2$.

\subsection{Drinfeld--Jimbo definition and new realization}
\label{subsec:isoDJD}

The {\em quantum affine algebra} $U_q(\wh{\g})$ is
a unital associative algebra over $\CC(q)$ with generators $E_{\al_i}$, $F_{\al_i}$ and
$k_i^{\pm 1}$ with $i=0,1,\dots,n$, subject to the defining relations:
\ben
k_ik_i^{-1}=k_i^{-1}k_i=1,\qquad k_ik_j=k_ik_j,
\een
\ben
k_iE_{\al_j}k_i^{-1}=q^{A_{ij}}_iE_{\al_j},\qquad k_iF_{\al_j}k_i^{-1}=q^{-A_{ij}}_iF_{\al_j},
\een
\ben
[E_{\al_i},F_{\al_j}]=\delta_{ij}\frac{k_i-k_i^{-1}}{q_i-q_i^{-1}},
\een
\ben
\bal
\sum_{r=0}^{1-A_{ij}}(-1)^r{{1-A_{ij}}\brack{r}}_{q_i}
(E_{\al_i})^rE_{\al_j}(E_{\al_i})^{1-A_{ij}-r}&=0,\qquad\text{if}\quad i\ne j,\\
\sum_{r=0}^{1-A_{ij}}(-1)^r{{1-A_{ij}}\brack{r}}_{q_i}
(F_{\al_i})^rF_{\al_j}(F_{\al_i})^{1-A_{ij}-r}&=0,\qquad\text{if}\quad i\ne j.
\eal
\een

By using the braid group action, the set of generators of the algebra $U_q(\wh\g)$
can be extended to the set of affine root vectors of the form $E_{\al+k\de}$, $F_{\al+k\de}$,
$E_{(k\de,i)}$ and $F_{(k\de,i)}$, where $\al$ runs over the positive roots of $\g$,
and $\de$ is the basic imaginary root; see \cite{b:bg, bcp:ac} for details.
The root vectors are used in the explicit isomorphism between the
Drinfeld--Jimbo presentation of the algebra $U_q(\wh\g)$
and the ``new realization" of Drinfeld which goes back to \cite{d:nr}, while detailed arguments
were given by Beck~\cite{b:bg}; see also \cite{bcp:ac}.
In particular, for the Drinfeld presentation of the algebra $U_q(\wh{\spa}_{2n})$
given in the Introduction, we find that
the isomorphism between these presentations is given by
\begin{alignat}{2}
\non
x^{+}_{ik}&\mapsto o(i)^kE_{\al_i+k\de}, \qquad &x^{-}_{i,-k}
&\mapsto o(i)^kF_{\al_i+k\de},\qquad\qquad  k\geqslant 0,\\
\non
x^{+}_{i,-k}&\mapsto -o(i)^kF_{-\al_i+k\de}\ts k_i^{-1}q^{kc},
\qquad &x^{-}_{i,k}
&\mapsto -o(i)^kq^{-kc}\ts k_i\ts E_{-\al_i+k\de},\quad  k> 0,\\
a_{i,k}&\mapsto o(i)^kq^{-kc/2}E_{(k\de,i)}, \qquad &a_{i,-k}
&\mapsto o(i)^k\tss F_{(k\de,i)}q^{kc/2},\qquad\qquad k> 0,
\non
\end{alignat}
where $o: \{1,2,\dots,n\}\rightarrow \{\pm 1\}$ is a map such that $o(i)=-o(j)$ whenever $A_{ij}<0$.

\subsection{Extended quantum affine algebra}
\label{subsec:eqaa}

As with the embedding of quantum affine algebras $U_q(\wh{\sll}_N)\hra U_q(\wh{\gl}_N)$
considered in \cite{df:it} and \cite{fm:ha}, it will be convenient to embed
the algebra $U_q(\wh{\spa}_{2n})$ into an extended quantum affine algebra which we denote by
$U^{\ext}_q(\wh{\spa}_{2n})$.

Beside the scalar function $f(u)$ defined by \eqref{furel} and \eqref{fu} we will also use the function
\beql{gu}
g(u)=f(u)(u-q^{-2})(u-\xi).
\eeq

\bde\label{def:eqaa}
The {\em extended quantum affine algebra} $U^{\ext}_q(\wh{\spa}_{2n})$ is an associative algebra
with generators $X^{\pm}_{i,k}$, $h^{+}_{j,m}$, $h^{-}_{j,-m}$ and $q^{c/2}$, where the subscripts take
values $i=1,\dots,n$ and $k\in\ZZ$, while $j=1,\dots,n+1$ and $m\in \ZZ_{+}$.
The defining relations are written with the use of generating functions in a formal
variable $u$:
\ben
X^{\pm}_i(u)=\sum_{k\in\ZZ}X^{\pm}_{i,k}\ts u^{-k},\qquad
h^{\pm}_i(u)=\sum_{m=0}^{\infty}h^{\pm}_{i,\pm m}\ts u^{\pm m},
\een
they take the following form. The element $q^{c/2}$ is
central and invertible,
\ben
h_{i,0}^{+}h_{i,0}^{-}=h_{i,0}^{-}h_{i,0}^{+}=1\Fand h^{+}_{n,0}h^{+}_{n+1,0}=1,
\een
for the relations involving $h^{\pm}_i(u)$ we have
\begin{align}\label{hihj}
h^{\pm}_i(u)h^{\pm}_j(v)&=h^{\pm}_j(v)h^{\pm}_i(u),\\[0.4em]
\label{hihipm}
g\big((uq^c/v)^{\pm 1}\big)\ts h^{\pm}_i(u)h^{\mp}_i(v)&
=g\big((uq^{-c}/v)^{\pm 1}\big)\ts h^{\mp}_i(v)h^{\pm}_i(u),\\
\label{hihjpm}
g\big((uq^c/v)^{\pm 1}\big)\ts\frac{u_{\pm}-v_{\mp}}{qu_{\pm}-q^{-1}v_{\mp}}h^{\pm}_i(u)h^{\mp}_j(v)&=
g\big((uq^{-c}/v)^{\pm 1}\big)\ts\frac{u_{\mp}-v_{\pm}}{qu_{\mp}-q^{-1}v_{\pm}}h^{\mp}_j(v)h^{\pm}_i(u)
\end{align}
for $i<j$ and $i\neq n$, and
\begin{multline}
\label{hnhn+1pm}
g\big((uq^c/v)^{\pm 1}\big)\ts
\frac{u_{\pm}-v_{\mp}}{q^2u_{\pm}-q^{-2}v_{\mp}}h^{\pm}_n(u)h^{\mp}_{n+1}(v)\\
{}=g\big((uq^{-c}/v)^{\pm 1}\big)\ts
\frac{u_{\mp}-v_{\pm}}{q^2u_{\mp}-q^{-2}v_{\pm}}h^{\mp}_{n+1}(v)h^{\pm}_n(u),
\end{multline}
where we use the notation $u_{\pm}=u\tss q^{\pm c/2}$. The relations
involving $h^{\pm}_i(u)$ and $X_{j}^{\pm}(v)$ are
\begin{align}\label{hiXjp}
h_{i}^{\pm}(u)X_{j}^{+}(v)&
=\frac{u-v_{\pm}}{q^{(\ep_i,\alpha_j)}u-q^{-(\ep_i,\alpha_j)}
v_{\pm}} X_{j}^{+}(v)h_{i}^{\pm}(u),\\
\label{hiXjm}
h_{i}^{\pm}(u)X_{j}^{-}(v)&
=\frac{q^{(\ep_i,\alpha_j)}u_{\pm}-q^{-(\ep_i,\alpha_j)}v}{u_{\pm}-v} X_{j}^{-}(v)
h_{i}^{\pm}(u),
\end{align}
for $i\neq n+1$, together with
\begin{align}\label{hn+1Xnp}
h_{n+1}^{\pm}(u)X_n^{+}(v)&=
\frac{u_{\mp}-v}{q^{-2}u_{\mp}-q^{2}v}X^{+}_n(v)h^{\pm}_{n+1}(u),\\
\label{hn+1Xnm}
h_{n+1}^{\pm}(u)X_n^{-}(v)&=
\frac{q^{-2}u_{\pm}-q^{2}v}{u_{\pm}-v}X_n^{-}(v)h_{n+1}^{\pm}(u),
\end{align}
and
\begin{align}\label{hn+1Xn-1p}
{h}_{n+1}^{\pm}(u)^{-1}{X}_{n-1}^{+}(v){h}_{n+1}^{\pm}(u)
&=\frac{q^{-1}u-qv_{\pm}}{q^{-2}u-q^2v_{\pm}}{X}_{n-1}^{+}(v),\\
\label{hn+1Xn-1m}
{h}_{n+1}^{\pm}(u){X}_{n-1}^{-}(v){h}_{n+1}^{\pm}(u)^{-1}
&=\frac{q^{-1}u-qv_{\mp}}{q^{-2}u-q^2v_{\mp}}{X}_{n-1}^{-}(v),
\end{align}
while
\begin{align}\label{hn+1Xip}
{h}_{n+1}^{\pm}(u){X}_{i}^{+}(v)
&={X}_{i}^{+}(v)\ts {h}_{n+1}^{\pm}(u),\\[0.4em]
\label{hn+1Xim}
h_{n+1}^{\pm}(u){X}_{i}^{-}(v)
&={X}_{i}^{-}(v)\ts {h}_{n+1}^{\pm}(u),
\end{align}
for $1\leqslant i\leqslant n-2$. For the relations involving $X^{\pm}_i(u)$ we have
\ben
(u-q^{\pm (\alpha_i,\alpha_j)}v)X_{i}^{\pm}(uq^i)X_{j}^{\pm}(vq^j)
=(q^{\pm (\alpha_i,\alpha_j)}u-v) X_{j}^{\pm}(vq^j)X_{i}^{\pm}(uq^i)
\een
for $i,j=1,\dots, n-1$,
\ben
(u-q^{\pm (\alpha_i,\alpha_n)}v)X_{i}^{\pm}(uq^i)X_{n}^{\pm}(vq^{n+1})
=(q^{\pm (\alpha_i,\alpha_n)}u-v) X_{n}^{\pm}(vq^{n+1})X_{i}^{\pm}(uq^i)
\een
for $i=1,\dots,n-1$,
\ben
(u-q^{\pm (\alpha_n,\alpha_n)}v)X_{n}^{\pm}(u)X_{n}^{\pm}(v)
=(q^{\pm (\alpha_n,\alpha_n)}u-v) X_{n}^{\pm}(v)X_{n}^{\pm}(u),
\een
and
\begin{multline}
\non
[X_i^{+}(u),X_j^{-}(v)]=\delta_{ij}(q_i-q_i^{-1})\\[0.4em]
{}\times\Big(\delta({u}q^{-c}/v)h_i^{-}(v_+)^{-1}h_{i+1}^{-}(v_+)
-\delta({u}q^{c}/{v})h_i^{+}(u_+)^{-1}h_{i+1}^{+}(u_+)\Big),
\end{multline}
together with the
{\em Serre relations}
\beql{serre}
\sum_{\pi\in \Sym_{r}}\sum_{l=0}^{r}(-1)^l{{r}\brack{l}}_{q_i}
  X^{\pm}_{i}(u_{\pi(1)})\dots X^{\pm}_{i}(u_{\pi(l)})
  X^{\pm}_{j}(v)\tss X^{\pm}_{i}(u_{\pi(l+1)})\dots X^{\pm}_{i}(u_{\pi(r)})=0,
\eeq
which hold for all $i\neq j$ and we set $r=1-A_{ij}$. Here we used the notation
\beql{delta-f}
\de(u)=\sum_{r\in\ZZ}\tss u^r
\eeq
for the {\em formal $\delta$-function}.
\qed
\ede

Introduce two formal power series $z^+(u)$ and $z^-(u)$ in $u$ and $u^{-1}$, respectively,
with coefficients in the algebra $U^{\ext}_q(\wh{\spa}_{2n})$ by
\beql{zpm}
z^{\pm}(u)=\prod_{i=1}^{n-1}h^{\pm}_{i}(u\tss\xi\tss q^{2i})^{-1}
\prod_{i=1}^{n}h^{\pm}_{i}(u\tss\xi\tss q^{2i-2})h^{\pm}_{n+1}(u),
\eeq
where we keep using the notation $\xi=q^{-2n-2}$.
Note that by \eqref{hihj} the ordering of the factors in the products is irrelevant.

\bpr\label{prop:zu}
The coefficients of $z^{\pm}(u)$ are central elements of $U^{\ext}_q(\wh{\spa}_{2n})$.
\epr

\bpf
We will outline
the arguments for $z^+(u)$; the case of
$z^-(u)$ is quite similar. By \eqref{hihj} we obviously have
$z^{+}(u)\ts h_j^{+}(v)=h_j^{+}(v)\ts z^{+}(u)$ for all $j=1,\dots,n+1$.
It is straightforward to deduce
from the defining relations in Definition~\ref{def:eqaa}
that $z^{+}(u)\ts X_j^{\pm}(v)=X_j^{\pm}(v)\ts z^{+}(u)$ for $j=1,\dots,n$.
We give more details to check that
$z^{+}(u)\ts h_{n+1}^{-}(v)=h_{n+1}^{-}(v)\ts z^{+}(u)$
as this involves the function \eqref{gu}.
The remaining calculations are performed in the same way.
Applying \eqref{hihipm}
we get
\ben
\bal
z^{+}(u)h^{-}_{n+1}(v)&=g(uq^{-c}/v)\ts g(uq^c/v)^{-1}\\
&\times \prod_{i=1}^{n-1}h^{+}_{i}(u\xi q^{2i})^{-1}
\prod_{i=1}^{n-1}h^{+}_{i}(u\xi q^{2i-2})\ts
h_{n}^{+}(uq^{-4})h^{-}_{n+1}(v)h_{n+1}^{+}(u).
\eal
\een
Furthermore,
\eqref{hnhn+1pm} implies
\begin{multline}
\non
z^{+}(u)h^{-}_{n+1}(v)=g(uq^{-c-4}/v)\ts g(uq^{c-4}/v)^{-1}\ts g(uq^{-c}/v)\ts g(uq^c/v)^{-1}\\[0.4em]
{}\times \frac{q^{-2}u_{-}-q^2v_{+}}{u_{-}-v_{+}}\frac{u_{+}-v_{-}}{q^{-2}u_{+}-q^2v_{-}}\\
{}\times\prod_{i=1}^{n-1}h^{+}_{i}(u\xi q^{2i})^{-1}
\prod_{i=1}^{n-1}h^{+}_{i}(u\xi q^{2i-2})
h^{-}_{n+1}(v)h_{n}^{+}(uq^{-4})h_{n+1}^{+}(u).
\end{multline}
Due to \eqref{hihj}, the last product can be rearranged as
\ben
\prod_{i=1}^{n-2}h^{+}_{i}(u\xi q^{2i})^{-1}h^{+}_{i+1}(u\xi q^{2i})
h^{+}_{n-1}(uq^{-4})^{-1}h^{+}_{1}(u\xi)
h^{-}_{n+1}(v)h_{n}^{+}(uq^{-4})h_{n+1}^{+}(u).
\een
Now, applying \eqref{hihjpm} repeatedly, we come to the relation
\ben
\bal
z^{+}(u)h^{-}_{n+1}(v)&=
g(u\xi q^{-c}/v)\ts g(u\xi q^c/v)^{-1}\ts g(uq^{-c}/v)\ts g(uq^c/v)^{-1}\\[0.5em]
&\times \frac{u_{-}\xi-v_{+}}{u_{-}\xi q-v_{+}q^{-1}}\frac{u_{+}\xi q-v_{-}q^{-1}}{u_{+}\xi-v_{-}}
 \frac{q^{-1}u_{-}-qv_{+}}{u_{-}-v_{+}}\frac{u_{+}-v_{-}}{q^{-1}u_{+}-qv_{-}}
 \ts h^{-}_{n+1}(v)z^{+}(u).
\eal
\een
Replace $g(u)$ by \eqref{gu} to get
\ben
\bal
z^{+}(u)h^{-}_{n+1}(v)&=f(u\xi q^{-c}/v)\ts f(u\xi q^c/v)^{-1}\ts f(uq^{-c}/v)\ts f(uq^c/v)^{-1}\\[0.5em]
&\times \frac{(u_{-}/v_{+}-q^{-2})(u_{-}/v_{+}-\xi)(u_{-}/v_{+}-q^{2})(u_{-}/v_{+}\xi-1)}
{(u_{+}/v_{-}-q^{-2})(u_{+}/v_{-}-\xi)(u_{+}/v_{-}-q^{2})(u_{+}/v_{-}\xi-1)}\ts
h^{-}_{n+1}(v)z^{+}(u).
\eal
\een
Since
\ben
f(u)f(u\xi)=\frac{1}{(1-uq^{-2})(1-uq^{2})(1-u\xi)(1-u\xi^{-1})},
\een
we can conclude that $z^{+}(u)h^{-}_{n+1}(v)=h^{-}_{n+1}(v)z^{+}(u)$.
\epf

\bpr\label{prop:embed}
The maps $q^{c/2}\mapsto q^{c/2}$,
\ben
\bal
x^{\pm}_{i}(u)&\mapsto (q_i-q_i^{-1})^{-1}X^{\pm}_i(uq^i),\\[0.4em]
\psi_{i}(u)&\mapsto h^{-}_{i+1}(uq^i)\ts h^{-}_{i}(uq^i)^{-1},\\[0.4em]
\varphi_{i}(u)&\mapsto h^{+}_{i+1}(uq^i)\ts h^{+}_{i}(uq^i)^{-1},\\[0.4em]
k_i &\mapsto h^{+}_{i,0}(h^{+}_{i+1,0})^{-1}
\eal
\een
for $i=1,\dots,n-1$, and
\ben
\bal
x^{\pm}_{n}(u)&\mapsto (q_n-q_n^{-1})^{-1}X^{\pm}_n(uq^{n+1}),\\[0.4em]
\psi_{n}(u)&\mapsto h^{-}_{n+1}(uq^{n+1})\ts h^{-}_{n}(uq^{n+1})^{-1},\\[0.4em]
\varphi_{n}(u)&\mapsto h^{+}_{n+1}(uq^{n+1})\ts h^{+}_{n}(uq^{n+1})^{-1},\\[0.4em]
k_n &\mapsto (h^{+}_{n,0})^2,
\eal
\een
define an embedding $\varsigma:U_q(\wh{\spa}_{2n})\hra U^{\ext}_q(\wh{\spa}_{2n})$.
\epr

\bpf
Writing the defining relations of the quantum affine algebra
$U_q(\wh{\spa}_{2n})$ in terms of the generating series $x^{\pm}_{i}(u)$,
$\psi_{i}(u)$ and $\varphi_{i}(u)$, it is straightforward to check that
the maps define a homomorphism. To show that its kernel
is zero, we will construct another
homomorphism $\varrho:U^{\ext}_q(\wh{\spa}_{2n})\to U_q(\wh{\spa}_{2n})$ such that
the composition $\varrho\circ\varsigma$ is the identity homomorphism on $U_q(\wh{\spa}_{2n})$.
We will extend $U_q(\wh{\spa}_{2n})$ by adjoining the square roots $k_n^{\pm 1/2}$
and keep the same notation
for the extended algebra for the rest of the argument. There exist power series
$\ze^{\pm}(u)$ with coefficients in the center of $U^{\ext}_q(\wh{\spa}_{2n})$
such that $\ze^{\pm}(u)\ts\ze^{\pm}(u\tss\xi)=z^{\pm}(u)$. Explicitly,
\ben
\ze^{\pm}(u)=\prod_{m=0}^{\infty} z^{\pm}(u\xi^{-2m-1})z^{\pm}(u\xi^{-2m-2})^{-1}.
\een
The mappings
$X^{\pm}_i(u)\mapsto X^{\pm}_i(u)$ for $i=1,\dots, n$ and
$h^{\pm}_j(u)\mapsto h^{\pm}_j(u)\ts \ze^{\pm}(u)$
for $j=1,\dots, n+1$
define a homomorphism from the algebra $U^{\ext}_q(\wh{\spa}_{2n})$ to itself.
The definition of the series $\ze^{\pm}(u)$ implies that for images of $h^{\pm}_i(u)$
we have the relation
\ben
h^{\pm}_i(u)\ts \ze^{\pm}(u)\ts h^{\pm}_i(u\xi)\ts \ze^{\pm}(u\xi)=
h^{\pm}_i(u)\ts h^{\pm}_i(u\xi)\ts z^{\pm}(u).
\een
Hence the property $\varrho\circ\varsigma={\rm id}$
will be satisfied if we define the map $\varrho:U^{\ext}_q(\wh{\spa}_{2n})\to U_q(\wh{\spa}_{2n})$ by
\ben
X^{\pm}_i(u)\mapsto
(q_i-q_i^{-1})\ts x^{\pm}_{i}(uq^{-i})\qquad\text{for}\quad i=1,\dots,n-1,
\een
and
\ben
X^{\pm}_n(u)\mapsto
(q_n-q_n^{-1})\ts x^{\pm}_{n}(uq^{-n-1}),
\een
while
\ben
h^{\pm}_i(u)\mapsto \al^{\pm}_i(u)\qquad\text{for}\quad i=1,\dots,n+1,
\een
where the series $\al^{+}_i(u)$ are defined by the relations
\ben
\al^{+}_i(u)\ts \al^{+}_i(u\tss\xi)=\vp_{n}(u\tss q^{-n-1})^{-1}
\ts\prod_{k=1}^{n-1}\vp_{k}(u\tss\xi\tss  q^k)^{-1}
\ts\prod_{k=1}^{i-1}\vp_{k}(u\tss\xi\tss  q^{-k})\ts\prod_{k=i}^{n-1}\vp_{k}(u\tss q^{-k})^{-1}
\een
for $i=1,\dots,n$, and
\ben
\al^{+}_{n+1}(u)\ts \al^{+}_{n+1}(u\tss\xi)=\vp_{n}(u\tss\xi\tss q^{-n-1})
\ts\prod_{k=1}^{n-1}\vp_{k}(u\tss\xi\tss  q^k)^{-1}
\ts\prod_{k=1}^{n-1}\vp_{k}(u\tss\xi\tss  q^{-k}).
\een
Explicitly, by setting $\tilde{\varphi}_{j}(u)=k_j\tss\varphi_{j}(u)$, we have
\begin{multline}
\non
\al^{+}_1(u)=\prod_{m=0}^{\infty} \prod_{j=1}^{n-1}\tilde{\varphi}_{j}(u\xi^{-2m}q^j)^{-1}
\ts\tilde{\varphi}_{j}(u\xi^{-2m-1}q^j)\ts
\tilde{\varphi}_{j}(u\xi^{-2m-1}q^{-j})^{-1}\ts\tilde{\varphi}_{j}(u\xi^{-2m-2}q^{-j})\\
{}\times \prod_{m=0}^{\infty} \tilde{\varphi}_{n}(u\xi^{-2m}q^{n+1})^{-1}
\tilde{\varphi}_{n}(u\xi^{-2m-1}q^{n+1})\times \prod_{j=1}^{n-1}k_ik^{1/2}_n,
\end{multline}
\begin{multline}
\non
\al^{+}_i(u)= \prod_{m=0}^{\infty} \prod_{j=1}^{n-1}\tilde{\varphi}_{j}(u\xi^{-2m}q^j)^{-1}
\ts\tilde{\varphi}_{j}(u\xi^{-2m-1}q^j)\ts
\tilde{\varphi}_{j}(u\xi^{-2m-1}q^{-j})^{-1}\tilde{\varphi}_{j}(u\xi^{-2m-2}q^{-j})\\
{}\times \prod_{m=0}^{\infty} \tilde{\varphi}_{n}(u\xi^{-2m}q^{n+1})^{-1}\ts
\tilde{\varphi}_{n}(u\xi^{-2m-1}q^{n+1})\times \prod_{j=1}^{i-1}\tilde{\varphi}_{j}(uq^{-j})
\times \prod_{j=i}^{n-1}k_ik^{1/2}_n,
\end{multline}
for $i=2,\dots, n-1$,
\begin{multline}
\non
\al^{+}_n(u)= \prod_{m=0}^{\infty} \prod_{j=1}^{n-1}\tilde{\varphi}_{j}(u\xi^{-2m}q^j)^{-1}
\ts\tilde{\varphi}_{j}(u\xi^{-2m-1}q^j)\ts
\tilde{\varphi}_{j}(u\xi^{-2m-1}q^{-j})^{-1}\tilde{\varphi}_{j}(u\xi^{-2m-2}q^{-j})\\
{}\times \prod_{m=0}^{\infty} \tilde{\varphi}_{n}(u\xi^{-2m}q^{n+1})^{-1}\ts
\tilde{\varphi}_{n}(u\xi^{-2m-1}q^{n+1})\times \prod_{j=1}^{n-1}
\tilde{\varphi}_{j}(uq^{-j})\times k^{1/2}_n
\end{multline}
and
\begin{multline}
\non
\al^{+}_{n+1}(u)= \prod_{m=0}^{\infty} \prod_{j=1}^{n-1}
\tilde{\varphi}_{j}(u\xi^{-2m}q^j)^{-1}\ts
\tilde{\varphi}_{j}(u\xi^{-2m-1}q^j)\ts
\tilde{\varphi}_{j}(u\xi^{-2m-1}q^{-j})^{-1}\tilde{\varphi}_{j}(u\xi^{-2m-2}q^{-j})\\
{}\times \prod_{m=0}^{\infty} \tilde{\varphi}_{n}(u\xi^{-2m}q^{n+1})^{-1}\ts
\tilde{\varphi}_{n}(u\xi^{-2m-1}q^{n+1})\times \prod_{j=1}^{n-1}
\tilde{\varphi}_{j}(uq^{-j})\ts\tilde{\varphi}_{n}(uq^{-(n+1)})\times k^{-1/2}_n.
\end{multline}
The relations defining $\al^{-}_i(u)$ are obtained from those above
by the respective replacements $\al^{+}_i(u)\to \al^{-}_i(u)$, $k_i\to k_i^{-1}$
and $\vp_{k}(u)\to\psi_k(u)$.

As with the map $\varsigma$, it is possible to verify directly that
the map $\varrho$ defines a homomorphism. Alternatively, this fact also follows
from calculations with Gaussian generators; see Remark~\ref{rem:prop2.3} below.
\epf

By Proposition~\ref{prop:embed}, we may regard $U_q(\wh{\spa}_{2n})$ as a subalgebra
of the extended quantum affine algebra $U^{\ext}_q(\wh{\spa}_{2n})$.
As in the proof of the proposition,
we will assume that $U_q(\wh{\spa}_{2n})$ is extended by adjoining the square roots $k_n^{\pm 1/2}$,
and let $\Cc$ be the subalgebra of $U^{\ext}_q(\wh{\spa}_{2n})$
generated by the coefficients
of the series $z^{\pm}(u)$.

\bco\label{cor:decomp}
We have the tensor product decomposition
\ben
U^{\ext}_q(\wh{\spa}_{2n}) = U_q(\wh{\spa}_{2n})\otimes \Cc.
\een
\eco

\bpf
From the proof of Proposition \ref{prop:embed}, we can define a
surjective homomorphism from $U^{\ext}_q(\wh{\spa}_{2n})$ to $U_q(\wh{\spa}_{2n})\otimes \Cc$ by
\ben
\bal
X^{\pm}_i(u)&\mapsto
(q_i-q_i^{-1})\ts x^{\pm}_{i}(uq^{-i})\ot 1\qquad\qquad\text{for}\quad i=1,\dots,n-1,\\[0.4em]
X^{\pm}_n(u)&\mapsto
(q_n-q_n^{-1})\ts x^{\pm}_{n}(uq^{-n-1})\ot 1,
\eal
\een
and
\ben
h^{\pm}_i(u)\mapsto \varrho\big(h^{\pm}_i(u)\big)\ot \ze^{\pm}(u)
\een
for $i=1,\dots,n+1$. The inverse map
$U_q(\wh{\spa}_{2n})\otimes \Cc\to U^{\ext}_q(\wh{\spa}_{2n})$
is defined by extending the homomorphism $\varsigma$ so that
$1\ot\ze^{\pm}(u)\mapsto \ze^{\pm}(u)$.
\epf

\section{$R$-matrix presentations}
\label{sec:nd}


\subsection{The algebras $U(R)$ and $U(\overline{R})$}

Recall from the Introduction that the algebra $U(R)$
is generated by an invertible central element $q^{c/2}$ and
elements  ${l}^{\pm}_{ij}[\mp m]$ with $1\leqslant i,j\leqslant 2n$ and $m\in \ZZ_{+}$
such that
\ben
{l}_{ij}^{+}[0]={l}^{-}_{ji}[0]=0\quad\text{for}\quad i<j,
\qquad {l}^{+}_{n,n}[0]\ts {l}_{n+1,n+1}^{+}[0]=1
\fand {l}^{+}_{ii}[0]\ts {l}_{ii}^{-}[0]={l}_{ii}^{-}[0]\ts{l}^{+}_{ii}[0]=1,
\een
and the remaining relations \eqref{rllss} and \eqref{rllmp} (omitting \eqref{unitary})
written in terms of the formal power series \eqref{liju}.
We will need another algebra $U(\overline{R})$ which is defined in a very similar way,
except that it is associated with a different $R$-matrix $\overline{R}(u)$ instead of \eqref{ru}.
Namely, the two $R$-matrices are related by $R(u)=g(u)\overline{R}(u)$
with $g(u)$ defined in \eqref{gu}, so that
\beql{rbar}
\overline{R}(u)=\frac{u-1}{u\tss q-q^{-1}}\ts R
+\frac{q-q^{-1}}{u\tss q-q^{-1}}\ts P-\frac{(q-q^{-1})(u-1)\xi}{(u\tss q-q^{-1})(u-\xi)}\ts Q.
\eeq
Note the {\em unitarity property}
\beql{unitarityrbar}
\overline{R}_{12}(u)\ts \overline{R}_{21}(u^{-1})=1,
\eeq
satisfied by
this $R$-matrix, where $\overline{R}_{12}(u)=\overline{R}(u)$ and
$\overline{R}_{21}(u)=P\overline{R}(u)P$.

The {\em algebra $U(\overline{R})$ over $\CC(q)$}
is generated by an invertible central element $q^{c/2}$ and
elements  ${\ell}^{\tss\pm}_{ij}[\mp m]$ with $1\leqslant i,j\leqslant 2n$ and $m\in \ZZ_{+}$
such that
\ben
{\ell}_{ij}^{+}[0]={\ell}^{-}_{ji}[0]=0\quad\text{for}
\quad i<j,\quad{\ell}^{+}_{n,n}[0]\ts {\ell}_{n+1,n+1}^{+}[0]=1
\fand {\ell}^{+}_{ii}[0]\ts {\ell}_{ii}^{-}[0]={\ell}_{ii}^{-}[0]\ts{\ell}^{+}_{ii}[0]=1.
\een
Introduce the formal power series
\beql{lijubar}
{\ell}^{\tss\pm}_{ij}(u)=\sum_{m=0}^{\infty}{\ell}^{\tss\pm}_{ij}[\mp m]\ts u^{\pm m}
\eeq
which we combine into the respective matrices
\ben
{\Lc}^{\pm}(u)=\sum\limits_{i,j=1}^{2n}e_{ij}\ot {\ell}^{\tss\pm}_{ij}(u)\in
\End\CC^{2n}\ot U(\overline{R})[[u,u^{-1}]].
\een
The remaining defining relations of the algebra $U(\overline{R})$ take the form
\begin{align}
\label{gen rel1}
\overline{R}(u/v)\ts\Lc^{\pm}_{1}(u)\ts\Lc^{\pm}_2(v)&
=\Lc^{\pm}_2(v)\ts\Lc^{\pm}_{1}(u)\ts\overline{R}(u/v),\\[0.4em]
\overline{R}(u\tss q^c/v)\ts\Lc^{+}_1(u)\ts\Lc^{-}_2(v)&
=\Lc^{-}_2(v)\ts\Lc^{+}_1(u)\ts\overline{R}(u\tss q^{-c}/v),
\label{gen rel2}
\end{align}
where the subscripts have the same meaning as in \eqref{lonetwo}.
The unitarity property \eqref{unitarityrbar} implies that
relation \eqref{gen rel2} can be written in the equivalent form
\beql{gen rel3}
\overline{R}(u\tss q^{-c}/v)\ts\Lc^{-}_1(u)\ts\Lc^{+}_2(v)
=\Lc^{+}_2(v)\ts\Lc^{-}_1(u)\ts\overline{R}(u\tss q^{c}/v).
\eeq

\bre\label{rem:gln}
The defining relations satisfied by the series
$\ell^{\tss\pm}_{ij}(u)$
with $1\leqslant i,j\leqslant n$ coincide with those for the quantum affine algebra
$U_q(\wh{\gl}_n)$ in \cite{df:it}.
\qed
\ere

Now we will follow \cite{df:it} to describe a relationship between the algebras
$U(R)$ and $U(\overline{R})$.
Introduce a Heisenberg algebra $\Hc_q(n)$
with generators $q^c$ and $\be_r$ with $r\in\ZZ\setminus \{0\}$.
The defining relations of $\Hc_q(n)$ have the form
\ben
\big[\be_r,\be_s\big]=\de_{r, -s}\ts \al_r,\qquad r\geqslant 1,
\een
and $q^c$ is central and invertible.
The elements $\al_r$ are defined by the expansion
\ben
\exp\ts\sum_{r=1}^{\infty}\al_r u^r=\frac{g(u\tss q^{-c})}{g(u\tss q^{c})}.
\een
So we have the identity
\ben
g(u\ts q^{c}/v)\ts \exp\ts\sum_{r=1}^{\infty}\be_r u^r\cdot
\ts \exp\ts\sum_{s=1}^{\infty}\be_{-s} v^{-s}
=g(u\ts q^{-c}/v)\ts \exp\ts\sum_{s=1}^{\infty}\be_{-s} v^{-s}
\cdot \exp\ts\sum_{r=1}^{\infty}\be_r u^r.
\een

\bpr\label{prop:homheis}
The mappings
\beql{heisext}
\Lc^{\ts +}(u)\mapsto \exp\ts\sum_{r=1}^{\infty}\be_{-r} u^{-r}\cdot
L^+(u),\qquad
\Lc^{\ts -}(u)\mapsto \exp\ts\sum_{r=1}^{\infty}\be_{r} u^{r}\cdot
L^-(u),
\eeq
define a homomorphism $U(\overline{R})\to \Hc_q(n)\ot_{\CC[q^c,\ts q^{-c}]}U(R)$.
\qed
\epr

We will need to apply the matrix transposition defined in
\eqref{unitary} to certain copies of the endomorphism algebra
$\End\CC^{2n}$ in multiple tensor products.
The corresponding partial transposition applied to the $a$-th copy will be denoted by $\tra_a$.
We point out the following {\em crossing symmetry} relations
satisfied by the $R$-matrices:
\begin{align}\label{RDRD}
 \overline{R}(u)D_1\overline{R}(u\tss\xi)^{\tra_1}D_1^{-1}
 &=\frac{(u-q^2)(u\tss\xi-1)}{(1-u)(1-u\tss\xi q^2)},\\[0.4em]
 {R}(u)D_1{R}(u\tss\xi)^{\tra_1}D_1^{-1}&=\xi^2 q^{-2},
 \label{crsymr}
\end{align}
where the diagonal matrix $D$ is defined in \eqref{D} and the meaning of the subscripts
is the same as in \eqref{lonetwo}.

\bpr\label{prop:central}
In the algebras $U(R)$ and $U(\overline{R})$ we have
the relations
\beql{DLDL}
D{L}^{\pm}(u\tss\xi)^{\tra}D^{-1}{L}^{\pm}(u)
={L}^{\pm}(u)D{L}^{\pm}(u\tss\xi)^{\tra}D^{-1}={z}^{\pm}(u)\ts1,
\eeq
and
\beql{DLbarDLbar}
D\Lc^{\pm}(u\tss\xi)^{\tra}D^{-1}\Lc^{\pm}(u)
=\Lc^{\pm}(u)D\Lc^{\pm}(u\tss\xi)^{\tra}D^{-1}=\z^{\pm}(u)\ts1,
\eeq
for certain series ${z}^{\pm}(u)$ and $\z^{\pm}(u)$ with coefficients
in the respective algebra.
\epr

\bpf
The proof is the same in both cases so we only consider
the algebra $U(\overline{R})$.
Multiply both sides of \eqref{gen rel1} by $u/v-\xi$ and set
$u/v=\xi$ to get
\beql{QL1L2}
Q\tss\Lc^{\pm}_1(u\tss\xi)\Lc^{\pm}_2(u)=\Lc^{\pm}_2(u)\Lc^{\pm}_1(u\tss\xi)\tss Q.
\eeq
By the definition of the element $Q$, we can write
$Q=D_1^{-1}P^{\ts\tra_1}D_1$.
Therefore, \eqref{QL1L2} takes the form
\beql{pdld}
P^{\ts\tra_1}D_1\Lc_1^{\pm}(u\tss\xi)D_1^{-1}\Lc_2^{\pm}(u)
=\Lc_2^{\pm}(u)D_1\Lc_1^{\pm}(u\tss\xi)D_1^{-1}P^{\ts\tra_1}.
\eeq
The image of the operator $P^{\ts\tra_1}$ in $\End(\CC^{2n})^{\ot 2}$
is one-dimensional,
so that each side of this equality
must be equal to $P^{\ts\tra_1}$ times a certain series $\z^{\pm}(u)$
with coefficients in $U(\overline{R})$. Observe that
$P^{\ts\tra_1}D_1=P^{\ts\tra_1}D^{-1}_2$
and $P^{\ts\tra_1}\Lc_1^{\pm}(u\tss\xi)=P^{\ts\tra_1}\Lc_2^{\pm}(u\tss\xi)^{\tra}$
and so we get
\ben
P^{\ts\tra_1}D_2\Lc_2^{\pm}(u\tss\xi)^{\tra}D_2^{-1}\Lc_2^{\pm}(u)
=\Lc_2^{\pm}(u)D_2\Lc_2^{\pm}(u\tss\xi)^{\tra}D_2^{-1}P^{\ts\tra_1}=\z^{\pm}(u)\tss P^{\ts\tra_1}.
\een
The required relations now follow by taking trace of the first copy of $\End\CC^{2n}$.
\epf

\bpr\label{prop:centr}
All coefficients of the series $z^{+}(u)$ and $z^{-}(u)$
belong to the center of the algebra $U(R)$.
\epr

\bpf
We will verify that $z^{+}(u)$ commutes with all series $l_{ij}^-(v)$; the remaining
cases follow by similar or simpler arguments.
By the defining relations \eqref{rllmp} we can write
\ben
\bal
D_1L^{+}_1(u\tss\xi)^{\tra}D_1^{-1}L^{+}_1(u)L^{-}_2(v)
=D_1L^{+}_1(u\tss\xi)^{\tra}D_1^{-1}{R}(uq^c/v)^{-1}L^{-}_2(v)L^{+}_1(u){R}(uq^{-c}/v).
\eal
\een
By \eqref{crsymr} the right hand side equals
\ben
\xi^{-2}q^2D_1L^{+}_1(u\tss\xi)^{\tra}{R}(u\tss\xi q^c/v)^{\tra_1}L^{-}_2(v)
D_1^{-1}L^{+}_1(u){R}(uq^{-c}/v).
\een
Applying the patrial transposition $\tra_1$ to both sides in
\eqref{rllmp} we get the relation
\ben
L^{+}_1(u\tss\xi)^{\tra}{R}(u\tss\xi q^c/v)^{\tra_1}L^{-}_2(v)
=L^{-}_2(v){R}(u\tss\xi q^{-c}/v)^{\tra_1}L^{+}_1(u\tss\xi)^{\tra}.
\een
Hence, using \eqref{crsymr} and \eqref{DLDL} we obtain
\ben
\bal
z^{+}(u)L^{-}_2(v)&=D_1L^{+}_1(u\tss\xi)^{\tra}D_1^{-1}L^{+}_1(u)L^{-}_2(v)\\[0.4em]
&=\xi^{-2}q^2 L^{-}_2(v)D_1 {R}(u\tss\xi q^{-c}/v)^{\tra_1}
D_1^{-1}D_1L^{+}_1(u\tss\xi)^{\tra}D_1^{-1}L^{+}_1(u){R}(uq^{-c}/v)\\[0.4em]
&=\xi^{-2}q^2 L^{-}_2(v)D_1 {R}(u\tss\xi q^{-c}/v)^{\tra_1}
D_1^{-1}z^{+}(u){R}(uq^{-c}/v)=L^{-}_2(v)z^{+}(u),
\eal
\een
as required.
\epf

\bre\label{rem:noncent}
The crossing symmetry properties \eqref{crsymr} of the $R$-matrix $R(u)$ were essential
for Proposition~\ref{prop:centr} to hold.
Although the coefficients of the series $\z^{+}(u)$ and $\z^{-}(u)$ are central
in the respective subalgebras of $U(\overline{R})$ generated by the coefficients of the
series $\ell^+_{ij}(u)$ and $\ell^-_{ij}(u)$, they are not central in the entire
algebra $U(\overline{R})$.
\qed
\ere

\subsection{Quasideterminants and quantum minors}

Let $A=[a_{ij}]$ be an $N\times N$ matrix over a ring with $1$.
Denote by $A^{ij}$ the matrix obtained from $A$
by deleting the $i$-th row
and $j$-th column. Suppose that the matrix
$A^{ij}$ is invertible.
The $ij$-{\em th quasideterminant of} $A$
is defined by the formula
\ben
|A|_{ij}=a_{ij}-r^{\tss j}_i(A^{ij})^{-1}\ts c^{\tss i}_j,
\een
where $r^{\tss j}_i$ is the row matrix obtained from the $i$-th
row of $A$ by deleting the element $a_{ij}$, and $c^{\tss i}_j$
is the column matrix obtained from the $j$-th
column of $A$ by deleting the element $a_{ij}$; see
\cite{gr:dm}, \cite{gr:tn}.
The quasideterminant $|A|_{ij}$ is also denoted
by boxing the entry $a_{ij}$ in the matrix $A$.

Throughout the rest of this section we will regard elements of the
tensor product algebra $\End(\CC^{2n})^{\ot m}\ot U(\overline{R})$
as operators on the space $(\CC^{2n})^{\ot m}$ with coefficients in $U(\overline{R})$.
Accordingly, for such an element
\ben
X=\sum_{a_i,b_i}\ts e_{a_1b_1}\ot\dots\ot e_{a_mb_m}
\ot X^{a_1\dots a_m}_{\ts b_1\dots b_m}
\een
we will use a standard notation
\beql{matrelem}
X^{a_1\dots a_m}_{\ts b_1\dots b_m}=\langle a_1,\dots, a_m\ts |\ts X\ts |\ts b_1,\dots, b_m\ts\rangle
\eeq
and its counterparts $X\ts |\ts b_1,\dots, b_m\ts\rangle$ and $\langle a_1,\dots, a_m\ts |\ts X$.

Consider the algebra $U(\overline{R})$ and for any $2\leqslant i,j\leqslant 2\pr$
introduce the quasideterminant
\ben
s^{\pm}_{ij}(u)=\left|\begin{matrix}
\ell^{\tss\pm}_{11}(u)&\ell^{\tss\pm}_{1j}(u)\\[0.2em]
\ell^{\tss\pm}_{i1}(u)&\boxed{\ell^{\tss\pm}_{ij}(u)}
\end{matrix}\right|
=\ell^{\tss\pm}_{ij}(u)-\ell^{\tss\pm}_{i1}(u)\tss
\ell^{\tss\pm}_{11}(u)^{-1}\tss \ell^{\tss\pm}_{1j}(u).
\een
Let the power series
$\ell^{\tss\pm\tss a_1 a_2}_{~~\tss b_1 b_2}(u)$ ({\em quantum minors})
with coefficients in ${U}(\overline{R})$ be defined by
\beql{quamintau}
\ell^{\tss\pm\tss a_1 a_2}_{~~\tss b_1 b_2}(u)=
\langle a_1,a_2\ts |\ts \wh{R}(q^{-2})\ts \Lc^{\pm}_1(u)\ts \Lc^{\pm}_2(u\tss q^{2})\ts |
\ts b_1,b_2\rangle,
\eeq
where $a_i,b_i\in\{1,\dots,2n\}$ and we set
\beql{rhat}
\wh{R}(u)=\frac{uq-q^{-1}}{u-1}\ts\overline{R}(u).
\eeq
The following symmetry properties are straightforward to verify.

\ble\label{lem:skewsymm}\quad
(i)\quad
If $a_1\ne a'_2$ and $a_1<a_2$ then $\ell^{\tss\pm\tss a_1 a_2}_{~~\tss b_1 b_2}(u)
=-q^{-1}\ts\ell^{\tss\pm\tss a_2 a_1}_{~~\tss b_1 b_2}(u)$.
\medskip

\noindent
(ii)\quad
If $b_1\ne b'_2$ and $b_1<b_2$ then $\ell^{\tss\pm\tss a_1 a_2}_{~~\tss b_1 b_2}(u)
=-q\ts \ell^{\tss\pm\tss a_1 a_2}_{~~\tss b_2 b_1}(u)$.
\qed
\ele

\ble\label{lem:toneone}
For any $2\leqslant i,j\leqslant 2\pr$ we have
\beql{quasi-minor}
s^{\pm}_{ij}(u)=\ell^{\tss\pm}_{11}(uq^{-2})^{-1}\tss \ell^{\tss\pm\tss 1 i}_{~~\tss 1 j}(uq^{-2}).
\eeq
Moreover,
\beql{commtoo}
\big[\ell^{\tss\pm}_{11}(u), \ell^{\tss\pm\tss 1 i}_{~~\tss 1 j}(v)\big]=0
\eeq
and
\beql{pm}
\frac{q^{-1}u_{\pm}-q\tss v_{\mp}}{u_{\pm}-v_{\mp}}\ts
\ell_{11}^{\tss\pm}(u)\ts\ell^{\tss\mp\tss 1i}_{~~1j}(v)
=\frac{q^{-1}u_{\mp}-q\tss v_{\pm}}{u_{\mp}-v_{\pm}}\ts
\ell^{\tss\mp\tss 1i}_{~~1j}(v)\ts\ell_{11}^{\tss\pm}(u).
\eeq
\ele

\bpf
By the definition of quantum minors,
\beql{l1i1j}
\ell^{\tss\pm\tss 1 i}_{~~\tss 1 j}(u)
=\langle 1,i\ts|\ts\wh{R}(q^{-2})\Lc^{\pm}_1(u)\Lc^{\pm}_2(uq^2)\ts|\ts 1,j\rangle
=\ell_{11}^{\tss\pm}(u)\ell_{ij}^{\tss\pm}(uq^2)-q^{-1}\ell_{i1}^{\tss\pm}(u)\ell_{1j}^{\tss\pm}(uq^2).
\eeq
The defining relations of the algebra $U(\overline{R})$ give
\ben
\langle 1,i\ts|\ts\overline{R}(u/v)\Lc^{\pm}_1(u)
\Lc^{\pm}_2(v)\ts|\ts1,1\rangle
=\langle 1,i\ts|\ts\Lc^{\pm}_2(v)\Lc^{\pm}_1(u)\overline{R}(u/v)\ts|\ts1,1\rangle,
\een
and so
\ben
\Big(\frac{u}{v}-1\Big)\ell^{\tss\pm}_{11}(u)\ell^{\tss\pm}_{i1}(v)
+(q-q^{-1})\frac{u}{v}\ts\ell^{\tss\pm}_{i1}(u)\ell^{\tss\pm}_{11}(v)
=\Big(\frac{u}{v}q-q^{-1}\Big)\ell_{i1}^{\tss\pm}(v)\ell^{\tss\pm}_{11}(u).
\een
In particular,
\beql{l11li1}
\ell_{11}^{\tss\pm}(uq^{-2})\ell_{i1}^{\tss\pm}(u)=q^{-1}\ell_{i1}^{\tss\pm}(uq^{-2})\ell^{\tss\pm}_{11}(u).
\eeq
Relations \eqref{l1i1j} and \eqref{l11li1} imply
\ben
\bal
\ell^{\tss\pm\tss 1 i}_{~~\tss 1 j}(uq^{-2})&
=\ell_{11}^{\tss\pm}(uq^{-2})\ell_{ij}^{\tss\pm}(u)-q^{-1}\ell_{i1}^{\tss\pm}(uq^{-2})\ell_{1j}^{\tss\pm}(u)\\
&=\ell_{11}^{\tss\pm}(uq^{-2})\ell_{ij}^{\tss\pm}(u)-q^{-1}\ell_{i1}^{\tss\pm}(uq^{-2})
\ell^{\tss\pm}_{11}(u)\ell^{\tss\pm}_{11}(u)^{-1}\ell_{1j}^{\tss\pm}(u)\\
&=\ell_{11}^{\tss\pm}(uq^{-2})\ell_{ij}^{\tss\pm}(u)-\ell_{11}^{\tss\pm}(uq^{-2})
\ell^{\tss\pm}_{i1}(u)\ell^{\tss\pm}_{11}(u)^{-1}\ell_{1j}^{\tss\pm}(u)
=\ell_{11}^{\tss\pm}(uq^{-2})s^{\pm}_{ij}(u),
\eal
\een
thus proving \eqref{quasi-minor}. To verify
\eqref{pm}, note that by the Yang--Baxter equation
\eqref{ybe} and relations \eqref{gen rel2} and \eqref{gen rel3} we have
\ben
\begin{aligned}
&\langle 1,1,i\ts|\ts\overline{R}_{01}({u_{\pm}}/{v_{\mp}})
\overline{R}_{02}({u_{\pm}\ts q^{-2}}/{v_{\mp}})
\Lc_0^{\pm}(u)\wh{R}_{12}(q^{-2})\Lc_1^{\mp}(v)\Lc_2^{\mp}(vq^2)\ts|\ts1,1,j\rangle=\\[0.4em]
&\langle1,1,i\ts|\ts\wh{R}_{12}(q^{-2})\Lc_1^{\mp}(v)\Lc_2^{\mp}(vq^2)
\Lc_0^{\pm}(u)\overline{R}_{02}({u_{\mp}\ts q^{-2}}/{v_{\pm}})
\overline{R}_{01}({u_{\pm}}/{v_{\mp}})\ts|\ts1,1,j\rangle
\end{aligned}
\een
which gives \eqref{pm}. The calculation for
\eqref{commtoo} is quite similar.
\epf

We point out the following consequences of Lemma \ref{lem:toneone}:
for $2\leqslant i,j\leqslant 2^{\tss\prime}$ we have
\beql{LS}
[\ell^{\tss\pm}_{11}(u),s^{\pm}_{ij}(v)]=0
\eeq
and
\beql{LPMSMP}
\frac{u_{\pm}-v_{\mp}}{qu_{\pm}-q^{-1}v_{\mp}}\ell_{11}^{\tss\pm}(u)s^{\mp}_{ij}(v)
=\frac{u_{\mp}-v_{\pm}}{qu_{\mp}-q^{-1}v_{\pm}}s^{\mp}_{ij}(v)\ts\ell_{11}^{\tss\pm}(u).
\eeq

\subsection{Homomorphism theorems}

Now we aim to make a connection between the algebras $U(\overline{R})$
associated with the Lie algebras $\spa_{2n-2}$ and $\spa_{2n}$.
Since the rank $n$ will vary, we will indicate the dependence on $n$ by
adding a subscript $[n]$ to the $R$-matrices.
Consider the algebra $U(\overline{R}^{\tss[n-1]})$ and let
the indices
of the generators $\ell^{\tss\pm}_{ij}[\mp m]$ range over the sets
$2\leqslant i,j\leqslant 2\pr$ and $m=0,1,\dots$, where
$i\pr=2n-i+1$, as before.

\bth\label{thm:embed}
The mappings $q^{\pm c/2}\mapsto q^{\pm c/2}$ and
\beql{embedgen}
\ell^{\tss\pm}_{ij}(u)\mapsto \left|\begin{matrix}
\ell^{\tss\pm}_{11}(u)&\ell^{\tss\pm}_{1j}(u)\\[0.2em]
\ell^{\tss\pm}_{i1}(u)&\boxed{\ell^{\tss\pm}_{ij}(u)}
\end{matrix}\right|,\qquad 2\leqslant i,j\leqslant 2\pr,
\eeq
define a homomorphism ${U}(\overline{R}^{\tss[n-1]})\to {U}(\overline{R}^{\tss[n]})$.
\eth

\bpf
Consider the tensor product algebra $\End(\CC^{2n})^{\ot 4}\ot U(\overline{R}^{\tss[n]})$.
We begin with calculations of certain matrix elements of operators which are
straightforward from the definition of the $R$-matrix \eqref{rbar}.
We will use notation \eqref{rhat} and
suppose that $2\leqslant i,j\leqslant 2^{\tss\prime}$. Then
\ben
\bal
\overline{R}^{\tss[n]}_{13}(a)\overline{R}^{\tss[n]}_{23}(aq^2)\ts|\ts1,i,1,j\rangle
&=\frac{aq-q^{-1}}{aq^2-q^{-2}}\ts|\ts1,i,1,j\rangle
+\frac{(q-q^{-1})(a-1)aq^2}{(aq^2-1)(aq^2-q^{-2})}\ts|\ts1,1,i,j\rangle\\
&+\frac{(q-q^{-1})^2aq^2}{(aq^2-1)(aq^2-q^{-2})}\ts|\ts i,1,1,j\rangle
\eal
\een
and
\beql{R12R13R231i1j}
\wh{R}^{\tss[n]}_{12}(q^{-2})\overline{R}^{\tss[n]}_{13}(a)
\overline{R}^{\tss[n]}_{23}(aq^2)\ts|\ts1,i,1,j\rangle
=\frac{a-1}{aq-q^{-1}}\wh{R}^{\tss[n]}_{12}(q^{-2})\ts|\ts1,i,1,j\rangle.
\eeq
Furthermore, we have
\begin{align}\label{R34R141i1j}
\wh{R}^{\tss[n]}_{34}(q^{-2})\overline{R}^{\tss[n]}_{14}(aq^{-2})\ts|\ts1,i,1,j\rangle
&=\frac{aq^{-1}-q}{a-1}\wh{R}^{\tss[n]}_{34}(q^{-2})\ts|\ts1,i,1,j\rangle,\\
\label{1N11}
\wh{R}^{\tss[n]}_{12}(q^{-2})\wh{R}^{\tss[n]}_{34}(q^{-2})
\overline{R}^{\tss[n]}_{14}(aq^{-2})\ts|\ts1,1',1,1\rangle&=0
\end{align}
and
\begin{align}\non
\wh{R}^{\tss[n]}_{12}(q^{-2})\wh{R}^{\tss[n]}_{34}(q^{-2})
\overline{R}^{\tss[n]}_{14}(aq^{-2})\ts|\ts1,1,1,1'\rangle&=
-\frac{(q^{-2}-1)(aq^{-1}-q)}{(a-1)(a\xi^{-1}q^{-2}-1)}
\wh{R}^{\tss[n]}_{12}(q^{-2})\wh{R}^{\tss[n]}_{34}(q^{-2})\\
&\qquad\times
\sum_{a=2}^{2^{\tss\prime}}\varepsilon_{1'}\varepsilon_aq^{\bar{a}-\bar{1'}}\ts|\ts1,a',1,a\rangle.
\label{111N}
\end{align}
These observations together with the formula for
matrix elements of $\overline{R}^{\tss[n]}_{24}(a)$ given by
\ben
\begin{aligned}
\overline{R}^{\tss[n]}_{24}(a)&\ts|\ts1,i,1,j\rangle
=\frac{a-1}{aq-q^{-1}}\big\{\big(\delta_{ij}q^{1-\delta_{ij'}}
+(1-\delta_{ij})(1-\delta_{ij'})
+(1-\delta_{i,i'})q^{-1}\delta_{ij'}\big)\ts|\ts1,i,1,j\rangle\\
&+\delta_{i>j}(q-q^{-1})\ts|\ts1,j,1,i\rangle-(q-q^{-1})\delta_{ij'}
\sum_{a>j}\varepsilon_a\varepsilon_j q^{\bar{a}-\bj}\ts|\ts1,a',1,a\rangle\\
&+\frac{q-q^{-1}}{a-1}\ts|\ts1,j,1,i\rangle+\frac{q^{-1}-q}{a\xi^{-1}-1}\delta_{ij'}
\sum_{a=2}^{2^{\tss\prime}}\varepsilon_a\varepsilon_j q^{\bar{a}-\bj}\ts|\ts1,a',1,a\rangle\\
&+\frac{q^{-1}-q}{a\xi^{-1}-1}\delta_{ij'}
\varepsilon_1\varepsilon_j q^{\bar{1}-\bj}\ts|\ts1,1',1,1\rangle
+\frac{(q^{-1}-q)a\xi^{-1}}{a\xi^{-1}-1}\delta_{ij'}
\varepsilon_{1'}\varepsilon_j q^{\bar{1'}-\bj}\ts|\ts1,1,1,1'\rangle\big\}
\end{aligned}
\een
lead to the relation
\begin{multline}\label{R12R34R14R241i1j}
\wh{R}^{\tss[n]}_{12}(q^{-2})\wh{R}^{\tss[n]}_{34}(q^{-2})
\overline{R}^{\tss[n]}_{14}(aq^{-2})\overline{R}^{\tss[n]}_{24}(a)\ts|\ts1,i,1,j\rangle\\
=\frac{aq^{-1}-q}{a-1}\wh{R}^{\tss[n]}_{12}(q^{-2})
\wh{R}^{\tss[n]}_{34}(q^{-2})\overline{R}^{\tss[n-1]}_{24}(a)\ts|\ts1,i,1,j\rangle.
\end{multline}
Now applying the Yang--Baxter equation \eqref{ybe} and relations
\eqref{R12R13R231i1j} and \eqref{R12R34R14R241i1j}
we deduce the following matrix element formulas:
\begin{multline}\label{R23ran}
\wh{R}^{\tss[n]}_{12}(q^{-2})\wh{R}^{\tss[n]}_{34}(q^{-2})\overline{R}^{\tss[n]}_{14}(aq^{-2})
\overline{R}^{\tss[n]}_{24}(a)\overline{R}^{\tss[n]}_{13}(a)
\overline{R}^{\tss[n]}_{23}(aq^2)\ts|\ts1,i,1,j\rangle\\[0.4em]
{}=\frac{aq^{-1}-q}{aq-q^{-1}}\ts\wh{R}^{\tss[n]}_{12}(q^{-2})
\wh{R}^{\tss[n]}_{34}(q^{-2})\overline{R}^{\tss[n-1]}_{24}(a)\ts|\ts1,i,1,j\rangle
\end{multline}
and
\begin{multline}\label{lanR23}
\langle 1,i,1,j\ts|\ts\overline{R}^{\tss[n]}_{23}(aq^2)\overline{R}^{\tss[n]}_{13}(a)
\overline{R}^{\tss[n]}_{24}(a)\overline{R}^{\tss[n]}_{14}(aq^{-2})
\wh{R}^{\tss[n]}_{12}(q^{-2})\wh{R}^{\tss[n]}_{34}(q^{-2})\\[0.4em]
{}=\frac{aq^{-1}-q}{aq-q^{-1}}\ts\langle 1,i,1,j\ts|\ts\overline{R}^{\tss[n-1]}_{24}(a)
\wh{R}^{\tss[n]}_{12}(q^{-2})\wh{R}^{\tss[n]}_{34}(q^{-2}).
\end{multline}

To complete the proof of the theorem, introduce the matrices
\ben
\Gamma^{\pm}(u)=\sum\limits_{i,j=2}^{2^{\tss\prime}} e_{ij}\otimes
\ell^{\tss\pm\tss 1 i}_{~~ 1 j}(u)\in \End \CC^{2n}\ot U(\overline{R}^{\tss[n]}).
\een
Our next step is to verify that the following relations hold in
the algebra $U(\overline{R}^{\tss[n]})$:
\begin{align}\non
\overline{R}^{\tss[n-1]}(u/v)\Gamma^{\pm}_1(u)\Gamma^{\pm}_2(v)
&=\Gamma^{\pm}_2(v)\Gamma^{\pm}_1(u)\overline{R}^{\tss[n-1]}(u/v),\\[0.4em]
\frac{q^{-1}{u_+}-qv_{-}}{q{u_+}-q^{-1}v_{-}}\ts
\overline{R}^{\tss[n-1]}(uq^c/v)\Gamma^{+}_1(u)\Gamma^{-}_2(v)&
=\frac{q^{-1}{u_-}-qv_{+}}{q{u_-}-q^{-1}v_{+}}
\ts\Gamma^{-}_2(v)\Gamma^{+}_1(u)\overline{R}^{\tss[n-1]}(uq^{-c}/v).
\non
\end{align}
The calculations are quite similar in both cases so we only give details for the first relation.
The Yang--Baxter equation and the defining relations for the algebra $U(\overline{R}^{\tss[n]})$ give
\begin{multline}\non
\overline{R}^{\tss[n]}_{23}\big(\frac{u q^2}{v}\big)\overline{R}^{\tss[n]}_{13}\big(\frac uv\big)
\overline{R}^{\tss[n]}_{24}\big(\frac{u}{v}\big)\overline{R}^{\tss[n]}_{14}\big(\frac{u}{v q^{2}}\big)
\wh{R}^{\tss[n]}_{12}(q^{-2})\Lc^{\pm}_1(u)\ts \Lc^{\pm}_2(uq^{2})
\wh{R}^{\tss[n]}_{34}(q^{-2})\Lc^{\pm}_3(v)\ts \Lc^{\pm}_4(vq^{2})\\[0.4em]
=\wh{R}^{\tss[n]}_{34}(q^{-2})\Lc^{\pm}_3(v)\ts \Lc^{\pm}_4(vq^{2})
\wh{R}^{\tss[n]}_{12}(q^{-2})\Lc^{\pm}_1(u)\ts \Lc^{\pm}_2(uq^{2})
\overline{R}^{\tss[n]}_{14}\big(\frac{u}{v q^{2}}\big)\overline{R}^{\tss[n]}_{24}\big(\frac{u}{v}\big)
\overline{R}^{\tss[n]}_{13}\big(\frac uv\big)\overline{R}^{\tss[n]}_{23}\big(\frac{u q^2}{v}\big).
\end{multline}
Hence, assuming that $2\leqslant i,j,k,l\leqslant 2^{\tss\prime}$
and applying \eqref{R23ran} and \eqref{lanR23}
we get
\begin{multline}\non
\langle 1,k,1,l\ts|\ts\overline{R}^{\tss[n-1]}_{24}\big(\frac{u}{v}\big)
\wh{R}^{\tss[n]}_{12}(q^{-2})\Lc^{\pm}_1(u)\ts \Lc^{\pm}_2\big(uq^{2}\big)
\wh{R}^{\tss[n]}_{34}(q^{-2})\Lc^{\pm}_3(v)\ts \Lc^{\pm}_4(vq^{2})\ts|\ts1,i,1,j\rangle\\[0.4em]
=\langle 1,k,1,l\ts|\ts\wh{R}^{\tss[n]}_{34}(q^{-2})\Lc^{\pm}_3(v)\ts \Lc^{\pm}_4(vq^{2})
\wh{R}^{\tss[n]}_{12}(q^{-2})\Lc^{\pm}_1(u)\ts \Lc^{\pm}_2(uq^{2})
\overline{R}^{\tss[n]}_{24}\big(\frac{u}{v}\big)\ts|\ts1,i,1,j\rangle,
\end{multline}
which is equivalent to
\ben
\overline{R}^{\tss[n-1]}_{24}(u/v)\Gamma^{\pm}_2(u)\Gamma^{\pm}_4(v)
=\Gamma^{\pm}_4(v)\Gamma^{\pm}_2(u)\overline{R}^{\tss[n-1]}_{24}(u/v),
\een
as required.
Finally, set
\ben
S^{\pm}(u)=\sum\limits_{2\leqslant i,j\leqslant 2\pr}e_{ij}\ot s^{\pm}_{ij}(u).
\een
By Lemma \ref{lem:toneone},
\ben
S^{\pm}(u)=\ell^{\tss\pm}_{11}(uq^{-2})^{-1}\Ga^{\pm}(uq^{-2})
\een
and
\ben
\frac{q^{-1}u_{\pm}-qv_{\mp}}{u_{\pm}-v_{\mp}}\ts\ell_{11}^{\tss\pm}(u)\ts\Ga^{\mp}(v)
=\frac{q^{-1}u_{\mp}-qv_{\tss\pm}}{u_{\mp}-v_{\pm}}\ts\Ga^{\mp}(v)\ts\ell_{11}^{\tss\pm}(u).
\een
The above relations for the matrices $\Ga^{\pm}(u)$ imply
\ben
\bal
\overline{R}^{\tss[n-1]}(u/v)S^{\pm}_1(u)S^{\pm}_2(v)
&=S^{\pm}_2(v)S^{\pm}_1(u)\overline{R}^{\tss[n-1]}(u/v),\\[0.4em]
\overline{R}^{\tss[n-1]}(uq^{\pm c}/v)S^{\pm}_1(u)S^{\mp}_2(v)
&=S^{\mp}_2(v)S^{\pm}_1(u)\overline{R}^{\tss[n-1]}(uq^{\mp c}/v),
\eal
\een
thus completing the proof.
\epf

The following is a generalization of
Theorem~\ref{thm:embed} which is immediate from
the Sylvester theorem for quasideterminants
\cite{gr:dm}, \cite{kl:mi}; cf. the proof of its Yangian counterpart
given in \cite[Thm~3.7]{jlm:ib}.
Fix a positive integer $m$ such that
$m< n$.
Suppose that the generators $\ell_{ij}^{\tss\pm}(u)$ of the algebra $U(\overline{R}^{\tss[n-m]})$ are
labelled by the indices
$m+1\leqslant i,j\leqslant (m+1)\pr$ with $i\pr=2\tss n-i+1$ as before.

\bth\label{thm:red}
The mapping
\beql{redu}
\ell^{\tss\pm}_{ij}(u)\mapsto \left|\begin{matrix}
\ell^{\tss\pm}_{11}(u)&\dots&\ell^{\tss\pm}_{1m}(u)&\ell^{\tss\pm}_{1j}(u)\\
\dots&\dots&\dots&\dots\\
\ell^{\tss\pm}_{m1}(u)&\dots&\ell^{\tss\pm}_{mm}(u)&\ell^{\tss\pm}_{mj}(u)\\[0.2em]
\ell^{\tss\pm}_{i1}(u)&\dots&\ell^{\tss\pm}_{im}(u)&\boxed{\ell^{\tss\pm}_{ij}(u)}
\end{matrix}\right|,\qquad m+1\leqslant i,j\leqslant (m+1)\pr,
\eeq
defines a homomorphism $\psi_m:U(\overline{R}^{\tss[n-m]})\to U(\overline{R}^{\tss[n]})$.
\qed
\eth

As another application of the Sylvester theorem for quasideterminants, we get
a consistence property of the homomorphisms \eqref{redu};
cf.~\cite[Prop.~3.8]{jlm:ib} and \cite[eq.~(1.85)]{m:yc} for its Yangian counterparts.
We will write $\psi_m=\psi^{(n)}_m$ to indicate the dependence of $n$.
For a parameter $l$ we have the corresponding homomorphism
\ben
\psi^{(n-l)}_m:U(\overline{R}^{\tss[n-l-m]})\to U(\overline{R}^{\tss[n-l]})
\een
provided by \eqref{redu}. Then we have the equality of maps
\beql{consist}
\psi^{(n)}_l\circ\psi^{(n-l)}_m=\psi^{(n)}_{l+m}.
\eeq

Suppose that $\{a_1,\dots,a_k\}$ and $\{b_1,\dots,b_k\}$ are subsets of $\{1,\dots,2n\}$,
assuming that $a_1<a_2<\dots<a_k$ and $b_1<b_2<\dots<b_k$, such that
$a_i\neq a_j'$ and $b_i\neq b_j'$ for all $i,j$. Introduce the corresponding
{\em type $A$ quantum minors} as the matrix elements \eqref{matrelem}:
\begin{multline}\non
\ell^{\tss\pm~a_1,\dots,a_k}_{\ts\ts~~b_1,\dots,b_k}(u)
=\langle a_1,\dots,a_k\ts|\ts \wh{R}_{k-1,k}(\wh{R}_{k-2,k}\wh{R}_{k-2,k-1})
\dots(\wh{R}_{1,k}\dots\wh{R}_{1,2})\\[0.3em]
{}\times\Lc_1^{\pm}(u)\Lc_2^{\pm}(uq^2)\dots\Lc_k^{\pm}(uq^{2k-2})\ts|\ts b_1,\dots,b_k\rangle,
\end{multline}
where $\wh{R}_{ij}=\wh{R}_{ij}(q^{2(i-j)})$.  They are
given by the following formulas:
\ben
\bal
\ell^{\tss\pm~a_1,\dots,a_k}_{\ts\ts~~b_1,\dots,b_k}(u)
&=\sum_{\sigma\in \Sym_k}(-q)^{-l(\sigma)}\ell^{\tss\pm}_{a_{\sigma(1)}b_1}(u)
\dots\ell^{\tss\pm}_{a_{\sigma(k)}b_k}(uq^{2k-2})\\
&=\sum_{\sigma\in \Sym_k}(-q)^{l(\sigma)}
\ell^{\tss\pm}_{a_{k}b_{\sigma(k)}}(uq^{2k-2})\dots\ell^{\tss\pm}_{a_{1}b_{\sigma(1)}}(u),
\eal
\een
where $l(\sigma)$ denotes the number of inversions of the permutation $\sigma\in \Sym_k$.
The assumptions on the indices $a_i$ and $b_i$ imply that
certain relations for these
quantum minors take the same form
as those for the quantum affine algebra $U_q(\wh\gl_n)$. Such relations for the latter can be
deduced by applying $R$-matrix calculations which are quite analogous to the Yangian case;
cf.~\cite{hm:qa}, \cite{kl:mi} and \cite[Ch.~1]{m:yc}.
In particular, for $1\leqslant i,j\leqslant k$ we have
\ben
\bal[]
[\ell^{\tss\pm}_{a_ib_j}(u),\ell^{\tss\pm~a_1,a_2,\dots,a_k}_{\ts\ts~~b_1,b_2,\dots,b_k}(v)]&=0,\\[0.4em]
 \displaystyle\prod_{a=1}^{k-1}\frac{u_{\pm}q^{-k}-v_{\mp}q^{k}}{u_{\pm}q^{1-k}-v_{\mp}q^{k-1}}
\ell^{\tss\pm}_{a_ib_j}(u)\ell^{\tss\mp~a_1,a_2,\dots,a_k}_{\ts\ts~~b_1,b_2,\dots,b_k}(v)&=
\displaystyle\prod_{a=1}^{k-1}\frac{u_{\mp}q^{-k}-v_{\pm}q^{k}}{u_{\mp}q^{1-k}-v_{\pm}q^{k-1}}
\ell^{\tss\mp~a_1,a_2,\dots,a_k}_{\ts\ts~~b_1,b_2,\dots,b_k}(v)\ell^{\tss\pm}_{a_ib_j}(u).
\eal
\een

\bco\label{cor:commu}
Under the assumptions of Theorem~\ref{thm:red} we have
\ben
\bal
\big[\ell^{\tss\pm}_{ab}(u),\psi_m(\ell^{\tss\pm}_{ij}(v))\big]&=0,\\[0.4em]
\frac{u_{\pm}-v_{\mp}}{qu_{\pm}-q^{-1}v_{\mp}}\ell^{\tss\pm}_{ab}(u)\psi_m(\ell^{\mp}_{ij}(v))
{}&=\frac{u_{\mp}-v_{\pm}}{qu_{\mp}-q^{-1}v_{\pm}}\psi_m(\ell^{\mp}_{ij}(v))\ell^{\tss\pm}_{ab}(u),
\eal
\een
for all $1\leqslant a,b\leqslant m$ and $m+1\leqslant i,j\leqslant (m+1)\pr$.
\eco

\bpf
Both formulas are verified with the use of the relations between the
quasideterminants and quantum minors:
\ben
\left|\begin{matrix}
\ell^{\tss\pm}_{11}(u)&\dots&\ell^{\tss\pm}_{1m}(u)&\ell^{\tss\pm}_{1j}(u)\\
\dots&\dots&\dots&\dots\\
\ell^{\tss\pm}_{m1}(u)&\dots&\ell^{\tss\pm}_{mm}(u)&\ell^{\tss\pm}_{mj}(u)\\[0.2em]
\ell^{\tss\pm}_{i1}(u)&\dots&\ell^{\tss\pm}_{im}(u)&\boxed{\ell^{\tss\pm}_{ij}(u)}
\end{matrix}\right|={\ell\ts}^{\pm\tss1\dots\ts m}_{~~1\dots\ts m}(uq^{-2m})^{-1}\cdot
{\ell\ts}^{\pm\tss1\dots\ts m\ts i}_{~~1\dots\ts m\ts j}(uq^{-2m}).
\een
They are consequences of the type $A$ relations; cf. the Yangian case in
\cite[Sec.~3]{jlm:ib}.
\epf

\section{Gauss decomposition}\label{sec:gauss decom}

We will apply the Gauss decompositions \eqref{gaussdec} to the generators matrices
$L^{\pm}(u)$ and $\Lc^{\pm}(u)$ for the respective algebras $U(R^{[n]})$
and $U(\overline{R}^{\tss[n]})$. Each of these algebras
is generated by the coefficients of the matrix elements of the triangular and diagonal matrices
which we will refer to as the {\em Gaussian generators}. Our goal in this section
is to produce necessary relations satisfied by these generators to be able to get presentations
of the $R$-matrix algebras $U(R^{[n]})$ and $U(\overline{R}^{\tss[n]})$.

\subsection{Gaussian generators}

The entries of the matrices $F^{\pm}(u)$, $H^{\pm}(u)$ and $E^{\pm}(u)$ which occur
in the decompositions \eqref{gaussdec} can be described by the universal quasideterminant
formulas
as follows \cite{gr:dm}, \cite{gr:tn}:
\beql{hmqua}
h^{\pm}_i(u)=\begin{vmatrix} l^{\pm}_{1\tss 1}(u)&\dots&l^{\pm}_{1\ts i-1}(u)&l^{\pm}_{1\tss i}(u)\\
                          \vdots&\ddots&\vdots&\vdots\\
                         l^{\pm}_{i-1\ts 1}(u)&\dots&l^{\pm}_{i-1\ts i-1}(u)&l^{\pm}_{i-1\ts i}(u)\\[0.2em]
                         l^{\pm}_{i\tss 1}(u)&\dots&l^{\pm}_{i\ts i-1}(u)&\boxed{l^{\pm}_{i\tss i}(u)}\\
           \end{vmatrix},\qquad i=1,\dots,2n,
\eeq
whereas
\beql{eijmlqua}
e^{\pm}_{ij}(u)=h^{\pm}_i(u)^{-1}\ts\begin{vmatrix} l^{\pm}_{1\tss 1}(u)&\dots&
l^{\pm}_{1\ts i-1}(u)&l^{\pm}_{1\ts j}(u)\\
                          \vdots&\ddots&\vdots&\vdots\\
                         l^{\pm}_{i-1\ts 1}(u)&\dots&l^{\pm}_{i-1\ts i-1}(u)&l^{\pm}_{i-1\ts j}(u)\\[0.2em]
                         l^{\pm}_{i\tss 1}(u)&\dots&l^{\pm}_{i\ts i-1}(u)&\boxed{l^{\pm}_{i\tss j}(u)}\\
           \end{vmatrix}
\eeq
and
\beql{fijlmqua}
f^{\pm}_{ji}(u)=\begin{vmatrix} l^{\pm}_{1\tss 1}(u)&\dots&l^{\pm}_{1\ts i-1}(u)&l^{\pm}_{1\tss i}(u)\\
                          \vdots&\ddots&\vdots&\vdots\\
                         l^{\pm}_{i-1\ts 1}(u)&\dots&l^{\pm}_{i-1\ts i-1}(u)&l^{\pm}_{i-1\ts i}(u)\\[0.2em]
                         l^{\pm}_{j\ts 1}(u)&\dots&l^{\pm}_{j\ts i-1}(u)&\boxed{l^{\pm}_{j\tss i}(u)}\\
           \end{vmatrix}\ts h^{\pm}_i(u)^{-1}
\eeq
for $1\leqslant i<j\leqslant 2n$. The same formulas hold for the expressions
of the entries of the respective triangular matrices $\Fc^{\pm}(u)$ and
$\Ec^{\pm}(u)$ and the diagonal matrices $\Hc^{\pm}(u)=\diag\ts[\h^{\pm}_{i}(u)]$
in terms of the formal series $\ell^{\tss\pm}_{ij}(u)$, which arise from the Gauss decomposition
\ben
\Lc^{\pm}(u)=\Fc^{\pm}(u)\ts\Hc^{\pm}(u)\ts\Ec^{\pm}(u)
\een
for the algebra $U(\overline{R}^{\tss[n]})$. We will denote by
$\e_{ij}(u)$ and $\f_{ji}(u)$ the entries of the respective matrices
$\Ec^{\pm}(u)$ and $\Fc^{\pm}(u)$ for $i<j$.

The following Laurent series with coefficients in
the respective algebras $U(R^{[n]})$ and $U(\overline{R}^{\tss[n]})$
will be used frequently:
\begin{align}
\label{Xi}
X^{+}_i(u)&=e^{+}_{ii+1}(u_{+})-e_{ii+1}^{-}(u_{-}),
\qquad X^{-}_i(u)=f^{+}_{i+1,i}(u_{-})-f^{-}_{i+1,i}(u_{+}),\\[0.3em]
\Xc^{+}_i(u)&=\e^{+}_{ii+1}(u_{+})-\e_{ii+1}^{-}(u_{-}),
\qquad \Xc^{-}_i(u)=\f^{+}_{i+1,i}(u_{-})-\f^{-}_{i+1,i}(u_{+}).
\label{Xbar}
\end{align}

\bpr\label{prop:corrgauss}
Under the homomorphism $U(\overline{R})\to \Hc_q(n)\ot_{\CC[q^c,\ts q^{-c}]}U(R)$
provided by Proposition~\ref{prop:homheis} we have
\ben
\bal
  \e^{\pm}_{ij}(u) & \mapsto e^{\pm}_{ij}(u), \\
  \f^{\pm}_{ij}(u) & \mapsto f^{\pm}_{ij}(u), \\
  \h^{\pm}_{i}(u) &\mapsto \exp\sum\limits _{k=1}^{\infty}\be_{\mp k}u^{\mp k}\cdot h^{\pm}_{i}(u).
\eal
\een
\epr

\bpf
This is immediate from the formulas for the Gaussian generators.
\epf

\subsection{Images of the generators under the homomorphism $\psi_m$}
Suppose that $0\leqslant m< n$.
We will use the superscript $[n-m]$ to indicate
square submatrices corresponding to rows and columns labelled by
$m+1,m+2,\dots,(m+1)'$. In particular, we set
\ben
\Fc^{\pm[n-m]}(u)=\begin{bmatrix}
1&0&\dots&0\ts\\
\f^{\pm}_{m+2\ts m+1}(u)&1&\dots&0\\
\vdots&\ddots&\ddots&\vdots\\
\f^{\pm}_{(m+1)'\tss m+1}(u)&\dots&\f^{\pm}_{(m+1)'\ts (m+2)'}(u)&1
\end{bmatrix},
\een
\ben
\Ec^{\pm[n-m]}(u)=\begin{bmatrix} 1&\e^{\pm}_{m+1\tss m+2}(u)&\ldots&\e^{\pm}_{m+1\tss (m+1)'}(u)\\
                        0&1&\ddots &\vdots\\
                         \vdots&\vdots&\ddots&\e^{\pm}_{(m+2)'\tss(m+1)'}(u)\\
                         0&0&\ldots&1\\
           \end{bmatrix}
\een
and $\Hc^{\pm[n-m]}(u)=\diag\ts\big[\h^{\pm}_{m+1}(u),\dots,\h^{\pm}_{(m+1)'}(u)\big]$. Furthermore,
introduce the products of these matrices by
\ben
\Lc^{\pm[n-m]}(u)=\Fc^{\pm[n-m]}(u)\ts \Hc^{\pm[n-m]}(u)\ts \Ec^{\pm[n-m]}(u).
\een
The entries of $\Lc^{\pm[n-m]}(u)$ will be denoted by $\ell^{\tss\pm[n-m]}_{ij}(u)$.

\bpr\label{prop:gauss-consist}
The series $\ell^{\tss\pm[n-m]}_{ij}(u)$ coincides with the image of
the generator series $\ell^{\tss\pm}_{ij}(u)$
of the extended quantum affine algebra $U(\overline{R}^{\tss[n-m]})$
under the homomorphism \eqref{redu},
\ben
\ell^{\pm[n-m]}_{ij}(u)=\psi_m\big(\ell^{\pm}_{ij}(u)\big),\qquad m+1\leqslant i,j\leqslant (m+1)'.
\een
\epr

\bpf
This follows by the same argument as for the Yangian case; see \cite[Prop.~4.1]{jlm:ib}.
\epf

\bco\label{cor:guass-embed}
The following relations hold in $U(\overline{R}^{\tss[n]})${\rm :}
\beql{subRTT}
\overline{R}^{\tss[n-m]}_{12}(u/v)\ts \Lc^{\pm[n-m]}_1(u)\ts \Lc^{\pm[n-m]}_2(v)
=\Lc^{\pm[n-m]}_2(v)\ts \Lc^{\pm[n-m]}_1(u)\ts \overline{R}^{\tss[n-m]}_{12}(u/v),
\eeq
\beql{subRTT2}
\overline{R}^{\tss[n-m]}_{12}(u_{+}/v_{-})\ts \Lc^{+[n-m]}_1(u)\ts \Lc^{-[n-m]}_2(v)
=\Lc^{-[n-m]}_2(v)\ts \Lc^{+[n-m]}_1(u)\ts \overline{R}^{\tss[n-m]}_{12}(u_{-}/v_{+}).
\eeq
\eco

\bpf
This is immediate from Proposition~\ref{prop:gauss-consist}.
\epf

\bpr\label{prop:relmone}
Suppose that
$m+1\leqslant j,k,l\leqslant (m+1)'$ and $j\neq l'$. Then the following relations hold
in $U(\overline{R}^{\tss[n]})${\rm :} if $j=l$ then
\begin{align}\label{ELMPj=l}
\e_{mj}^{\pm}(u)\ell^{\mp [n-m]}_{kl}(v)
&=\frac{qu_{\mp}-q^{-1}v_{\pm}}{u_{\mp}-v_{\pm}}\ell^{\mp [n-m]}_{kj}(v)\e_{ml}^{\pm}(u)
-\frac{(q-q^{-1})u_{\mp}}{u_{\mp}-v_{\pm}}\ell^{\mp [n-m]}_{kj}(v)\e_{mj}^{\mp}(v),\\
\e_{mj}^{\pm}(u)\ell^{\tss\pm [n-m]}_{kl}(v)
&=\frac{qu-q^{-1}v}{u-v}\ell^{\tss\pm [n-m]}_{kj}(v)\e_{ml}^{\pm}(u)
-\frac{(q-q^{-1})u}{u-v}\ell^{\tss\pm [n-m]}_{kj}(v)\e_{mj}^{\pm}(v);
\non
\end{align}
if $j<l$ then
\begin{align}\label{ELMPj<l}
[\e^{\pm}_{mj}(u),\ell^{\mp [n-m]}_{kl}(v)]
&=\frac{(q-q^{-1})v_{\pm}}{u_{\mp}-v_{\pm}}\ell^{\mp [n-m]}_{kj}(v)\e^{\pm}_{ml}(u)-
\frac{(q-q^{-1})u_{\mp}}{u_{\mp}-v_{\pm}}\ell^{\mp [n-m]}_{kj}(v)\e^{\mp}_{ml}(v),\\
[\e^{\pm}_{mj}(u),\ell^{\tss\pm [n-m]}_{kl}(v)]
&=\frac{(q-q^{-1})v}{u-v}\ell^{\tss\pm [n-m]}_{kj}(v)\e^{\pm}_{ml}(u)-
\frac{(q-q^{-1})u}{u-v}\ell^{\tss\pm [n-m]}_{kj}(v)\e^{\pm}_{ml}(v);
\non
\end{align}
if $j>l$ then
\ben
\bal[]
[\e^{\pm}_{mj}(u),\ell^{\mp [n-m]}_{kl}(v)]
&=\frac{(q-q^{-1})u_{\mp}}{u_{\mp}-v_{\pm}}\ell^{\mp [n-m]}_{kj}(v)(\e^{\pm}_{ml}(u)-\e_{ml}^{\mp}(v)),\\
[\e^{\pm}_{mj}(u),\ell^{\tss\pm [n-m]}_{kl}(v)]
&=\frac{(q-q^{-1})u}{u-v}\ell^{\tss\pm [n-m]}_{kj}(v)(\e^{\pm}_{ml}(u)-\e_{ml}^{\pm}(v)).
\eal
\een
\epr

\bpf
It is sufficient to verify the relations for $m=1$;
the general case will then follow by the application of the homomorphism $\psi_m$.
The calculations are similar for all the relations so we only
verify \eqref{ELMPj<l}.
By the defining relations,
\begin{multline}\label{l1jlkl}
\frac{u_{\pm}-v_{\mp}}{qu_{\pm}-q^{-1}v_{\mp}}\ell_{1j}^{\pm}(u)\ell_{kl}^{\mp}(v)+
\frac{(q-q^{-1})u_{\pm}}{qu_{\pm}-q^{-1}v_{\mp}}\ell_{kj}^{\pm}(u)\ell_{1l}^{\mp}(v)\\
=\frac{u_{\mp}-v_{\pm}}{qu_{\mp}-q^{-1}v_{\pm}}\ell_{kl}^{\mp}(v)\ell_{1j}^{\pm}(u)
+\frac{(q-q^{-1})v_{\pm}}{qu_{\mp}-q^{-1}v_{\pm}}\ell_{kj}^{\mp}(v)\ell_{1l}^{\pm}(u).
\end{multline}
Since $\ell_{kl}^{\mp}(v)=\ell^{\mp [n-1]}_{kl}(v)+\f_{k1}^{\mp}(v)\ts\h_1^{\mp}(v)\ts\e_{1l}^{\mp}(v)$,
we can write the left hand side of \eqref{l1jlkl} as
\begin{multline}\non
\frac{u_{\pm}-v_{\mp}}{qu_{\pm}-q^{-1}v_{\mp}}\ell_{1j}^{\pm}(u)\ell_{kl}^{\mp [n-1]}(v)
+\frac{u_{\pm}-v_{\mp}}{qu_{\pm}-q^{-1}v_{\mp}}
\ell_{1j}^{\pm}(u)\f^{\mp}_{k1}(v)\ts\h_1^{\mp}(v)\ts\e_{1l}^{\mp}(v)\\
{}+\frac{(q-q^{-1})v_{\mp}}{qu_{\pm}-q^{-1}v_{\mp}}\ell_{kj}^{\pm}(u)\ell_{1l}^{\mp}(v).
\end{multline}
By the defining relations, we have
\begin{multline}\non
\frac{u_{\pm}-v_{\mp}}{qu_{\pm}-q^{-1}v_{\mp}}\ell_{1j}^{\pm}(u)\ell^{\mp}_{k1}(v)+
\frac{(q-q^{-1})u_{\pm}}{qu_{\pm}-q^{-1}v_{\mp}}\ell_{kj}^{\pm}(u)\ell_{11}^{\mp}(v)\\
{}=\frac{u_{\mp}-v_{\pm}}{qu_{\mp}-q^{-1}v_{\pm}}\ell^{\mp}_{k1}(v)\ell_{1j}^{\pm}(u)+
\frac{(q-q^{-1})u_{\mp}}{qu_{\mp}-q^{-1}v_{\pm}}\ell_{kj}^{\mp}(v)\ell_{11}^{\pm}(u).
\end{multline}
Hence, the left hand side of \eqref{l1jlkl} equals
\begin{multline}\non
\frac{u_{\pm}-v_{\mp}}{qu_{\pm}-q^{-1}v_{\mp}}\ell_{1j}^{\pm}(u)\ell_{kl}^{\mp [n-1]}(v)
+\frac{u_{\mp}-v_{\pm}}{qu_{\mp}-q^{-1}v_{\pm}}\f^{\mp}_{k1}(v)\ell^{\mp}_{11}(v)
\ell_{1j}^{\pm}(u)\e_{1l}^{\mp}(v)\\
+\frac{(q-q^{-1})u_{\mp}}{qu_{\mp}-q^{-1}v_{\pm}}\ell_{kj}^{\mp}(v)\ell_{11}^{\pm}(u)\e_{1l}^{\mp}(v).
\end{multline}
Furthermore, using the relation
\ben
\ell^{\tss\pm}_{1j}(u)\ell^{\mp}_{11}(v)=\frac{u_{\mp}-v_{\pm}}{qu_{\mp}-q^{-1}v_{\pm}}\ell^{\mp}_{11}(v)
\ell^{\tss\pm}_{1j}(u)+\frac{(q-q^{-1})u_{\mp}}{qu_{\mp}-q^{-1}v_{\pm}}\ell^{\mp}_{1j}(v)
\ell^{\tss\pm}_{11}(u),
\een
we can bring the left hand side of \eqref{l1jlkl} to the form
\ben
\frac{u_{\pm}-v_{\mp}}{qu_{\pm}-q^{-1}v_{\mp}}\ell_{1j}^{\pm}(u)
\ell_{kl}^{\mp [n-1]}(v)+\f^{\mp}_{k1}(v)\ell_{1j}^{\pm}(u)\ell_{1l}^{\mp}(v)
+\frac{(q-q^{-1})u_{\mp}}{qu_{\mp}-q^{-1}v_{\pm}}
\ell_{kj}^{\mp [n-1]}(v)\ell_{11}^{\pm}(u)\e_{1l}^{\mp}(v).
\een
For $j<l$ we have
\ben
\ell_{1j}^{\pm}(u)\ell_{1l}^{\mp}(v)
=\frac{u_{\mp}-v_{\pm}}{qu_{\mp}-q^{-1}v_{\pm}}\ell_{1l}^{\mp}(v)\ell_{1j}^{\pm}(u)
+\frac{(q-q^{-1})v_{\pm}}{qu_{\mp}-q^{-1}v_{\pm}}\ell_{1j}^{\mp}(v)\ell_{1l}^{\pm}(u),
\een
so that the left hand side of \eqref{l1jlkl} becomes
\begin{multline}\non
\frac{u_{\pm}-v_{\mp}}{qu_{\pm}-q^{-1}v_{\mp}}\ell_{1j}^{\pm}(u)\ell_{kl}^{\mp [n-1]}(v)
-\frac{u_{\mp}-v_{\pm}}{qu_{\mp}-q^{-1}vu_{\pm}}
\ell_{kl}^{\mp [n-1]}(v)\ell_{1j}^{\pm}(u)\\
{}=\frac{(q-q^{-1})v_{\pm}}{qu_{\mp}-q^{-1}v_{\pm}}
\ell_{kj}^{\mp [n-1]}(v)\ell_{11}^{\pm}(u)\e_{1l}^{\pm}(u)-
\frac{(q-q^{-1})u_{\mp}}{qu_{\mp}-q^{-1}v_{\pm}}
\ell_{kj}^{\mp [n-1]}(v)\ell_{11}^{\pm}(u)\e_{1l}^{\mp}(v).
\end{multline}
Finally, Corollary~\ref{cor:commu} implies
\ben
\bal
\frac{u_{\pm}-v_{\mp}}{qu_{\pm}-q^{-1}v_{\mp}}\ell_{11}^{\pm}(u)\ell^{\mp [n-1]}_{kl}(v)
{}&=\frac{u_{\mp}-v_{\pm}}{qu_{\mp}-q^{-1}v_{\pm}}\ell^{\mp [n-1]}_{kl}(v)\ell_{11}^{\pm}(u),\\
\frac{u_{\pm}-v_{\mp}}{qu_{\pm}-q^{-1}v_{\mp}}\ell_{11}^{\pm}(u)\ell^{\mp [n-1]}_{kj}(v)
{}&=\frac{u_{\mp}-v_{\pm}}{qu_{\mp}-q^{-1}v_{\pm}}\ell^{\mp [n-1]}_{kj}(v)\ell_{11}^{\pm}(u),
\eal
\een
thus completing the proof of \eqref{ELMPj<l}.
\epf

Quite similar arguments prove the following counterpart of Proposition~\ref{prop:relmone}
involving the generator series $\f^{\pm}_{ji}(u)$.

\bpr\label{prop:relmonf}
Suppose that
$m+1\leqslant j,k,l\leqslant (m+1)'$ and $j\neq k'$. Then the following relations hold
in $U(\overline{R}^{\tss[n]})${\rm :} if $j=k$ then
\ben
\bal
\f_{jm}^{\pm}(u)\ell_{jl}^{\mp [n-m]}(v)
&=\frac{u_{\pm}-v_{\mp}}{qu_{\pm}-q^{-1}v_{\mp}}\ell_{jl}^{\mp [n-m]}(v)\f_{jm}^{\pm}(u)
+\frac{(q-q^{-1})v_{\mp}}{qu_{\pm}-q^{-1}v_{\mp}}\ts\f_{jm}^{\mp}(v)\ell_{jl}^{\mp [n-m]}(v),\\
\f_{jm}^{\pm}(u)\ell_{jl}^{\pm [n-m]}(v)
&=\frac{u-v}{qu-q^{-1}v}\ell_{jl}^{\pm [n-m]}(uv)\ts\f_{jm}^{\pm}(u)
+\frac{(q-q^{-1})v}{qu-q^{-1}v}\f_{jm}^{\pm}(v)\ell_{jl}^{\pm [n-m]}(v);
\eal
\een
if $j<k$ then
\ben
\bal[]
[\f_{jm}^{\pm}(u),\ell_{kl}^{\mp [n-m]}(v)]
&=\frac{(q-q^{-1})v_{\mp}}{u_{\pm}-v_{\mp}}\f_{km}^{\mp}(v)\ell_{jl}^{\mp [n-m]}(v)
-\frac{(q-q^{-1})u_{\pm}}{u_{\pm}-v_{\mp}}\f_{km}^{\pm}(u)\ell_{jl}^{\mp [n-m]}(v),\\
[\f_{jm}^{\pm}(u),\ell_{kl}^{\pm [n-m]}(v)]
&=\frac{(q-q^{-1})v}{u-v}\f_{km}^{\pm}(v)\ell_{jl}^{\pm [n-m]}(v)
-\frac{(q-q^{-1})u}{u-v}\f_{km}^{\pm}(u)\ell_{jl}^{\pm [n-m]}(v);
\eal
\een
if $j>k$ then
\ben
\bal[]
[\f_{jm}^{\pm}(u),\ell_{kl}^{\mp [n-m]}(v)]
&=\frac{(q-q^{-1})v_{\mp}}{u_{\pm}-v_{\mp}}(\f_{km}^{\mp}(v)-\f_{km}^{\pm}(u))
\ell_{jl}^{\mp [n-m]}(v),\\
[\f_{jm}^{\pm}(u),\ell_{kl}^{\pm [n-m]}(v)]
&=\frac{(q-q^{-1})v}{u-v}(\f_{km}^{\pm}(v)-\f_{km}^{\pm}(u))\ell_{jl}^{\pm [n-m]}(v).
\eal
\vspace{-0.7cm}
\een
\qed
\epr

\subsection{Type $A$ relations}
\label{subsec:typea}

Due to the observation made in Remark~\ref{rem:gln}
and the quasideterminant formulas \eqref{hmqua}, \eqref{eijmlqua}
and \eqref{fijlmqua}, some of the relations between the gaussian generators will follow
from those for the quantum affine algebra
$U_q(\wh{\mathfrak{gl}}_n)$; see \cite{df:it}. To reproduce them, set
\ben
\Lc^{A\pm}(u)=\sum_{i,j=1}^n e_{ij}\otimes \ell^{\tss\pm}_{ij}(u)
\een
and consider the $R$-matrix used in \cite{df:it} which is given by
\begin{multline}
R_A(u)=\sum_{i=1}^n e_{ii}\otimes e_{ii}+\frac{u-1}{qu-q^{-1}}
\sum_{i\neq j} e_{ii}\otimes e_{jj}\\
+\frac{q-q^{-1}}{qu-q^{-1}}\sum_{i>j}e_{ij}\otimes e_{ji}
+\frac{(q-q^{-1})u}{qu-q^{-1}}\sum_{i<j}e_{ij}\otimes e_{ji}.
\label{rtypea}
\end{multline}
By comparing it with the $R$-matrix \eqref{rbar}, we come to
the relations in the algebra
$U(\overline{R}^{\tss[n]})${\rm:}
\ben
\bal
R_A(u/w)\Lc^{A\pm}_{1}(u)\Lc^{A\pm}_2(v)&=\Lc^{A\pm}_2(v)\Lc^{A\pm}_{1}(u){R}_A(u/v),\\[0.3em]
R_A(uq^c/v)\Lc^{A+}_1(u)\Lc^{A-}_2(v)&=\Lc^{A-}_2(v)\Lc^{A+}_1(u){R}_A(uq^{-c}/v).
\eal
\een
Hence we get the following relations for the Gaussian generators
$\h^{\pm}_i(u)$ with $i=1,\dots,n$ and $\Xc_i^{\pm}(u)$ with $i=1,\dots,n-1$,
which were verified in
\cite{df:it}; see \eqref{Xbar}.

\bpr\label{TypeArelation}
In
the algebra $U(\overline{R}^{\tss[n]})$ we have
\ben
\h^{\pm}_i(u)\h^{\pm}_j(v)=\h^{\pm}_j(v)\h^{\pm}_i(u),\qquad
\h^{\pm}_i(u)\h^{\mp}_i(v)=\h^{\mp}_i(v)\h^{\pm}_i(u),
\een
\ben
\frac{u_{\pm}-v_{\mp}}{qu_{\pm}-q^{-1}v_{\mp}}\h^{\pm}_i(u)\h^{\mp}_j(v)=
\frac{u_{\mp}-v_{\pm}}{qu_{\mp}-q^{-1}v_{\pm}}\h^{\mp}_j(v)\h^{\pm}_i(u) \qquad\text{for}\quad i<j.
\een
Moreover,
\ben
\bal
\h_{i}^{\pm}(u)\Xc_{j}^{+}(v)
&=\frac{u_{\mp}-v}{q^{(\ep_i,\alpha_j)}u_{\mp}-q^{-(\ep_i,\alpha_j)}v} \Xc_{j}^{+}(v)\h_{i}^{\pm}(u),\\
\h_{i}^{\pm}(u)\Xc_{j}^{-}(v)
&=\frac{q^{(\ep_i,\alpha_j)}u_{\pm}-q^{-(\ep_i,\alpha_j)}v}{u_{\pm}-v} \Xc_{j}^{-}(v)\h_{i}^{\pm}(u),
\eal
\een
and
\ben
(u-q^{\pm (\alpha_i,\alpha_j)}v)\Xc_{i}^{\pm}(uq^i)\Xc_{j}^{\pm}(vq^j)
=(q^{\pm (\alpha_i,\alpha_j)}u-v) \Xc_{j}^{\pm}(vq^j)\Xc_{i}^{\pm}(uq^i),
\een
\ben
[\Xc_i^{+}(u),\Xc_j^{-}(v)]=\delta_{ij}(q-q^{-1})
\Big(\delta\big(u\ts q^{-c}/v\big)\ts\h_i^{-}(v_+)^{-1}\h_{i+1}^{-}(v_+)
-\delta\big(u\ts q^{c}/v\big)\ts\h_i^{+}(u_+)^{-1}\h_{i+1}^{+}(u_+)\Big)
\een
together with the Serre relations \eqref{serre} for the series $\Xc_i^{\pm}(u)$.
\qed
\epr

\bre\label{re:anotherA}
Consider the inverse matrices $\Lc^{\pm}(u)^{-1}=[\ell^{\tss\pm}_{ij}(u)']_{i,j=1,\dots,2n}$.
By the defining relations
\eqref{gen rel1} and \eqref{gen rel2}, we have
\ben
\bal
\Lc^{\pm}_{1}(u)^{-1}\Lc^{\pm}_2(v)^{-1}\overline{R}^{\tss[n]}(u/v)
&=\overline{R}^{\tss[n]}(u/v)\Lc^{\pm}_2(v)^{-1}\Lc^{\pm}_{1}(u)^{-1},\\
\Lc^{-}_2(v)^{-1}\Lc^{+}_1(u)^{-1}\overline{R}^{\tss[n]}({u}q^{c}/v)
&=\overline{R}^{\tss[n]}({u}q^{-c}/v)\Lc^{-}_2(v)^{-1}\Lc^{+}_1(u)^{-1}.
\eal
\een
So we can get another family of generators of the algebra $U(\overline{R}^{\tss[n]})$
which satisfy the defining relations of $U_q(\wh{\gl}_n)$. Namely, these relations
are satisfied by
the coefficients of the series
$\ell^{\tss\pm}_{ij}(u)'$ with $i,j=n',\dots,1'$. In particular,
by taking the inverse matrices, we get a Gauss decomposition
for the matrix $[\ell^{\tss\pm}_{ij}(u)']_{i,j=n',\dots,1'}$ from
the Gauss decomposition of the matrix $\Lc^{\pm}(u)$.
\qed
\ere

\subsection{Relations for the long root generators}
\label{subsec:rlr}

In the particular case $n=1$ the $R$-matrix \eqref{rbar} takes the form
\begin{multline}
\overline{R}^{\tss[1]}(u)=\sum_{i=n}^{n+1} e_{ii}\ot e_{ii}+\frac{u-1}{q^2u-q^{-2}}
\sum_{i\neq j}e_{ii}\ot e_{jj}\\
{}+\frac{(q^2-q^{-2})u}{q^2u-q^{-2}}e_{n,n+1}\ot e_{n+1,n}
+\frac{q^2-q^{-2}}{q^2u-q^{-2}} e_{n+1,n}\ot e_{n,n+1}
\non
\end{multline}
and so it coincides with
the $R$-matrix associated with $U_{q^2}(\wh{\gl}_2)$; cf.~\eqref{rtypea}.
Therefore, a set of relations involving the long root generators are implied by
Corollary~\ref{cor:guass-embed} and Proposition~\ref{TypeArelation}.

\bpr\label{LowC}
The following relations hold
in the algebra $U(\overline{R}^{\tss[n]})$:
\ben
\bal
\h^{\pm}_i(u)\h^{\pm}_j(v)&=\h^{\pm}_j(v)\h^{\pm}_i(u),\qquad i,j=n,n+1,\\[0.3em]
\h^{\pm}_i(u)\h^{\mp}_i(v)&=\h^{\mp}_i(v)\h^{\pm}_i(u),\qquad i=n,n+1,
\eal
\een
\ben
\frac{u_{\pm}-v_{\mp}}{q^2u_{\pm}-q^{-2}v_{\mp}}\h^{\pm}_n(u)\h^{\mp}_{n+1}(v)
=\frac{u_{\mp}-v_{\pm}}{q^2u_{\mp}-q^{-2}v_{\pm}}\h^{\mp}_{n+1}(v)\h^{\pm}_n(u).
\een
Moreover,
\ben
\bal
\h_{n}^{\pm}(u)\Xc_{n}^{+}(v)&=
\frac{u_{\mp}-v}{q^{2}u_{\mp}-q^{-2}v} \Xc_{n}^{+}(v)\h_{n}^{\pm}(u),\\
\h_{n+1}^{\pm}(u)\Xc_{n}^{+}(v)&=
\frac{u_{\mp}-v}{q^{-2}u_{\mp}-q^{2}v} \Xc_{n}^{+}(v)\h_{n+1}^{\pm}(u),\\
\h_{n}^{\pm}(u)\Xc_{n}^{-}(v)&=
\frac{q^{2}u_{\pm}-q^{-2}v}{u_{\pm}-v} \Xc_{n}^{-}(v)\h_{n}^{\pm}(u),\\
\h_{n+1}^{\pm}(u)\Xc_{n}^{-}(v)&=
\frac{q^{-2}u_{\pm}-q^{2}v}{u_{\pm}-v} \Xc_{n}^{-}(v)\h_{n+1}^{\pm}(u),
\eal
\een
and
\ben
\bal
&(u-q^{\pm (\alpha_n,\alpha_n)}v)\Xc_{n}^{\pm}(u)\Xc_{n}^{\pm}(v)
=(q^{\pm (\alpha_n,\alpha_n)}u-v) \Xc_{n}^{\pm}(v)\Xc_{n}^{\pm}(u),\\[0.4em]
[\Xc_n^{+}(u),\Xc_n^{-}(v)]&=(q^2-q^{-2})
\Big(\delta\big(\frac{u}{v\tss q^{c}}\big)\ts\h_n^{-}(v_+)^{-1}\h_{n+1}^{-}(v_+)
-\delta\big(\frac{u\tss q^{c}}{v}\big)\ts\h_n^{+}(u_+)^{-1}\h_{n+1}^{+}(u_+)\Big).
\eal
\een
\epr

\subsection{Formulas for the series $z^{\pm}(u)$ and $\z^{\pm}(u)$}
\label{subsec:fsz}

Recall that the series $z^{\pm}(u)$ and $\z^{\pm}(u)$
were defined in Proposition~\ref{prop:central}. We will now indicate the dependence on $n$
by adding the corresponding superscript.
Write relation \eqref{DLbarDLbar} in the form
\beql{DLbarD}
D\Lc^{\pm}(u\xi)^{\tra}D^{-1}=\Lc^{\pm}(u)^{-1}\z^{\pm\ts[n]}(u).
\eeq
Using the Gauss decomposition for $\Lc^{\pm}(u)$ and taking the $(2n,2n)$-entry
on both sides of \eqref{DLbarD} we get
\beql{h1h1pr}
\h_1^{\pm}(u\xi)=\h_{1'}^{\pm}(u)^{-1}\z^{\pm\ts[n]}(u).
\eeq
By Proposition \ref{prop:gauss-consist}, in the subalgebra generated by the
coefficients of the series $\ell^{\pm[1]}_{ij}(u)$
we have
\ben
\z^{\pm\ts[1]}(u)=\h_n^{\pm}(uq^{-4})\ts\h_{n'}^{\pm}(u).
\een

\ble\label{lem:eiprei}
The following relations hold in the algebra $U(\overline{R}^{\tss[n]})$ for
$1\leqslant i\leqslant n-1${\rm:}
\beql{ei'ei}
\e^{\pm}_{(i+1)'\ts i'}(u)=-\e_{i,i+1}^{\pm}(u\xi q^{2i})
\Fand
\f^{\pm}_{i'\ts (i+1)'}(u)=-\f_{i+1,i}^{\pm}(u\xi q^{2i}).
\eeq
\ele

\bpf
By Propositions~\ref{prop:central} and \ref{prop:gauss-consist}, for any $1\leqslant i\leqslant n-1$
we have
\beql{Lbar[n-i+1] z}
\Lc^{\pm\ts[n-i+1]}(u)^{-1}\z^{\pm\ts[n-i+1]}(u)=D^{[n-i+1]}
\Lc^{\pm\ts[n-i+1]}(u\xi q^{2i-2})'(D^{[n-i+1]})^{-1},
\eeq
where
\ben
D^{[n-i+1]}=\diag [\ts q^{n-i+1},\dots,q,q^{-1},\dots,q^{-n+i-1}\tss].
\een
By taking the $(i',i')$ and $((i+1)',i')$ entries on both sides of \eqref{Lbar[n-i+1] z} we get
\beql{hihipr}
\h_i^{\pm}(u\xi q^{2i-2})=\h_{i'}^{\pm}(u)^{-1}\z^{\pm\ts[n-i+1]}(u)
\eeq
and
\ben
-\e^{\pm}_{(i+1)',i'}(u)\ts\h^{\pm}_{i'}(u)^{-1}\ts\z^{\pm\ts[n-i+1]}(u)
=q\ts \h_i^{\pm}(u\xi q^{2i-2})\e_{i,i+1}^{\pm}(u\xi q^{2i-2}).
\een
Due to \eqref{hihipr}, this formula can written as
\beql{ei'hi}
-\e^{\pm}_{(i+1)',i'}(u)\ts\h_i^{\pm}(u\xi q^{2i-2})
=q\ts \h_i^{\pm}(u\xi q^{2i-2})\ts\e_{i,i+1}^{\pm}(u\xi q^{2i-2}).
\eeq
Furthermore, by the results of \cite{df:it},
\ben
q\ts\h_i^{\pm}(u)\ts\e_{i,i+1}^{\pm}(u)=\e_{i,i+1}^{\pm}(uq^2)\ts\h_i^{\pm}(u),
\een
so that \eqref{ei'hi} is equivalent to
\beq
-\e^{\pm}_{(i+1)',i'}(u)\ts\h_i^{\pm}(u\xi q^{2i-2})
=\e_{i,i+1}^{\pm}(u\xi q^{2i})\ts\h_i^{\pm}(u\xi q^{2i-2}),
\eeq
thus proving the first relation in \eqref{ei'ei}. The second relation is verified
in a similar way.
\epf

\bpr\label{prop:formze}
In the algebras $U(\overline{R}^{\tss[n]})$ and $U({R}^{[n]})$ we have
the respective formulas:
\ben
\bal
\z^{\pm\ts[n]}(u)&=\prod_{i=1}^{n-1}\h^{\pm}_{i}(u\tss\xi q^{2i})^{-1}
\prod_{i=1}^{n}\h^{\pm}_{i}(u\tss\xi q^{2i-2})\ts\h^{\pm}_{n+1}(u),\\
{z}^{\pm\ts[n]}(u)&=\prod_{i=1}^{n-1}{h}^{\pm}_{i}(u\xi q^{2i})^{-1}
\prod_{i=1}^{n}{h}^{\pm}_{i}(u\xi q^{2i-2})\ts{h}^{\pm}_{n+1}(u).
\eal
\een
\epr

\bpf
The arguments for both formulas are quite similar so we only give the proof of the first one.
Taking the $(2^{\tss\prime},2^{\tss\prime})$-entry on both sides of
\eqref{Lbar[n-i+1] z} and expressing the
entries of the matrices $\Lc^{\pm\ts[n]}(u)^{-1}$ and $\Lc^{\pm\ts[n]}(u\xi)^{\tra}$
in terms of the Gauss generators,
we get
\ben
\h^{\pm}_2(u\xi)+\f^{\pm}_{21}(u\xi)\ts\h^{\pm}_{1}(u\xi)\ts\e^{\pm}_{12}(u\xi)
=\big(\h^{\pm}_{2^{\tss\prime}}(u)^{-1}
+\e^{\pm}_{2^{\tss\prime},1'}(u)\h^{\pm}_{1'}(u)^{-1}
\f^{\pm}_{1',2^{\tss\prime}}(u)^{-1}\big)\z^{\pm\ts[n]}(u).
\een
As we pointed out in Remark~\ref{rem:noncent},
the coefficients of the series $\z^{\pm\ts[n]}(u)$ are central in the
respective subalgebras generated by the coefficients of $\ell^{\pm\ts[n]}_{ij}(u)$.
Therefore, using
\eqref{h1h1pr}, we can rewrite the above relation as
\ben
\h^{\pm}_{2^{\tss\prime}}(u)^{-1}\z^{\pm\ts[n]}(u)=
\h^{\pm}_2(u\xi)+\f^{\pm}_{21}(u\xi)\ts\h^{\pm}_{1}(u\xi)\ts\e^{\pm}_{12}(u\xi)
-\e^{\pm}_{2^{\tss\prime},1'}(u)\ts\h^{\pm}_{1}(u\xi)\ts\f^{\pm}_{1',2^{\tss\prime}}(u).
\een
Now apply Lemma~\ref{lem:eiprei} to obtain
\ben
\h^{\pm}_{2^{\tss\prime}}(u)^{-1}\z^{\pm\ts[n]}(u)=
\h^{\pm}_2(u\xi)+\f^{\pm}_{21}(u\xi)\ts\h^{\pm}_{1}(u\xi)\ts\e^{\pm}_{12}(u\xi)
-\e^{\pm}_{12}(u\xi q^2)\ts\h^{\pm}_{1}(u\xi)\ts\f^{\pm}_{21}(u\xi q^2).
\een
On the other hand, by the results of \cite{df:it} we have
\ben
\h^{\pm}_{1}(u)\ts\e^{\pm}_{12}(u)=q^{-1}\e^{\pm}_{12}(uq^2)\ts\h^{\pm}_{1}(u),
\qquad
\h^{\pm}_{1}(u)\ts\f^{\pm}_{21}(uq^2)=q^{-1}\f^{\pm}_{21}(u)\ts\h^{\pm}_{1}(u),
\een
and
\ben
[\e^{\pm}_{12}(u),\f^{\pm}_{21}(v)]=
\frac{u(q-q^{-1})}{u-v}\ts(\h^{\pm}_{2}(v)\h^{\pm}_{1}(v)^{-1}
-\h^{\pm}_{2}(u)\ts\h^{\pm}_{1}(u)^{-1}).
\een
This leads to the expression
\ben
\h^{\pm}_{2^{\tss\prime}}(u)^{-1}\z^{\pm\ts[n]}(u)=
\h^{\pm}_2(u\xi q^2)\ts\h^{\pm}_1(u\xi q^2)^{-1}\h^{\pm}_{1}(u\xi).
\een
Since $\z^{\pm\ts[n-1]}(u)=\h^{\pm}_{2^{\tss\prime}}(u)\h^{\pm}_2(u\xi q^2)$,
we get a recurrence formula
\ben
\z^{\pm\ts[n]}(u)=
\h^{\pm}_1(u\xi q^2)^{-1}\h^{\pm}_{1}(u\xi)\ts\z^{\pm\ts[n-1]}(u)
\een
thus completing the proof.
\epf

\subsection{Drinfeld-type relations in the algebras $U(\overline{R}^{\tss[n]})$ and $U(R^{[n]}_{})$}
\label{subsec:dtr}

We will now extend the sets of relations produced in
Secs~\ref{subsec:typea} and \ref{subsec:rlr} to obtain all necessary relations
in the algebras $U(\overline{R}^{\tss[n]})$ and $U(R^{[n]}_{})$ to be able to prove
the Main Theorem.

\bth\label{thm:relrbar}
The following relations between the gaussian generators hold in the algebra
$U(\overline{R}^{\tss[n]})$. For the relations involving $\h^{\pm}_i(u)$ we have
\beq
\h_{i,0}^{+}\h_{i,0}^{-}=\h_{i,0}^{-}\h_{i,0}^{+}=1, \tss \text{and}~~\tss \h^{+}_{n,0}\h^{+}_{n+1,0}=1,
\eeq
\begin{align}\non
\h^{\pm}_i(u)\h^{\pm}_j(v)&=\h^{\pm}_j(v)\h^{\pm}_i(u),\\[0.4em]
\h^{\pm}_i(u)\h^{\mp}_i(v)&=\h^{\mp}_i(v)\h^{\pm}_i(u),
\non\\
\frac{u_{\pm}-v_{\mp}}{qu_{\pm}-q^{-1}v_{\mp}}\h^{\pm}_i(u)\h^{\mp}_j(v)&=
\frac{u_{\mp}-v_{\pm}}{qu_{\mp}-q^{-1}v_{\pm}}\h^{\mp}_j(v)\h^{\pm}_i(u)
\label{hihjmp}
\end{align}
for $i<j$ and $i\neq n$, and
\ben
\frac{u_{\pm}-v_{\mp}}{q^2u_{\pm}-q^{-2}v_{\mp}}\h^{\pm}_n(u)\h^{\mp}_{n+1}(v)=
\frac{u_{\mp}-v_{\pm}}{q^2u_{\mp}-q^{-2}v_{\pm}}\h^{\mp}_{n+1}(v)\h^{\pm}_n(u).
\een
The relations
involving $\h^{\pm}_i(u)$ and $\Xc_{j}^{\pm}(v)$ are
\ben
\bal
\h_{i}^{\pm}(u)\Xc_{j}^{+}(v)&
=\frac{u-v_{\pm}}{q^{(\ep_i,\alpha_j)}u-q^{-(\ep_i,\alpha_j)}v_{\pm}}
\Xc_{j}^{+}(v)\h_{i}^{\pm}(u),\\[0.4em]
\h_{i}^{\pm}(u)\Xc_{j}^{-}(v)&
=\frac{q^{(\ep_i,\alpha_j)}u_{\pm}-q^{-(\ep_i,\alpha_j)}v}{u_{\pm}-v}
\Xc_{j}^{-}(v)\h_{i}^{\pm}(u)
\eal
\een
for $i\neq n+1$, together with
\ben
\bal
\h_{n+1}^{\pm}(u)\Xc_n^{+}(v)&=
\frac{u_{\mp}-v}{q^{-2}u_{\mp}-q^2v}\Xc^{+}_n(v)\h^{\pm}_{n+1}(u),\\
\h_{n+1}^{\pm}(u)\Xc_n^{-}(v)&=
\frac{q^{-2}u_{\pm}-q^2v}{u_{\pm}-v}\Xc_n^{-}(v)\h_{n+1}^{\pm}(u),
\eal
\een
and
\ben
\bal
\h_{n+1}^{\pm}(u)^{-1}\Xc_{n-1}^{+}(v)\h_{n+1}^{\pm}(u)&=
\frac{q^{-1}u_{\mp}-qv}{q^{-2}u_{\mp}-q^{2}v}\Xc^{+}_{n-1}(v),\\
\h_{n+1}^{\pm}(u)\Xc_{n-1}^{-}(v)\h_{n+1}^{\pm}(u)^{-1}&=
\frac{q^{-1}u_{\pm}-qv}{q^{-2}u_{\pm}-q^2v}\Xc_{n-1}^{-}(v),
\eal
\een
while
\ben
\bal
\h_{n+1}^{\pm}(u)\Xc_{i}^{+}(v)
&=\Xc_{i}^{+}(v)\h_{n+1}^{\pm}(u),\\[0.4em]
\h_{n+1}^{\pm}(u)\Xc_{i}^{-}(v)
&=\Xc_{i}^{-}(v)\h_{n+1}^{\pm}(u),
\eal
\een
for $1\leqslant i\leqslant n-2$. For the relations involving $\Xc^{\pm}_i(u)$ we have
\ben
(u-q^{\pm (\alpha_i,\alpha_j)}v)\Xc_{i}^{\pm}(uq^{i})\Xc_{j}^{\pm}(vq^j)
=(q^{\pm (\alpha_i,\alpha_j)}u-v) \Xc_{j}^{\pm}(vq^j)\Xc_{i}^{\pm}(uq^i)
\een
for $i,j=1,\dots,n-1$;
\ben
(u-q^{\pm (\alpha_i,\alpha_n)}v)\Xc_{i}^{\pm}(uq^i)\Xc_{n}^{\pm}(vq^{n+1})
=(q^{\pm (\alpha_i,\alpha_n)}u-v) \Xc_{n}^{\pm}(vq^{n+1})\Xc_{i}^{\pm}(uq^i)
\een
for $i=1,\dots,n-1$;
\ben
(u-q^{\pm (\alpha_n,\alpha_n)}v)\Xc_{n}^{\pm}(u)\Xc_{n}^{\pm}(v)
=(q^{\pm (\alpha_n,\alpha_n)}u-v) \Xc_{n}^{\pm}(v)\Xc_{n}^{\pm}(u)
\een

and
\ben
[\Xc_i^{+}(u),\Xc_j^{-}(v)]=
\delta_{ij}(q_i-q_i^{-1})\Big(\delta\big(u\ts q^{-c}/v\big)\h_i^{-}(v_+)^{-1}\h_{i+1}^{-}(v_+)
-\delta\big(u\ts q^{c}/v\big)\h_i^{+}(u_+)^{-1}\h_{i+1}^{+}(u_+)\Big)
\een
together with the
Serre relations
\beql{serrex}
\sum_{\pi\in \Sym_{r}}\sum_{l=0}^{r}(-1)^l{{r}\brack{l}}_{q_i}
  \Xc^{\pm}_{i}(u_{\pi(1)})\dots \Xc^{\pm}_{i}(u_{\pi(l)})
  \Xc^{\pm}_{j}(v)\tss \Xc^{\pm}_{i}(u_{\pi(l+1)})\dots \Xc^{\pm}_{i}(u_{\pi(r)})=0,
\eeq
which hold for all $i\neq j$ and we set $r=1-A_{ij}$.
\eth

\bpf
We only need to verify the relations complementary to those produced in
Secs~\ref{subsec:typea} and \ref{subsec:rlr}. The additional relations are
verified by very similar arguments to those used
in \cite{df:it} and so we give only a few details illustrating the calculations
which are more specific to type $C$. We start with \eqref{hihjmp} and take $j=n+1$.
By using Corollary~\ref{cor:commu}, for $i<n$ we deduce
\begin{multline}
\label{hiln+1C}
\frac{u_{\pm}-v_{\mp}}{qu_{\pm}-q^{-1}v_{\mp}}\h^{\pm}_i(u)\big(\h^{\mp}_{n+1}(v)
+\f_{n+1,n}^{\mp}(v)\h_{n}^{\mp}(v)\e_{n,n+1}^{\mp}(v)\big)\\
{}=\frac{u_{\mp}-v_{\pm}}{qu_{\mp}-q^{-1}v_{\pm}}\big(\h^{\pm}_{n+1}(v)
+\f_{n+1,n}^{\mp}(v)\h_{n}^{\mp}(v)\e_{n,n+1}^{\mp}(v)\big)\h^{\pm}_i(u)
\end{multline}
and
\ben
\begin{aligned}
\frac{u_{\pm}-v_{\mp}}{qu_{\pm}-q^{-1}v_{\mp}}\h^{\pm}_i(u)\f_{n+1,n}^{\mp}(v)\h_{n}^{\mp}(v)=
\frac{u_{\mp}-v_{\pm}}{qu_{\mp}-q^{-1}v_{\pm}}\f_{n+1,n}^{\mp}(v)\h_{n}^{\mp}(v)\h^{\pm}_i(u).
\end{aligned}
\een
Hence, the left hand side of \eqref{hiln+1C} equals
\beql{lghhh}
 \frac{u_{\pm}-v_{\mp}}{qu_{\pm}-q^{-1}v_{\mp}}\h^{\pm}_i(u)\h^{\mp}_{n+1}(v)
+\frac{u_{\mp}-v_{\pm}}{qu_{\mp}-q^{-1}v_{\pm}}\f_{n+1,n}^{\mp}(v)
\h_{n}^{\mp}(v)\h^{\pm}_i(u)\e_{n,n+1}^{\mp}(v)).
\eeq
By Corollary \ref{cor:commu},
\ben
\begin{aligned}
\frac{u_{\pm}-v_{\mp}}{qu_{\pm}-q^{-1}v_{\mp}}\h^{\pm}_i(u)\h_{n}^{\mp}(v)=
\frac{u_{\mp}-v_{\pm}}{qu_{\mp}-q^{-1}v_{\pm}}\h_{n}^{\mp}(v)\h^{\pm}_i(u),
\end{aligned}
\een
so that \eqref{lghhh} can be written as
\ben
 \frac{u_{\pm}-v_{\mp}}{qu_{\pm}-q^{-1}v_{\mp}}\h^{\pm}_i(u)\h^{\mp}_{n+1}(v)
+\frac{u_{\pm}-v_{\mp}}{qu_{\pm}-q^{-1}v_{\mp}}\f_{n+1,n}^{\mp}(v)\h^{\pm}_i(u)
\h_{n}^{\mp}(v)\e_{n,n+1}^{\mp}(v)).
\een
Using Corollary \ref{cor:commu} once again, we find
\ben
\begin{aligned}
\frac{u_{\pm}-v_{\mp}}{qu_{\pm}-q^{-1}v_{\mp}}\h^{\pm}_i(u)\h_{n}^{\mp}(v)\e_{n,n+1}^{\mp}(v)=
\frac{u_{\mp}-v_{\pm}}{qu_{\mp}-q^{-1}v_{\pm}}\h_{n}^{\mp}(v)\e_{n,n+1}^{\mp}(v)\h^{\pm}_i(u)
\end{aligned}
\een
for $i=1,2,\dots,n-1$,
and so the left hand side of \eqref{hiln+1C} takes the form
\ben
 \frac{u_{\pm}-v_{\mp}}{qu_{\pm}-q^{-1}v_{\mp}}\h^{\pm}_i(u)\h^{\mp}_{n+1}(v)+
 \frac{u_{\mp}-v_{\pm}}{qu_{\mp}-q^{-1}v_{\pm}}\f_{n+1,n}^{\mp}(v)
 \h_{n}^{\mp}(v)\e_{n,n+1}^{\mp}(v))\h^{\pm}_i(u),
\een
which implies \eqref{hihjmp} with $j=n+1$.

The relations
involving $\h^{\pm}_{n+1}(v)$ and $\Xc_{i}^{\pm}(u)$ with $i=1,2,\dots,n-2$
are implied by the following:
\begin{alignat}{2}
\label{eihn+1pm}
\e_{i,i+1}^{\pm}(u)\h_{n+1}^{\pm}(v)&=\h_{n+1}^{\pm}(v)\e_{i,i+1}^{\pm}(u),\qquad
\e_{i,i+1}^{\pm}(u)\h_{n+1}^{\mp}(v)&&=\h_{n+1}^{\mp}(v)\e_{i,i+1}^{\pm}(u),\\[0.4em]
\f_{i+1,i}^{\pm}(u)\h_{n+1}^{\pm}(v)&=\h_{n+1}^{\pm}(v)\f_{i+1,i}^{\pm}(u),\qquad
\f_{i+1,i}^{\pm}(u)\h_{n+1}^{\mp}(v)&&=\h_{n+1}^{\mp}(v)\f_{i+1,i}^{\pm}(u).
\non
\end{alignat}
We will verify
the second relation in \eqref{eihn+1pm}.
Corollary~\ref{cor:commu} implies
\begin{multline}
\label{lii+1ln+1C}
\frac{u_{\pm}-v_{\mp}}{qu_{\pm}-q^{-1}v_{\mp}}\h_i^{\pm}(u)\e_{i,i+1}^{\pm}(u)\big(\h^{\mp}_{n+1}(v)
+\f_{n+1,n}^{\mp}(v)\h_{n}^{\mp}(v)\e_{n,n+1}^{\mp}(v)\big)\\
=\frac{u_{\mp}-v_{\pm}}{qu_{\mp}-q^{-1}v_{\pm}}\big(\h^{\mp}_{n+1}(v)
+\f_{n+1,n}^{\mp}(v)\h_{n}^{\mp}(v)\e_{n,n+1}^{\mp}(v)\big)\h_i^{\pm}(u)\e_{i,i+1}^{\pm}(u).
\end{multline}
Moreover, by Corollary \ref{cor:commu}, we have the following relations:
\ben
\bal
\frac{u_{\pm}-v_{\mp}}{qu_{\pm}-q^{-1}v_{\mp}}\h_i^{\pm}(u)
\e_{i,i+1}^{\pm}(u)\f_{n+1,n}^{\mp}(v)\h_{n}^{\mp}(v)
{}&=\frac{u_{\mp}-v_{\pm}}{qu_{\mp}-q^{-1}v_{\pm}}
\f_{n+1,n}^{\mp}(v)\h_{n}^{\mp}(v)
\h_i^{\pm}(u)\e_{i,i+1}^{\pm}(u),\\
\frac{u_{\pm}-v_{\mp}}{qu_{\pm}-q^{-1}v_{\mp}}\h_i^{\pm}(u)\e_{i,i+1}^{\pm}(u)\h_{n}^{\mp}(v)
{}&=\frac{u_{\mp}-v_{\pm}}{qu_{\mp}-q^{-1}v_{\pm}}\h_{n}^{\mp}(v)
\h_i^{\pm}(u)\e_{i,i+1}^{\pm}(u),\\
\frac{u_{\pm}-v_{\mp}}{qu_{\pm}-q^{-1}v_{\mp}}\h_i^{\pm}(u)
\e_{i,i+1}^{\pm}(u)\h_{n}^{\mp}(v)\e_{n,n+1}^{\mp}(v)
&=\frac{u_{\mp}-v_{\pm}}{qu_{\mp}-q^{-1}v_{\pm}}\h_{n}^{\mp}(v)\e_{n,n+1}^{\mp}(v)
\h_i^{\pm}(u)\e_{i,i+1}^{\pm}(u).
\eal
\een
Thus, the left hand side of \eqref{lii+1ln+1C} equals
\ben
\frac{u_{\pm}-v_{\mp}}{qu_{\pm}-q^{-1}v_{\mp}}\h^{\mp}_{n+1}(v)+
\frac{u_{\mp}-v_{\pm}}{qu_{\mp}-q^{-1}v_{\pm}}\f_{n+1,n}^{\mp}(v)
 \h_{n}^{\mp}(v)\e_{n,n+1}^{\mp}(v)\h_i^{\pm}(u)\e_{i,i+1}^{\pm}(u)
\een
and so,
\ben
\frac{u_{\pm}-v_{\mp}}{qu_{\pm}-q^{-1}v_{\mp}}\h_i^{\pm}(u)\e_{i,i+1}^{\pm}(u)\h^{\mp}_{n+1}(v)
=\frac{u_{\mp}-v_{\pm}}{qu_{\mp}-q^{-1}v_{\pm}}\h^{\mp}_{n+1}(v)\h_i^{\pm}(u)\e_{i,i+1}^{\pm}(u).
\een
Using now \eqref{hihjmp} with $j=n+1$ we get the second relation in
\eqref{eihn+1pm}.

The remaining cases of the type $C$-specific relations
involving $\Xc_{i}^{\pm}(u)$ and $\h^{\pm}_{j}(v)$ are deduced
with the use of Remark~\ref{re:anotherA}, Lemma~\ref{lem:eiprei} and Corollary~\ref{cor:commu}.
In particular, Remark~\ref{re:anotherA} and the corresponding type $A$ relations in \cite{df:it}
imply
\ben
\h_{n'}^{\pm}(u)^{-1}\e_{n',(n-1)'}^{\pm}(v)\h_{n'}^{\pm}(u)
=\frac{qu-q^{-1}v}{u-v}\e_{n',(n-1)'}^{\pm}(v)
-\frac{(q-q^{-1})v}{u-v}\e_{n',(n-1)'}^{\pm}(u).
\een
By Lemma~\ref{lem:eiprei}, we can write this relation as
\ben
\h_{n'}^{\pm}(u)^{-1}\e_{n-1,n}^{\pm}(vq^{-4})\h_{n'}^{\pm}(u)
=\frac{qu-q^{-1}v}{u-v}\e_{n-1,n}^{\pm}(vq^{-4})
-\frac{(q-q^{-1})v}{u-v}\e_{n-1,n}^{\pm}(uq^{-4}),
\een
which leads to the relations
involving $\Xc_{n-1}^{\pm}(u)$ and $\h^{\pm}_{n+1}(v)$.

Now turn to the relations between the series $\Xc_{i}^{\pm}(u)$. For $i<n-1$ we have
\begin{alignat}{2}
\e_{i,i+1}^{\pm}(u)\e_{n,n+1}^{\pm}(v)&=\e_{n,n+1}^{\pm}(v)\e_{i,i+1}^{\pm}(u),
\qquad
\e_{i,i+1}^{\pm}(u)\e_{n,n+1}^{\mp}(v)&&=\e_{n,n+1}^{\mp}(v)\e_{i,i+1}^{\pm}(u),
\non\\[0.4em]
\f_{i+1,i}^{\pm}(u)\f_{n+1,n}^{\pm}(v)&=\f_{n+1,n}^{\pm}(v)\f_{i+1,i}^{\pm}(u),
\qquad
\f_{i+1,i}^{\pm}(u)\f_{n+1,n}^{\mp}(v)&&=\f_{n+1,n}^{\mp}(v)\f_{i+1,i}^{\pm}(u).
\non
\end{alignat}
This is verified with the same use of Corollary~\ref{cor:commu} as above.
Therefore, for all $i< n-1$ we get
$\Xc_i^{\pm}(u)\Xc^{\pm}_n(v)=\Xc^{\pm}_n(v)\Xc_i^{\pm}(u)$.
For the relation involving $\Xc_{n-1}^{\pm}(u)$ and $\Xc^{\pm}_n(v)$ it will be sufficient to
consider the case $n=2$ and then apply Theorem~\ref{thm:red}.
We find from the defining
relations \eqref{gen rel2} and \eqref{gen rel3} that
the expression
\beql{l12l23C}
\frac{u_{\pm}-v_{\mp}}{qu_{\pm}-q^{-1}v_{\mp}}\ell_{12}^{\pm}(u)\ell_{23}^{\mp}(v)
+\frac{(q-q^{-1})u_{\pm}}{qu_{\pm}-q^{-1}v_{\mp}}\ell_{22}^{\pm}(u)\ell_{13}^{\mp}(v)
\eeq
equals
\ben
\bal
\frac{(u_{\mp}-v_{\pm})q^{-1}(u_{\mp}-q^{-4}v_{\pm})}{(qu_{\mp}-q^{-1}v_{\pm})(u_{\mp}-q^{-6}v_{\pm})}
\ell_{23}^{\mp}(v)\ell_{12}^{\pm}(u)
&+\frac{(q-q^{-1})v_{\pm}}{qu_{\mp}-q^{-1}v_{\pm}}
\ell_{22}^{\mp}(v)\ell_{13}^{\pm}(u)\\
+\frac{(q-q^{-1})(u_{\mp}-v_{\pm})v_{\pm}q^{-3}}{(qu_{\mp}-q^{-1}v_{\pm})(u_{\mp}-q^{-6}v_{\pm})}
\ell_{21}^{\mp}(v)\ell_{14}^{\pm}(u)&+
\frac{(q-q^{-1})(u_{\mp}-v_{\pm})v_{\pm}q^{-4}}{(qu_{\mp}-q^{-1}v_{\pm})(u_{\mp}-q^{-6}v_{\pm})}
\ell_{22}^{\mp}(v)\ell_{13}^{\pm}(u)\\[0.3em]
&-\frac{(q-q^{-1})(u_{\mp}-v_{\pm})q^{-1}u_{\mp}}{(qu_{\mp}-q^{-1}v_{\pm})(u_{\mp}-q^{-6}v_{\pm})}
\ell_{24}^{\mp}(v)\ell_{11}^{\pm}(u).
\eal
\een
On the other hand, applying the formula for $\ell_{23}^{\mp}(v)$ arising from
the Gauss decomposition, we can write \eqref{l12l23C} as
\ben
\bal
\frac{u_{\pm}-v_{\mp}}{qu_{\pm}-q^{-1}v_{\mp}}
\ell_{12}^{\pm}(u)\h_{2}^{\mp}(v)\e_{23}^{\mp}(v)&+
\frac{u_{\pm}-v_{\mp}}{qu_{\pm}-q^{-1}v_{\mp}}
\ell_{12}^{\pm}(u)\f_{21}^{\mp}(v)
\h_{1}^{\mp}(v)\e_{13}^{\mp}(v)\\
{}&+\frac{(q-q^{-1})u_{\pm}}{qu_{\pm}-q^{-1}v_{\mp}}
\ell_{22}^{\pm}(u)\ell_{13}^{\mp}(v).
\eal
\een
By using the defining relations between the series $\ell_{12}^{+}(u)$ and $\ell_{21}^{-}(v)$
we can bring \eqref{l12l23C} to the form
\ben
\bal
\frac{u_{\pm}-v_{\mp}}{qu_{\pm}-q^{-1}v_{\mp}}\ell_{12}^{\pm}(u)\h_{2}^{\mp}(v)\e_{23}^{\mp}(v)&+
\frac{u_{\mp}-v_{\pm}}{qu_{\mp}-q^{-1}v_{\pm}}\ell_{21}^{\mp}(v)\ell_{12}^{\pm}(u)\e_{13}^{\mp}(v)\\
&+\frac{(q-q^{-1})u_{\mp}}{qu_{\mp}-q^{-1}v_{\pm}}\ell_{22}^{\mp}(v)\ell_{11}^{\pm}(u)\e_{13}^{\mp}(v).
\eal
\een
Further, by using the defining relations between $\ell_{12}^{\pm}(u)$ and $\ell_{11}^{\mp}(v)$
we can write \eqref{l12l23C} as
\ben
\frac{u_{\pm}-v_{\mp}}{qu_{\pm}-q^{-1}v_{\mp}}\ell_{12}^{\pm}(u)\h_{2}^{\mp}(v)\e_{23}^{\mp}(v)+
\f_{21}^{\mp}(v)\ell_{12}^{\pm}(u)\ell_{13}^{\mp}(v)
+\frac{(q-q^{-1})u_{\mp}}{qu_{\mp}-q^{-1}v_{\pm}}\h_{2}^{\mp}(v)\ell_{11}^{\pm}(u)\e_{13}^{\mp}(v).
\een
As a next step, apply
the relations between $\ell_{12}^{\pm}(u)$ and $\ell_{13}^{\mp}(v)$ to write the sum
\ben
\frac{u_{\pm}-v_{\mp}}{qu_{\pm}-q^{-1}v_{\mp}}\ell_{12}^{\pm}(u)\h_{2}^{\mp}(v)\e_{23}^{\mp}(v)
+\frac{(q-q^{-1})u_{\mp}}{qu_{\mp}-q^{-1}v_{\pm}}\h_{2}^{\mp}(v)\ell_{11}^{\pm}(u)\e_{13}^{\mp}(v)
\een
as
\begin{multline}
\frac{(u_{\mp}-v_{\pm})q^{-1}(u_{\mp}-q^{-4}v_{\pm})}{(qu_{\mp}-q^{-1}v_{\pm})(u_{\mp}-q^{-6}v_{\pm})}
\h_{2}^{\mp}(v)\e_{23}^{\mp}(v)\ell_{12}^{\pm}(u)+
\frac{(q-q^{-1})v_{\pm}}{qu_{\mp}-q^{-1}v_{\pm}}\h_{2}^{\mp}(v)\ell_{13}^{\pm}(u)\\
{}+\frac{(q-q^{-1})(u_{\mp}-v_{\pm})v_{\pm}q^{-4}}{(qu_{\mp}-q^{-1}v_{\pm})(u_{\mp}-q^{-6}v_{\pm})}
\h_{2}^{\mp}(v)\ell_{13}^{\pm}(u)\\
-\frac{(q-q^{-1})(u_{\mp}-v_{\pm})q^{-1}u_{\mp}}{(qu_{\mp}-q^{-1}v_{\pm})(u_{\mp}-q^{-6}v_{\pm})}
\h_{2}^{\mp}(v)\e_{24}(v)\ell_{11}^{\pm}(u).
\non
\end{multline}
Using \eqref{ELMPj=l} and relations between $\h_1^{\pm}(u)$ and $\h_2^{\mp}(u)$, we get the relations:
\ben
\bal
\ell_{12}^{\pm}(u)\h_2^{\mp}&(v)=\frac{qu_{\mp}-q^{-1}v_{\pm}}{u_{\mp}-v_{\pm}}
\h_1^{\pm}(u)\h_2^{\mp}(v)\e_{12}^{\pm}(u)
-\frac{(q-q^{-1})u_{\mp}}{u_{\mp}-v_{\pm}}\h_1^{\pm}(u)\h_2^{\mp}(v)\e_{12}^{\mp}(v)\\
&=\frac{qu_{\pm}-q^{-1}v_{\mp}}{u_{\pm}-v_{\mp}}\h_2^{\mp}(v)\h_1^{\pm}(u)\e_{12}^{\pm}(u)-
\frac{(q-q^{-1})u_{\mp}}{qu_{\mp}-q^{-1}v_{\pm}}\frac{qu_{\pm}-q^{-1}v_{\mp}}{u_{\pm}-v_{\mp}}
\h_2^{\mp}(v)\h_1^{\pm}(u)\e_{12}^{\mp}(v).
\eal
\een
Hence, since the series $\h_1^{\pm}(u)$ and $\h_2^{\mp}(v)$
are invertible, and $[\h_1^{\pm}(u),\e_{23}^{\mp}(v)]=0$,
we come to the relation
\begin{align}\label{e12e23C}
&\e_{12}^{\pm}(u)\e_{23}^{\mp}(v)
+\frac{(q-q^{-1})u_{\mp}}{qu_{\mp}-q^{-1}v_{\pm}}\e_{13}^{\mp}(v)
-\frac{(q-q^{-1})u_{\mp}}{qu_{\mp}-q^{-1}v_{\pm}}\e_{12}^{\mp}(v)\e_{23}^{\mp}(v)\\
&{}=
\frac{(u_{\mp}-v_{\pm})q^{-1}(u_{\mp}-q^{-4}v_{\pm})}{(qu_{\mp}-q^{-1}v_{\pm})(u_{\mp}-q^{-6}v_{\pm})}
\e_{23}^{\mp}(v)\e_{12}^{\pm}(u)+
\frac{(q-q^{-1})v_{\pm}}{qu_{\mp}-q^{-1}v_{\pm}}\e_{13}^{\pm}(u)
\non\\
&{}+\frac{(q-q^{-1})(u_{\mp}-v_{\pm})v_{\pm}q^{-4}}{(qu_{\mp}-q^{-1}
v_{\pm})(u_{\mp}-q^{-6}v_{\pm})}\e_{13}^{\pm}(u)
-\frac{(q-q^{-1})(u_{\mp}-v_{\pm})q^{-1}u_{\mp}}{(qu_{\mp}-q^{-1}v_{\pm})(u_{\mp}-q^{-6}v_{\pm})}
\h_1^{\pm}(u)^{-1}\e^{\mp}_{24}(v)\h_{1}^{\pm}(u).
\non
\end{align}
Similar arguments imply the relations
\ben
\begin{aligned}
\frac{u_{\pm}-v_{\mp}}{qu_{\pm}-q^{-1}v_{\mp}}\h_{1}^{\pm}(u)\h_{2}^{\mp}(v)\e_{24}^{\mp}(v)&=
\frac{(u_{\mp}-v_{\pm})q^{-1}(u_{\mp}-v_{\pm}q^{-4})}{(qu_{\mp}-q^{-1}v_{\pm})(u_{\mp}-q^{-6}v_{\pm})}
\h_{2}^{\mp}(v)\e_{24}^{\mp}(v)\h_{1}^{\pm}(u)\\
&+\frac{(q-q^{-1})(u_{\mp}-v_{\pm})v_{\pm}q^{-3}}{(qu_{\mp}-q^{-1}v_{\pm})(u_{\mp}-q^{-6}v_{\pm})}
\h_{2}^{\mp}(v)\h_{1}^{\pm}(u)\e_{13}^{\pm}(u)\\
&-\frac{(q-q^{-1})(u_{\mp}-v_{\pm})v_{\pm}q^{-5}}{(qu_{\mp}-q^{-1}v_{\pm})(u_{\mp}-q^{-6}v_{\pm})}
\h_{2}^{\mp}(v)\e^{\mp}_{23}(v)\h_{1}^{\pm}(u)\e_{12}^{\pm}(u).
\end{aligned}
\een
Due to the relation
\ben
\frac{u_{\pm}-v_{\mp}}{qu_{\pm}-q^{-1}v_{\mp}}\h_{1}^{\pm}(u)\h_{2}^{\mp}(v)
=\frac{u_{\mp}-v_{\pm}}{qu_{\mp}-q^{-1}v_{\pm}}\h_{2}^{\mp}(v)\h_{1}^{\pm}(u),
\een
we come to
\begin{multline}\non
\h_{1}^{\pm}(u)\e_{24}^{\mp}(v)=
\frac{q^{-1}(u_{\mp}-v_{\pm}q^{-4})}{u_{\mp}-q^{-6}v_{\pm}}\e_{24}^{\mp}(v)\h_{1}^{\pm}(u)\\
{}+\frac{(q-q^{-1})v_{\pm}q^{-3}}{u_{\mp}-q^{-6}v_{\pm}}\h_{1}^{\pm}(u)\e_{13}^{\pm}(u)
-\frac{(q-q^{-1})v_{\pm}q^{-5}}{u_{\mp}-q^{-6}v_{\pm}}\e^{\mp}_{23}(v)\h_{1}^{\pm}(u)\e_{12}^{\pm}(u).
\end{multline}
Together with the expression in \eqref{e12e23C}, this gives
\ben
\bal
\e_{12}^{\pm}(u)\e_{23}^{\mp}(v)&=\frac{1}{q^2u_{\mp}-q^{-2}v_{\pm}}
\big((q^2-q^{-2})v_{\pm}\e_{13}^{\pm}(u)+(u_{\mp}-v_{\pm})\e_{23}^{\mp}(v)\e_{12}^{\pm}(u)\\[0.5em]
&-(q^2-q^{-2})u_{\mp}\e_{12}^{\mp}(v)\e_{23}^{\mp}(v)+(q^2-q^{-2})u_{\mp}\e_{13}^{\mp}(v)\big)\\[0.5em]
&+\frac{(q-q^{-1})u_{\mp}(u_{\mp}-v_{\pm})}{(q^2u_{\mp}-q^{-2}v_{\pm})(qu_{\mp}-q^{-1}v_{\pm})}
\big(\e_{12}^{\mp}(v)\e_{23}^{\mp}(v)
-\e_{13}^{\mp}(v)-q^2\e_{24}^{\mp}(v)
\big).
\eal
\een
Multiply both sides by $qu_{\mp}-q^{-1}v_{\pm}$ and set $qu_{\mp}=q^{-1}v_{\pm}$
to see that the second summand
vanishes. Finally, by applying Theorem~\ref{thm:red},
we come to the relation
\ben
\bal
(q^2u_{\mp}-q^{-2}v_{\pm})&\e_{n-1,n}^{\pm}(u)\e_{n,n+1}^{\mp}(v)\\[0.4em]
&=(u_{\mp}-v_{\pm})\e_{n,n+1}^{\mp}(v)\e_{n-1,n}^{\pm}(u)
+(q^2-q^{-2})v_{\pm}\e_{n-1,n+1}^{\pm}(u)\\[0.4em]
&-(q^2-q^{-2})u_{\mp}\e_{n-1,n}^{\mp}(v)\e_{n,n+1}^{\mp}(v)
+(q^2-q^{-2})u_{\mp}\e_{n-1,n+1}^{\mp}(v).
\eal
\een
Quite similar calculations lead to its counterparts:
\ben
\bal
(q^2u-q^{-2}v)&\e_{n-1,n}^{\pm}(u)\e_{n,n+1}^{\pm}(v)\\[0.4em]
&=(u-v)\e_{n,n+1}^{\pm}(v)\e_{n-1,n}^{\pm}(u)+(q^2-q^{-2})v\e_{n-1,n+1}^{\pm}(u)\\[0.4em]
&-(q^2-q^{-2})u\e_{n-1,n}^{\pm}(v)\e_{n,n+1}^{\pm}(v)
+(q^2-q^{-2})u\e_{n-1,n+1}^{\pm}(v),
\eal
\een

\ben
\bal
(u_{\pm}-v_{\mp})&\f_{n,n-1}^{\pm}(u)\f_{n+1,n}^{\mp}(v)\\[0.4em]
&=(q^2u_{\pm}-q^{-2}v_{\mp})\f_{n+1,n}^{\mp}(v)\f_{n,n-1}^{\pm}(u)
-(q^2-q^{-2})v_{\mp}\f_{n+1,n-1}^{\mp}(v)
\\[0.4em]
&+(q^2-q^{-2})v_{\mp}\f_{n+1,n}^{\mp}(v)\f_{n,n-1}^{\mp}(v)+(q^2-q^{-2})u_{\pm}\f_{n+1,n-1}^{\pm}(u)
\eal
\een
and
\ben
\bal
(u-v)\f_{n,n-1}^{\pm}(u)\f_{n+1,n}^{\pm}(v)
&=(q^2u-q^{-2}v)\f_{n+1,n}^{\pm}(v)\f_{n,n-1}^{\pm}(u)-(q^2-q^{-2})v\f_{n+1,n-1}^{\pm}(v)\\[0.4em]
&+(q^2-q^{-2})v\f_{n+1,n}^{\pm}(v)\f_{n,n-1}^{\pm}(v)+(q^2-q^{-2})u\f_{n+1,n-1}^{\pm}(u).
\eal
\een
As a consequence, we have thus verified the relations
\ben
(uq^{2}-q^{-2}v)^{\pm 1}\Xc_{n-1}^{\pm}(u)\Xc^{\pm}_n(v)
=(u-v)^{\pm 1}\Xc^{\pm}_n(v)\Xc_{n-1}^{\pm}(u).
\een

Now suppose that $i\leqslant n-1$ and verify the relations
\begin{alignat}{2}\label{eifnpm}
\e_{i,i+1}^{\pm}(u)\f_{n+1,n}^{\pm}(v)&=\f_{n+1,n}^{\pm}(v)\e_{i,i+1}^{\pm}(u),\qquad
\e_{i,i+1}^{\pm}(u)\f_{n+1,n}^{\mp}(v)&&=\f_{n+1,n}^{\mp}(v)\e_{i,i+1}^{\pm}(u),\\
\f_{i+1,i}^{\pm}(u)\e_{n,n+1}^{\pm}(v)&=\e_{n,n+1}^{\pm}(v)\f_{i+1,i}^{\pm}(u),\qquad
\f_{i+1,i}^{\pm}(u)\e_{n,n+1}^{\mp}(v)&&=\e_{n,n+1}^{\mp}(v)\f_{i+1,i}^{\pm}(u).
\non
\end{alignat}
We only do the second relation in \eqref{eifnpm} as the arguments are quite similar.
If $i\leqslant n-2$ then Corollary \ref{cor:commu} gives
\ben
\frac{u_{\pm}-v_{\mp}}{qu_{\pm}-q^{-1}v_{\mp}}
\h_i^{\pm}(u)\e_{i,i+1}^{\pm}(u)\f_{n+1,n}^{\mp}(v)\h_{n}^{\mp}(v)
=\frac{u_{\mp}-v_{\pm}}{qu_{\mp}-q^{-1}v_{\pm}}
\f_{n+1,n}^{\mp}(v)\h_{n}^{\mp}(v)\h_i^{\pm}(u)\e_{i,i+1}^{\pm}(u)
\een
and
\ben
\frac{u_{\pm}-v_{\mp}}{qu_{\pm}-q^{-1}v_{\mp}}\h_i^{\pm}(u)\f_{n+1,n}^{\mp}(v)\h_{n}^{\mp}(v)
=\frac{u_{\mp}-v_{\pm}}{qu_{\mp}-q^{-1}v_{\pm}}\f_{n+1,n}^{\mp}(v)\h_{n}^{\mp}(v)\h_i^{\pm}(u).
\een
Therefore,
\ben
\h_i^{\pm}(u)\e_{i,i+1}^{\pm}(u)\f_{n+1,n}^{\mp}(v)\h_{n}^{\mp}(v)
=\h_i^{\pm}(u)\f_{n+1,n}^{\mp}(v)\h_{n}^{\mp}(v)\e_{i,i+1}^{\pm}(u).
\een
Since $[\h_{n}^{\mp}(v),\e_{i,i+1}^{\pm}(u)]=0$, the required relation follows.
Now let $i=n-1$. Due to \eqref{ELMPj=l}, we have
\ben
\bal
\e_{n-1,n}^{\pm}(u)\f_{n+1,n}^{\mp}(v)\h_n^{\mp}(v)
&=\frac{qu_{\mp}-q^{-1}v_{\pm}}{u_{\mp}-v_{\pm}}\f_{n+1,n}^{\mp}(v)\h_{n}^{\mp}(v)\e_{n-1,n}^{\pm}(u)\\
&+\frac{(q-q^{-1})v_{\pm}}{u_{\mp}-v_{\pm}}\f_{n+1,n}^{\mp}(v)\h_{n}^{\mp}(v)\e_{n-1,n}^{\mp}(v)
\eal
\een
and
\ben
\e_{n-1,n}^{\pm}(u)\h_n^{\mp}(v)=
\frac{qu_{\mp}-q^{-1}v_{\pm}}{u_{\mp}-v_{\pm}}\h_{n}^{\mp}(v)\e_{n-1,n}^{\pm}(u)
+\frac{(q-q^{-1})v_{\pm}}{u_{\mp}-v_{\pm}}\h_{n}^{\mp}(v)\e_{n-1,n}^{\mp}(v).
\een
Hence
\ben
\e_{n-1,n}^{\pm}(u)\f_{n+1,n}^{\mp}(v)\h_n^{\mp}(v)=
\f_{n+1,n}^{\mp}(v)\e_{n-1,n}^{\pm}(u)\h_n^{\mp}(v),
\een
so that the second relation in \eqref{eifnpm} is verified.
Thus, by applying Proposition~\ref{LowC} we thus derive all cases
for the commutator formula for the series $\Xc_i^{+}(u)$ and $\Xc_j^{-}(v)$.

To complete the proof of the theorem, we will now verify the Serre relations
\eqref{serrex}. By Proposition~\ref{prop:corrgauss}, these relations have the same form
for the algebras $U(R)$ and $U(\overline{R})$. We will work with the algebra $U(R)$
and introduce its elements
$x_{i,m}^{\pm}$ and $a_{i,l}$ for $i=1,\dots,n$ and
$m,l\in\ZZ$ with $l\ne 0$
by the formulas
\ben
\bal
x^{\pm}_{i}(u)&= (q_i-q_i^{-1})^{-1}X^{\pm}_i(uq^i),\\[0.4em]
\psi_{i}(u)&= h^{-}_{i+1}(uq^i)\ts h^{-}_{i}(uq^i)^{-1},\\[0.4em]
\varphi_{i}(u)&= h^{+}_{i+1}(uq^i)\ts h^{+}_{i}(uq^i)^{-1},
\eal
\een
for $i=1,\dots,n-1$, and
\ben
\bal
x^{\pm}_{n}(u)&= (q_n-q_n^{-1})^{-1}X^{\pm}_n(uq^{n+1}),\\[0.4em]
\psi_{n}(u)&= h^{-}_{n+1}(uq^{n+1})\ts h^{-}_{n}(uq^{n+1})^{-1},\\[0.4em]
\varphi_{n}(u)&= h^{+}_{n+1}(uq^{n+1})\ts h^{+}_{n}(uq^{n+1})^{-1},
\eal
\een
and the expansions \eqref{psiiu}, \eqref{phiiu} and \eqref{xpm},
where $k_i=h_{i\tss 0}^{-}\ts h_{i+1\tss 0}^{+}$ and $h^{\pm}_{i\tss 0}$ denotes the constant term
of the series \eqref{hmqua}. In terms of the elements $x_{i,m}^{\pm}$ the Serre relations
take the form
\beql{serrexx}
\sum_{\pi\in \Sym_{r}}\sum_{l=0}^{r}(-1)^l{{r}\brack{l}}_{q_i}
  x^{\pm}_{i,k_{\pi(1)}}\dots x^{\pm}_{i,k_{\pi(l)}}
  x^{\pm}_{j,s}\tss x^{\pm}_{i,k_{\pi(l+1)}}\dots x^{\pm}_{i,k_{\pi(r)}}=0,
\eeq
for any integers $k_1,\dots,k_r,s$. We will keep the indices $i\ne j$
fixed and denote the left hand side in \eqref{serrexx}
by $x^{\pm}(k_1,\dots,k_r;s)$. We will adapt an argument used in the Yangian context
by Levendorski~\cite{l:gd} to the quantum affine algebra case. We will prove
the relation $x^{\pm}(k_1,\dots,k_r;s)=0$ by using an induction argument
on the number of nonzero entries among the entries of
the tuples $(k_1,\dots,k_r;s)$. The induction base is the relation
$x^{\pm}(0,\dots,0;0)=0$. It holds because of the well-known equivalence between
the Drinfeld--Jimbo definition of the
the quantized enveloping algebra
$U_q({\spa}_{2n})$ and its $R$-matrix presentation; see \cite{rtf:ql}.
In our notation, the algebra $U_q({\spa}_{2n})$ can be identified with the subalgebra
of the quantum affine algebra $U_q(\wh{\spa}_{2n})$
obtained by restricting the range of the indices of the generators to the set $\{1,\dots,n\}$,
as defined in Section~\ref{subsec:isoDJD},
whereas its $R$-matrix presentation is the subalgebra of $U(R)$ generated by the zero mode
elements ${l}^{\pm}_{ij}[0]$ with $1\leqslant i,j\leqslant 2n$; see Section~\ref{sec:nd}.

The induction step will be based on the identities in the algebra $U(R)$ which are implied by
the previously verified relations,
\ben
\vp_{i}(u)\ts x^{\pm}_{j}(v)= \Big[\frac{u \tss q^{(\al_i,\al_j)\mp c/2} -v}
{u\tss q^{\mp c/2}-v\tss q^{(\al_i,\al_j)}}\Big]^{\pm 1}
x^{\pm}_{j}(v)\ts\vp_{i}(u)
\een
and
\ben
\psi_{i}(u)\ts x^{\pm}_{j}(v)= \Big[\frac{v \tss q^{(\al_i,\al_j)\mp c/2} -u}
{v\tss q^{\mp c/2}-u\tss q^{(\al_i,\al_j)}}\Big]^{\mp 1}
x^{\pm}_{j}(v)\ts\psi_{i}(u).
\een
By taking the coefficients of powers of $u$ and $v$ we derive that
\ben
[a_{i,k},x_{j,m}^{\pm}]=\pm\frac{[kA_{ij}]_{q_i}}{k}\ts q^{\mp |k|\tss c/2}\tss x_{j,k+m}^{\pm}.
\een
The rest of the argument is quite similar to \cite{l:gd}; it amounts
to calculating the commutators
\ben
\big[a_{i,k},x^{\pm}(k_1,\dots,k_p,0,\dots,0;s)\big]
\Fand \big[a_{j,k},x^{\pm}(k_1,\dots,k_p,0,\dots,0;s)\big]
\een
for a given $0\leqslant p<r$. By the induction hypothesis, both commutators
are zero which leads to
a system of two linear equations with
a nonzero determinant. Therefore, all elements
of the form $x^{\pm}(k_1,\dots,k_{p+1},0,\dots,0;s)$
are also equal to zero. This proves that
$x^{\pm}(k_1,\dots,k_r;s)=0$, as required.
\epf

Now recall the extended quantum affine algebra $U^{\ext}_q(\wh{\spa}_{2n})$
as introduced in Definition~\ref{def:eqaa}. By using Theorem~\ref{thm:relrbar}
and Proposition~\ref{prop:corrgauss} connecting the Gaussian generators of the
algebras
$U(R)$ and $U(\overline{R})$, we come to the following homomorphism theorem.

\bth\label{thm:homom}
The mapping
\begin{alignat}{2}\non
X^{+}_i(u)&\mapsto e^{+}_{ii+1}(u_{+})-e_{ii+1}^{-}(u_{-}),
\quad\qquad&&\text{for}\quad i=1,\dots,n,\\[0.4em]
\non
X^{-}_i(u)&\mapsto f^{+}_{i+1,i}(u_{-})-f^{-}_{i+1,i}(u_{+}),
\quad\qquad&&\text{for}\quad i=1,\dots,n,\\[0.4em]
h^{\pm}_j(u)&\mapsto h^{\pm}_j(u),\quad\qquad&&\text{for}\quad j=1,\dots,n+1.
\non
\end{alignat}
defines a homomorphism $DR: U^{\ext}_q(\wh{\spa}_{2n})\rightarrow U(R)$.
\eth

We will show in the next section that the homomorphism $DR$ provided by Theorem~\ref{thm:homom}
is in fact an isomorphism. To this end, we will construct an inverse map by employing
the universal $R$-matrix for the algebra $U_q(\wh{\spa}_{2n})$ in a way similar to
the type $A$ case; see \cite{fm:ha}.

\section{The universal $R$-matrix and inverse map}

We will use explicit formulas for the
universal $R$-matrix for the algebra $U_q(\wh{\g})$
obtained by Khoroshkin and Tolstoy~\cite{kt:ur}
and Damiani~\cite{d:rm,da:Un}.

Recall the Cartan matrix for $\g=\spa_{2n}$ defined in \eqref{cartan}
and consider the diagonal matrix $C=\diag[r_1,r_2,\dots,r_n]$ with
$r_i=(\al_i,\al_i)/2$. Then the matrix $B=[B_{ij}]:=CA$ is symmetric with $B_{ij}=(\al_i,\al_j)$.
We will use the notation $\tilde{B}=[\tilde{B}_{ij}]$ for the inverse matrix $B^{-1}$.
We will also need the $q$-deformed matrix $B(q)=[B_{ij}(q)]$ with $B_{ij}(q)=[B_{ij}]_q$
and its inverse $\tilde{B}(q)=[\tilde{B}_{ij}(q)]$; see \eqref{kq}.
It is clear that both matrices $\tilde{B}$ and $\tilde{B}(q)$ are symmetric.
The entries of $\tilde{B}$ are given by
\beql{Bij}
\tilde{B}_{ij}=\begin{cases}
n/4\qquad&\text{for}\quad i=j=n,\\
   j/2\qquad&\text{for}\quad i=n>j,\\
   j\qquad&\text{for}\quad n>i\geqslant j,
  \end{cases}
\eeq
while for any integer $k$ we have
\beql{Bijqk}
\tilde{B}_{ij}(q^{k})=
\begin{cases}\dfrac{[n]_{q^{k}}}{[2]_{q^k}[2]_{q^{k(n+1)}}}\qquad&\text{for}\quad i=j=n,\\[1.2em]
   \dfrac{[j]_{q^{k}}}{[2]_{q^{k(n+1)}}}\qquad&\text{for}\quad i=n>j,\\[1.2em]
   \dfrac{[2]_{q^{k(n+1-i)}}[j]_{q^{k}}}{[2]_{q^{k(n+1)}}}\qquad&\text{for}\quad n>i\geqslant j.
  \end{cases}
\eeq

With the presentation of the algebra $U_q(\wh{\g})$ used in Section~\ref{subsec:isoDJD},
consider the extended algebra $U_q(\tilde{\g})$ which is obtained
by adjoining an additional element $d$
with the relations
\ben
[d,k_i]=0, \qquad [d,E_{\al_i}]=\delta_{i,0}E_{\al_i},\qquad [d,F_{\al_i}]=-\delta_{i,0}F_{\al_i}.
\een
For a formal variable $u$ define an automorphism $D_u$ of
the algebra $U_q(\tilde{\g})\otimes \CC[u,u^{-1}]$
by
\ben
D_u(E_{\al_i})=u^{\delta_{i,0}}E_{\al_i},\qquad D_u(F_{\al_i})=u^{-\delta_{i,0}}F_{\al_i},
\qquad D_u(k_{i})=k_i,\qquad D_u(d)=d.
\een

The {\em universal $R$-matrix} is an element
$\mathfrak{R}\in U_q(\tilde{\g})\ts\hat{\otimes}\ts U_q(\tilde{\g})$
of a completed tensor product satisfying certain conditions; see~Drinfeld~\cite{d:ac}.
The conditions imply that this element is a solution of the Yang--Baxter equation
\ben
\mathfrak{R}_{12}\mathfrak{R}_{13}\mathfrak{R}_{23}=\mathfrak{R}_{23}\mathfrak{R}_{13}\mathfrak{R}_{12}.
\een
The explicit formula for $\mathfrak{R}$ uses the $\hbar$-adic settings so we will regard
the quantum affine algebra over $\CC[[\hbar]]$ and set $q=\exp(\hbar)\in \CC[[\hbar]]$.
Introduce elements $h_1,\dots,h_n$ by
setting $k_i=\exp(\hbar h_i)$. According to \cite{da:Un},
the universal $R$-matrix admits a triangular decomposition
\beql{dam}
\mathfrak{R}=\mathfrak{R}^{>0}\ts\mathfrak{R}^{0}\ts\mathfrak{R}^{<0}\ts\mathcal{K},
\eeq
where
\ben
\bal
\mathfrak{R}^{>0}&=\prod_{\alpha\in \Delta_+}\prod_{k\geqslant 0}
 \exp_{q_{\alpha}}\big((q_{\alpha}^{-1}-q_{\alpha})E_{\alpha+k\delta}\otimes F_{\alpha+k\delta}\big),\\
\mathfrak{R}^{<0}&=\prod_{\alpha\in \Delta_+}\prod_{k> 0}
 \exp_{q_{\alpha}}\big((q_{\alpha}^{-1}-q_{\alpha})E_{-\alpha+k\delta}\otimes F_{-\alpha+k\delta}\big),
\eal
\een
and
\ben
\mathcal{K}=T\ts q^{-(c\otimes d+d\otimes c)},\qquad
T=\exp(-\hbar \tilde{B}_{ij}\tss h_i\otimes h_j)
\een
where $q_{\alpha}=q^{(\alpha,\alpha)/2}$.

We will work with the parameter-dependent $R$-matrix defined by
\ben
\Rc(u)=(D_u\otimes \id)\ts\mathfrak{R}\ts q^{c\otimes d+d\otimes c}.
\een
It satisfies the Yang--Baxter equation in the form
\beql{MYBE}
\Rc_{12}(u)\Rc_{13}(uvq^{-c_2})\Rc_{23}(v)
=\Rc_{23}(v)\Rc_{13}(uvq^{c_2})\Rc_{12}(u)
\eeq
where $c_2=1\otimes c\otimes 1$; cf. \cite{fri:qa}.

A straightforward calculation verifies the following formulas for the vector
representation of the quantum affine algebra. As before, we denote by $e_{ij}\in\End\CC^{2n}$
the standard matrix units.

\bpr\label{prop:FFRep}
The mappings $q^{\pm c/2}\mapsto 1$,
\ben
\bal
    x^{+}_{ik} &\mapsto -q^{-ik}e_{i+1,i}+q^{-(2n+2-i)k}e_{i',(i+1)'}, \\
    x^{-}_{ik} &\mapsto -q^{-ik}e_{i,i+1}+q^{-(2n+2-i)k}e_{(i+1)',i'},\\
    a_{ik} &\mapsto \frac{[k]_{q_i}}{k}\big(q^{-ik}(q^{-k}e_{i+1,i+1}-q^{k}e_{ii})
    +q^{-(2n+2-i)k}(q^{-k}e_{i'i'}-q^ke_{(i+1)'(i+1)'})
    \big)\\
    k_i &\mapsto q(e_{i+1,i+1}+e_{i',i'})+q^{-1}(e_{ii}+e_{(i+1)',(i+1)'})
    +\sum\limits_{j\neq i,i+1,i',(i+1)'}e_{jj},
\eal
\een
    for $i=1,\dots,n-1$, and
\ben
\bal    x^{+}_{nk} &\mapsto -q^{-(n+1)k}e_{n+1,n}, \\
    x^{-}_{nk} &\mapsto -q^{-(n+1)k}e_{n,n+1},\\
    a_{nk} &\mapsto \frac{[k]_{q_n}}{k}\big(q^{-(n+1)k}(q^{-2k}e_{n+1,n+1}-q^{2k}e_{nn})
    \big)\\
    k_n &\mapsto q^2e_{n+1,n+1}+q^{-2}e_{n,n}+\sum\limits_{j\neq n,n+1}e_{jj},
\eal
\een
define a representation $\pi_V: U_q(\wh{\spa}_{2n})\to\End V$
of the algebra $U_q(\wh{\spa}_{2n})$ on the vector space $V=\mathbb{C}^{2n}$.
\qed
\epr

It follows from the results of \cite{fri:qa} that the $R$-matrix defined in \eqref{ru}
coincides with the image of the universal $R$-matrix:
\ben
R(u)=(\pi_{V}\otimes \pi_{V})\ts\Rc(u).
\een
Introduce the $L$-{\em operators} in $U_q(\wh{\spa}_{2n})$ by the formulas
\ben
\bal
  \tilde{L}^{+}(u) &= (\id\otimes \pi_{V})\ts \Rc_{21}(uq^{c/2}),\\
  \tilde{L}^{-}(u) &= (\id\otimes \pi_{V})\ts \Rc_{12}(u^{-1}q^{-c/2})^{-1}.
\eal
\een
Recall the series $z^{\pm}(u)$ defined in \eqref{zpm}. Their coefficients
are central in the algebra $U^{\ext}_q(\wh{\spa}_{2n})$; see Proposition~\ref{prop:zu}.
Therefore, the Yang--Baxter equation \eqref{MYBE}
implies the relations for the $L$-operators:
\ben
\bal
R(u/v)L^{\pm}_1(u)L^{\pm}_2(v) &= L^{\pm}_2(v)L^{\pm}_1(u)R(u/v), \\[0.4em]
  R(u_{+}/v_{-})L^{\pm}_1(u)L^{\mp}_2(v) &= L^{\mp}_2(v)L^{\pm}_1(u)R(u_{-}/v_{+}),
\eal
\een
where we set
\begin{align}
\label{lplus}
L^{+}(u) &= \tilde{L}^{+}(u)\prod_{m=0}^{\infty}
z^{+}(u\tss\xi^{-2m-1})\tss z^{+}(u\tss\xi^{-2m-2})^{-1},\\
\label{lminus}
L^{-}(u) &= \tilde{L}^{-}(u)\prod_{m=0}^{\infty}
z^{-}(u\tss\xi^{-2m-1})\tss z^{-}(u\tss\xi^{-2m-2})^{-1}.
\end{align}
Note that although these formulas for the entries of the matrices $L^{\pm}(u)$
involve a completion of the center of the algebra $U^{\ext}_q(\wh{\spa}_{2n})$,
it will turn out that the coefficients of the series in $u^{\pm1}$ actually
belong to $U^{\ext}_q(\wh{\spa}_{2n})$; see the proofs of
Propositions~\ref{prop:Hp} and \ref{prop:Lm} below. Thus, we may conclude that the mapping
\beql{inversemap}
RD:L^{\pm}(u)\mapsto L^{\pm}(u)
\eeq
defines a homomorphism $RD$ from the algebra $U(R)$
to a completed algebra $U^{\ext}_q(\wh{\spa}_{2n})$,
where we use the same notation for the corresponding elements of the algebras.

Returning to the universal $R$-matrix, observe that formula \eqref{dam} implies
the corresponding decomposition of the matrix $\Rc(u)$:
\beql{rdec}
\Rc(u)=\Rc^{>0}(u)\Rc^{0}(u)\Rc^{<0}(u),
\eeq
where
\ben
\bal
  \Rc^{>0}(u) &= \prod_{\alpha\in \Delta_+}\prod_{k\geqslant 0}
 \exp_{q_{\alpha}}\big((q_{\alpha}^{-1}-q_{\alpha})
 u^{k}E_{\alpha+k\delta}\otimes F_{\alpha+k\delta}\big), \\
 \Rc^{<0}(u) &= T^{-1}\prod_{\alpha\in \Delta_+}\prod_{k> 0}
 \exp_{q_{\alpha}}\big((q_{\alpha}^{-1}-q_{\alpha})
 u^{k}E_{-\alpha+k\delta}\otimes F_{-\alpha+k\delta}\big)\ts T
\eal
\een
and
\ben
  \Rc^{0}(u)=\exp\Big(
  \sum\limits_{k>0}\sum\limits_{i,j=1}^{n}\frac{(q_i^{-1}-q_i)(q_j^{-1}-q_j)}{q^{-1}-q}
  \frac{k}{[k]_q}\tilde{B}_{ij}(q^k)u^k q^{kc/2}a_{i,k}\otimes a_{j,-k}q^{-kc/2}\Big)\ts T.
\een
By using the vector representation $\pi_{V}$ defined in Proposition~\ref{prop:FFRep},
introduce the matrices $F^{+}(u)$, $E^{+}(u)$ and $H^{+}(u)$ by setting
\ben
\bal
F^{+}(u)&=(\id\otimes \pi_{V})\ts \Rc^{>0}_{21}(uq^{c/2})\\
&=(\id\otimes \pi_{V})\prod_{\alpha\in \Delta_+}\prod_{k\geqslant 0}
 \exp_{q_{\alpha}}\big((q_{\alpha}^{-1}-q_{\alpha})u^{k}q^{kc/2}
 F_{\alpha+k\delta}\otimes E_{\alpha+k\delta}\big),
\eal
\een
\ben
\bal
E^{+}(u)&=(\id\otimes \pi_{V})\ts \Rc^{<0}_{21}(uq^{c/2})\\
&=(\id\otimes \pi_{V})\Big(T_{21}^{-1}\prod_{\alpha\in \Delta_+}\prod_{k> 0}
 \exp_{q_{\alpha}}\big((q_{\alpha}^{-1}-q_{\alpha})u^{k}q^{kc/2}
 F_{-\alpha+k\delta}\otimes E_{-\alpha+k\delta}\big)\ts
 T_{21}\Big)
\eal
\een
and
\begin{multline}
\non
H^{+}(u)
=(\id\otimes \pi_{V})\\
\Big(\exp\Big(\sum\limits_{k>0}\sum\limits _{i,j=1}^{n}
\frac{(q_i^{-1}-q_i)(q_j^{-1}-q_j)}{q^{-1}-q}\frac{k}{[k]_q}
\tilde{B}_{ij}(q^k)u^kq^{kc/2} a_{j,-k}q^{-kc/2}\otimes q^{kc/2}a_{i,k} \Big)T_{21}\Big)\\
{}\times \prod_{m=0}^{\infty} z^{+}(u\xi^{-2m-1})z^{+}(u\xi^{-2m-2})^{-1}.
\end{multline}

The decomposition \eqref{rdec} implies the corresponding decomposition
for the matrix $L^{+}(u)$:
\ben
L^{+}(u)=F^{+}(u)H^{+}(u)E^{+}(u).
\een

Recall the Drinfeld generators $x^{\pm}_{i,k}$ of the algebra $U_q(\wh{\spa}_{2n})$,
as defined in the Introduction, and combine them into the formal series
\ben
\bal
x_i^{-}(u)^{\geqslant 0}&=\sum\limits _{k\geqslant 0}x^{-}_{i,-k}u^{k},
\qquad\quad x^{+}_{i}(u)^{> 0}=\sum\limits _{k> 0}x^{+}_{i,-k}u^{k},\\
x_i^{-}(u)^{< 0}&=\sum\limits _{k> 0}x^{-}_{i,k}u^{-k},
\qquad\quad x^{+}_{i}(u)^{\leqslant 0}=\sum\limits _{k\geqslant 0}x^{+}_{i,k}u^{-k}.
\eal
\een
Furthermore, set
\ben
\bal
f_i^{+}(u)&=(q_i -q_i^{-1})x_i^{-}(u_+q^{-i})^{\geqslant 0},\qquad\quad
e_i^{+}(u)=(q_i -q_i^{-1})x^{+}_{i}(u_{-}q^{-i})^{> 0},\\
f_{i}^{-}(u)&=(q_i^{-1}-q_i)x^{-}_{i}(u_{-}q^{-i})^{<0},\qquad\quad
e_i^{-}(u)=(q_i^{-1}-q_i)x^{+}_{i}(u_{+}q^{-i})^{\leqslant 0}
\eal
\een
for $i=1,\dots, n-1$, and
\ben
\bal
f_n^{+}(u)&=(q_n -q_n^{-1})x_n^{-}(u_+q^{-(n+1)})^{\geqslant 0},\qquad\quad
e_n^{+}(u)=(q_n -q_n^{-1})x_n^{+}(u_{-}q^{-(n+1)})^{> 0},\\
f_n^{-}(u)&=(q_n^{-1}-q_n)x^{-}_{n}(u_{-}q^{-(n+1)})^{<0},\qquad\quad
e_n^{-}(u)=(q_n^{-1}-q_n)x_n^{+}(u_{+}q^{-(n+1)})^{\leqslant 0}.
\eal
\een

\bpr\label{prop:Fp}
The matrix $F^{+}(u)$ is lower unitriangular and has the form
\ben
F^{+}(u)=
\begin{bmatrix}
  1 & &  & &  &  &  &  \\
  f_{1}^{+}(u) & \qquad 1 &  &  & & &\bigcirc  &  \\
   &  \qquad \ddots& & \ddots&  &  &  & \\
   &  &  & f^{+}_n(u)& 1&  &  &  \\[0.5em]
   &  &  &  & -f^{+}_{n-1}(u\xi q^{2(n-1)}) & 1 & &  \\
   &  \qquad \bigstar & &  &  &  \ddots& \ddots &  \\[0.2em]
   &  &  &  &  &  &  -f_{1}^{+}(u \xi q^{2})& 1
\end{bmatrix}.
\een
\epr

\bpf
By the construction of the root vectors $E_{\alpha+k\delta}$ and the formulas for the representation $\pi_V$,
we only need to evaluate the image of
the product
\ben
\prod_{k\geqslant 0}
 \exp_{q_i}\big((q_i^{-1}-q_i)(uq^{c/2})^{k} F_{\alpha_i+k\delta}\otimes E_{\alpha_i+k\delta}\big)
\een
for simple roots $\al_i$ with $i=1,\dots,n$.
Due to the isomorphism of Sec.~\ref{subsec:isoDJD}, we can rewrite
it in terms of Drinfeld generators as
\ben
\prod_{k\geqslant 0}
 \exp_{q_i}\big((q_i^{-1}-q_i)(uq^{c/2})^{k}x^{-}_{i,-k}\otimes x^{+}_{i,k}\big).
\een
Suppose first that $i\leqslant n-1$. Using the formulas for the action
of the generators $x^{+}_{i,k}$ from Proposition~\ref{prop:FFRep}, we get
\begin{align}
\label{itpiv}
&(\id\otimes \pi_V)\prod_{k\geqslant 0}
 \exp_{q_i}\big((q_i^{-1}-q_i)(uq^{c/2})^{k}x^{-}_{i,-k}\otimes x^{+}_{i,k}\big)\\
& =
 \prod_{k\geqslant 0}\exp_{q_i}\big(-(q_i^{-1}-q_i)(u_{+}q^{-i})^{k}x^{-}_{i,-k}\otimes e_{i+1,i}
 +(q_i^{-1}-q_i)(u_{+}q^{-(2n+2-i)})^{k}x^{-}_{i,-k}\otimes e_{i',(i+1)'}\big).
 \non
\end{align}
Expanding the $q$-exponent, we can write this expression in the form
\ben
\bal
 1&-(q_i^{-1}-q_i)\sum\limits_{k\geqslant 0}x^{-}_{i,-k}(u_{+}q^{-i})^{k}\otimes e_{i+1,i}
 +(q_i^{-1}-q_i)\sum\limits_{k\geqslant 0}x^{-}_{i,-k}(u_{+}q^{-(2n+2-i)})^{k}\otimes e_{i',(i+1)'}\\
 {}=1&+(q_i -q_i^{-1})x^{-}_{i}(u_{+}q^{-i})^{\geqslant 0}\otimes e_{i+1,i}
 -(q_i -q_i^{-1})x^{-}_{i}(u_{+}q^{-(2n+2-i)})^{\geqslant 0}\otimes e_{i',(i+1)'}
\eal
\een
which coincides with $1+f^{+}_{i}(u)\otimes e_{i+1,i}-f^{+}_{i}(u\xi q^{2i})\otimes e_{i',(i+1)'}$,
as required. A similar calculation shows that expression \eqref{itpiv} with $i=n$
simplifies to $1+f_n^{+}(u)\otimes e_{n+1,n}$.
\epf

As in Sec.~\ref{subsec:eqaa}, we will assume that the algebra $U_q(\wh{\spa}_{2n})$
is extended by adjoining the square roots $k_n^{\pm1/2}$.

\ble\label{lem:K} The image $(\id\ot \pi_V)\big(T_{21}\big)$ is the diagonal matrix
\begin{multline}
\diag\Big[k_1\dots k_{n-1}\tss k_n^{1/2},\ \  k_2\dots k_{n-1}\tss k_n^{1/2},
\ \ \dots,\ \ k_n^{1/2},\\ k_n^{-1/2},\ \
k_{n-1}^{-1}\tss k_n^{-1/2},\ \ \dots,\ \ k_1^{-1}\dots k_{n-1}^{-1}\tss k_n^{-1/2}\Big].
\non
\end{multline}
\ele

\bpf
By definition, we have
\ben
\bal
(\id\ot \pi_V)\big(T_{21}\big)&=
\exp\big(-\hbar \sum_{b=1}^n \sum_{a=1}^n \tilde{B}_{ab}h_b\otimes\pi_{V}(h_a)\big)\\
&=\exp\big(-\hbar\sum_{b=1}^n \sum_{a=1}^{n-1}
\tilde{B}_{ab}h_b\otimes (e_{a+1,a+1}-e_{a,a}-e_{(a+1)',(a+1)'}+e_{a',a'})\big)\\
&\qquad\qquad\qquad\qquad{}\times \exp\big(-2\hbar \sum_{b=1}^n
\tilde{B}_{nb}h_b\otimes (e_{n+1,n+1}-e_{n,n})\big)
\eal
\een
which equals
\ben
\bal
&\exp\Big(\hbar\sum_{b=1}^n \sum_{a=2}^{n-1} (\tilde{B}_{ab}-\tilde{B}_{a-1,b})h_j\otimes e_{a,a}+
\hbar\sum_{b=1}^n \tilde{B}_{1b}h_b\otimes e_{1,1}\\
&-\hbar\sum_{b=1}^n \sum_{a=2}^{n-1} (\tilde{B}_{ab}-\tilde{B}_{a-1,b})h_b\otimes e_{a',a'}+
\hbar\sum_{b=1}^n \tilde{B}_{1b}h_b\otimes e_{1',1'}\\
&+\hbar\sum_{b=1}^n (2\tilde{B}_{nb}-\tilde{B}_{n-1,b})h_j\otimes e_{n,n}
-\hbar\sum_{b=1}^n (2\tilde{B}_{nb}-\tilde{B}_{n-1,b})h_j\otimes e_{n',n'}\Big).
\eal
\een
The claim now follows by applying formula \eqref{Bij}
for the entries $\tilde{B}_{ij}$. For instance,
the $(1,1)$-entry of the diagonal matrix is found by
\ben
\exp\Big(\hbar\sum_{b=1}^n \tilde{B}_{1b}h_b\Big)
=\exp\Big(\hbar\sum_{b=1}^{n-1}h_b+\hbar\frac{1}{2}h_n)\Big)
=\prod_{b=1}^{n-1}k_b\ts k_n^{1/2},
\een
and the remaining entries are obtained by the same calculation.
\epf

\bpr\label{prop:Ep}
The matrix $E^{+}(u)$ is upper unitriangular and has the form
\ben
E^{+}(u)=
\begin{bmatrix}
  1 & e_1^{+}(u) &  &  &  &  &  &  \\
    & 1 & e_2^{+}(u) &  &  &  &  \bigstar& \\
   &  & \ddots & \ddots &  &  &  &  \\
   &  &  & 1 & e_n^{+}(u) &  &  &  \\
   & &  &  & \ddots & \ddots &  &  \\
   &  & &  &  & 1 &  -e_2^{+}(u\xi q^4)&  \\[0.5em]
   & \bigcirc &  &  &  &  &  1& -e_1^{+}(u\xi q^2) \\
  &  &  &  &  &  &  & 1
\end{bmatrix}.
\een
\epr

\bpf
It is sufficient to evaluate the image of the product
\ben
T_{21}^{-1}\prod_{k> 0}
 \exp_{q_i}\big((q_i^{-1}-q_i)u^{k}q^{kc/2}F_{-\alpha_i+k\delta}\otimes E_{-\alpha_i+k\delta}\big)
\ts T_{21}
\een
with respect to $\id\otimes \pi_{V}$
for simple roots $\al_i$ with $i=1,\dots,n$.
Using the isomorphism of Sec.~\ref{subsec:isoDJD}, we can rewrite the internal product
in terms of Drinfeld generators as
\ben
\prod_{k> 0}
 \exp_{q_i}\big((q_i^{-1}-q_i)(uq^{c/2})^{k}q^{-kc}x^{+}_{i,-k}k_i\otimes q^{kc}k_i^{-1}x^{-}_{i,k}\big).
\een
The remaining calculation is performed in the same way as in the proof
of Proposition~\ref{prop:Fp} with the use of Proposition~\ref{prop:FFRep},
Lemma~\ref{lem:K} and the relations
$k_ix_{j,k}^{\pm}k_i^{-1}=q_i^{\pm A_{ij}}x_{j,k}^{\pm}$.
\epf

In the next proposition we use the series $z^{\pm}(u)$
introduced in \eqref{zpm}. Their coefficients belong to the center of the algebra
$U^{\ext}_q(\wh{\spa}_{2n})$; see Proposition~\ref{prop:zu}.
For a nonnegative integer $m$ with
$m< n$ we will denote by $z^{\pm\ts[n-m]}(u)$ the respective series
for the subalgebra of $U^{\ext}_q(\wh{\spa}_{2n})$, whose generators
are all elements
$X^{\pm}_{i,k}$, $h^{\pm}_{j,k}$ and $q^{c/2}$ such that $i,j\geqslant m+1$;
see Definition~\ref{def:eqaa}. We also denote by $\xi^{[n-m]}$ the parameter
$\xi$ for this subalgebra so that $\xi^{[n-m]}=q^{-2n+2m-2}$.

\bpr\label{prop:Hp}
The matrix $H^{+}(u)$ is diagonal and has the form
\ben
H^{+}(u)=\diag\tss\big[h_1^{+}(u),\dots,h^{+}_{n}(u),
z^{+\tss [1]}(u)\tss h_{n}^{+}(u\xi^{[1]})^{-1},\dots,
z^{+\tss [n]}(u)\tss h_1^{+}(u\xi^{[n]})^{-1}\big].
\een
\epr

\bpf
 By definition,
\ben
\bal
H^{+}(u)&=\exp\Big(\sum\limits_{k>0}\sum\limits _{i,j=1}^{n}\frac{(q_i^{-1}-q_i)
(q_j^{-1}-q_j)}{q^{-1}-q}\frac{k}{[k]_q}
\tilde{B}_{ij}(q^k)u^kq^{kc/2} a_{j,-k}q^{-kc/2}\otimes\pi_V(q^{kc/2}a_{i,k}) \Big)\\
&\times (\id\ot \pi_V)(T_{21})\prod_{m=0}^{\infty} z^{+}(u\xi^{-2m-1})z^{+}(u\xi^{-2m-2})^{-1}.
\eal
\een
Using the formulas for $\pi_V(a_{i,k})$ from Proposition~\ref{prop:FFRep},
we can write the first factor as the exponent of the expression
\ben
\bal
\sum_{k>0}\sum_{j=1}^{n}\sum_{i=1}^{n-1}&(q_j-q_j^{-1})\tilde{B}_{ij}(q^{k})u^{k}a_{j,-k}\\
&\otimes
\big(q^{-(i-1)k}e_{i,i}-q^{-(i+1)k}e_{i+1,i+1}
-\xi^{k}q^{(i-1)k}e_{i',i'}+\xi^kq^{(i+1)k}e_{(i+1)',(i+1)'}\big)\\
{}+\sum_{k>0}\sum_{j=1}^{n}&(q_j-q_j^{-1})\tilde{B}_{nj}(q^{k})(q^k+q^{-k})u^{k}a_{j,-k}\otimes
\big(q^{-(n-1)k}e_{n,n}-q^{-(n+3)k}e_{n+1,n+1}\big).
\eal
\een
Consider the $(1,1)$-entry (the coefficient of $e_{1,1}$) in the first factor
in the formula for $H^{+}(u)$.
Using formula \eqref{Bijqk} for $\tilde{B}_{1,j}(q^{k})$ we get
\begin{multline}
\non
\exp\Big(\sum_{k>0}\sum_{j=1}^{n}(q_j-q_j^{-1})\tilde{B}_{1,j}(q^{k})u^{k}a_{j,-k}\Big)\\
{}=\exp\Big(\sum_{k>0}\sum_{j=1}^{n-1}(q-q^{-1})
\frac{q^{jk}+\xi^{-k}q^{-jk}}{1+\xi^{-k}}u^{k}a_{j,-k}\Big)\ts
\exp\Big(\sum_{k>0}(q_n-q^{-1}_{n})\frac{q^{(n+1)k}}{1+\xi^{-k}}u^{k}a_{n,-k}\Big).
\end{multline}
By expanding the fractions into power series, we can write this expression as
\begin{multline}
\non
\exp\Big(\sum_{k>0}\sum_{j=1}^{n-1}\sum_{m=0}^{\infty}(q-q^{-1})(-1)^m
\big(\xi^{-mk}q^{jk}+\xi^{-mk-k}q^{-jk}\big)u^{k}a_{j,-k}\Big)\\
{}\times \exp\Big(\sum_{k>0}
\sum_{m=0}^{\infty}(q_n-q^{-1}_n)(-1)^m\xi^{-mk}q^{-(n+1)k}u^{k}a_{n,-k}\Big).
\end{multline}
Using the definition \eqref{phiiu}
of the series $\varphi_{i}(u)$ and setting $\tilde{\varphi}_{i}(u)=k_i\tss\varphi_{i}(u)$,
we can bring the expression to the form
\begin{multline}
\non
\prod_{m=0}^{\infty} \prod_{j=1}^{n-1}\tilde{\varphi}_{j}(u\xi^{-2m}q^j)^{-1}
\tilde{\varphi}_{j}(u\xi^{-2m-1}q^j)
\tilde{\varphi}_{j}(u\xi^{-2m-1}q^{-j})^{-1}\tilde{\varphi}_{j}(u\xi^{-2m-2}q^{-j})\\
{}\times \prod_{m=0}^{\infty} \tilde{\varphi}_{n}(u\xi^{-2m}q^{n+1})^{-1}
\tilde{\varphi}_{n}(u\xi^{-2m-1}q^{n+1}).
\end{multline}
Setting $\tilde{h}^{+}_{i}(u)=t_i^{-1}\tss {h}^{+}_{i}(u)$ with
$t_i=h_{i,0}^{+}$ and
applying Proposition~\ref{prop:embed}, we can write this as
\begin{multline}
\non
\prod_{m=0}^{\infty} \prod_{j=1}^{n-1}\tilde{h}^{+}_{j}(u\xi^{-2m}q^{2j})
\tilde{h}^{+}_{j}(u\xi^{-2m-1}q^{2j})^{-1}
\times\prod_{m=0}^{\infty} \prod_{j=1}^{n}\tilde{h}^{+}_{j}(u\xi^{-2m}q^{2j-2})^{-1}
\tilde{h}^{+}_{j}(u\xi^{-2m-1}q^{2j-2})\\
\times\prod_{m=0}^{\infty}\tilde{h}^{+}_{n+1}(u\xi^{-2m-1})^{-1}
\tilde{h}^{+}_{n+1}(u\xi^{-2m-2})\times \tilde{h}^{+}_1(u).
\end{multline}
Now use definition \eqref{zpm} of the series $z^{\pm}(u)$ to conclude that
\ben
\exp\Big(\sum_{k>0}\sum_{j=1}^{n}(q_j-q_j^{-1})\tilde{B}_{1,j}(q^{k})u^{k}a_{j,-k}\Big)
=\prod_{m=0}^{\infty}z^{+}(u\xi^{-2m-1})^{-1}z^{+}(u\xi^{-2m-2})\times \tilde{h}^{+}_1(u).
\een
In particular, the coefficients of powers of $u$ in the infinite product
which occurs in \eqref{lplus} belong to the algebra $U^{\ext}_q(\wh{\spa}_{2n})$.

Furthermore, Lemma~\ref{lem:K} implies that
the $(1,1)$-entry of the matrix
$(\id\ot \pi_V)(T_{21})$ equals
$\prod_{j=1}^{n-1}k_j\times k_n^{1/2}=t_1$.
This proves that the $(1,1)$-entry of the matrix $H^{+}(u)$ is $h_1^{+}(u)$.

It is clear that the matrix $H^{+}(u)$ is diagonal, and we perform quite similar calculations
to evaluate the $(i,i)$-entries for $i=2,\dots,2n$. For instance, if
$i=2,\dots,n-1$ then
formula \eqref{Bijqk} for $\tilde{B}_{ij}(q^{k})$ implies that the exponent
\ben
\exp\Big(\sum_{k>0}\sum_{j=1}^{n}(q_j-q_j^{-1})\big(q^{-(i-1)k}\tilde{B}_{ij}(q^{k})-
q^{-ik}\tilde{B}_{i-1,j}(q^{k})\big)u^{k}a_{j,-k}\Big)
\een
can be written in terms of the series $\tilde{\varphi}_i(u)$ as
\ben
\bal
\prod_{m=0}^{\infty}\prod_{j=1}^{i-1}&\tilde{\varphi}_j(u\xi^{-2m}q^{j})^{-1}
\tilde{\varphi}_j(u\xi^{-2m}q^{-j})
\tilde{\varphi}_j(u\xi^{-2m-1}q^{j})\tilde{\varphi}_j(u\xi^{-2m-1}q^{-j})^{-1}\\
&\times \prod_{m=0}^{\infty}\prod_{j=i}^{n-1}\tilde{\varphi}_j(u\xi^{-2m}q^{j})^{-1}
\tilde{\varphi}_j(u\xi^{-2m-1}q^{-j})^{-1}
\tilde{\varphi}_j(u\xi^{-2m-1}q^{j})\tilde{\varphi}_j(u\xi^{-2m-2}q^{-j})\\
&\times \prod_{m=0}^{\infty}\tilde{\varphi}_n(u\xi^{-2m}q^{n+1})^{-1}
\tilde{\varphi}_n(u\xi^{-2m-1}q^{n+1}).
\eal
\een
Using the relations $\tilde{\varphi}_j(u)=\tilde{h}^{+}_j(uq^j)^{-1}\tilde{h}^{+}_{j+1}(uq^j)^{-1}$
and $\tilde{\varphi}_n(u)=\tilde{h}^{+}_n(uq^{n+1})^{-1}\tilde{h}^{+}_{n+1}(uq^{n+1})^{-1}$, we
can write this expression in the form
\ben
\bal
&\prod_{m=0}^{\infty}\Big(\prod_{j=1}^{n-1}\tilde{h}^{+}_j(u\xi^{-2m}q^{2j})
\prod_{j=1}^{n}\tilde{h}^{+}_{j}(u\xi^{-2m}q^{2j-2})^{-1}\times\tilde{h}^{+}_n(u\xi^{-2m-1})^{-1}\Big)\\
&\times
\prod_{m=0}^{\infty}
\Big(\prod_{j=1}^{n-1}\tilde{h}^{+}_j(u\xi^{-2m-1}q^{2j})^{-1}
\prod_{j=1}^{n}\tilde{h}^{+}_{j}(u\xi^{-2m-1}q^{2j-2})\times \tilde{h}^{+}_{n+1}(u\xi^{-2m-2})\Big)
\times \tilde{h}^{+}_i(u)
\eal
\een
which equals
\ben
\prod_{m=0}^{\infty}z^{+}(u\xi^{-2m-1})^{-1}z^{+}(u\xi^{-2m-2})\times\tilde{h}^{+}_i(u).
\een
Furthermore, by Lemma \ref{lem:K}, the coefficient of $e_{i,i}$ in
the image $(\id\ot \pi_V)(T_{21})$ coincides with
$\prod_{j=i}^{n-1}k_i k_n^{1/2}=t_i$. This shows that for $i=2,\dots,n-1$
the $(i,i)$-entry of the matrix $H^{+}_i(u)$ is $h^{+}_i(u)$.
We omit the calculations in the remaining cases which are quite similar.
\epf

\bre\label{rem:prop2.3}
The calculations used in the proof of Proposition~\ref{prop:Hp} also provide
an alternative argument to verify that the map $\varrho$ introduced
in the proof of Proposition~\ref{prop:embed} is a homomorphism.
Namely, the map $L^{\pm}(u)\mapsto \tilde{L}^{\pm}(u)$ defines a homomorphism
$U(R)\to U_q(\wh{\spa}_{2n})$. By the uniqueness of the Gauss decomposition
this implies that the images of the entries of the diagonal matrix $H^{\pm}(u)$
are the respective entries of the diagonal matrix $\tilde H^{\pm}(u)$
defined by the Gauss decomposition of $\tilde{L}^{\pm}(u)$. However,
the matrix $\tilde H^{\pm}(u)$ has the form
\ben
\tilde H^{\pm}(u)=\diag\tss\big[\al_1^{\pm}(u),\dots,\al^{\pm}_{2n}(u)\big],
\een
where the series $\al_i^{\pm}(u)$ were introduced in the proof of
Proposition~\ref{prop:embed}. On the other hand, Theorems~\ref{thm:relrbar} and \ref{thm:homom}
imply that the series $h^{\pm}_i(u)$ satisfy the relations described in Definition~\ref{def:eqaa}
and so do their respective $\varrho$-images $\al^{\pm}_i(u)$.
\qed
\ere

Now turn to the matrix $L^{-}(u)$. By definition \eqref{lminus} of $L^{-}(u)$,
we have
\ben
L^{-}(u)=F^{-}(u)H^{-}(u)E^{-}(u),
\een
where
\ben
E^{-}(u)=(\id\otimes \pi_V)\Rc^{>0}(u_+^{-1})^{-1},
\een
\ben
F^{-}(u)=(\id\otimes \pi_V)\Rc^{<0}(u_+^{-1})^{-1},
\een
and
\ben
H^{-}(u)=(\id\otimes \pi_V)(\Rc^{0}(u_+^{-1}))^{-1}
\prod_{m=0}^{\infty}z^{-}(u\tss\xi^{2m-1})^{-1}z^{-}(u\tss\xi^{2m-2}).
\een

\bpr\label{prop:Lm}
The matrix $E^{-}(u)$ is upper unitriangular and has the form
\ben
E^{-}(u)=
\begin{bmatrix}
  1 & e_1^{-}(u) &  &  &  &  &  &  \\
    & 1 & e_2^{-}(u) &  &  & & \bigstar &  \\
   & & \ddots & \ddots &  &  &  &  \\
   &  &  & 1 & e_n^{-}(u) &  &  &  \\
   & &  &  & \ddots & \ddots &  &  \\
   &  & &  &  & 1 &  -e_2^{-}(u\xi q^4)&  \\
   & \bigcirc &  &  &  &  &  1& -e_1^{-}(u\xi q^2) \\
  &  &  &  &  &  &  & 1
\end{bmatrix}.
\een
The matrix $F^{-}(u)$ is lower unitriangular and has the form
\ben
F^{-}(u)=
\begin{bmatrix}
  1 &  &  &  &  &  &  &  \\
  f_{1}^{-}(u) & 1 &  &  & & \bigcirc &  &  \\
   &  \ddots& \ddots & &  &  &  & \\
   &  &  f_{n-1}^{-}(u)& 1&  &  &  &  \\
   &  &  & f^{-}_n(u)& 1&  &  &  \\[0.5em]
   &  &  &  & -f^{-}_{n-1}(u\xi q^{2(n-1)}) & 1 & &  \\
   &  &\bigstar  &  &  &  \ddots& \ddots &  \\
   &  &  &  &  &  &  -f_{1}^{-}(u \xi q^{2})& 1
\end{bmatrix}.
\een
The matrix $H^{-}(u)$ is diagonal and has the form
\ben
H^{-}(u)=\diag\tss\big[h_1^{-}(u),\dots,h^{-}_{n}(u),
z^{-\tss [1]}(u)\tss h_{n}^{-}(u\xi^{[1]})^{-1},\dots,
z^{-\tss [n]}(u)\tss h_1^{-}(u\xi^{[n]})^{-1}\big].
\een
\epr

\bpf
As a first step, use the same arguments as in the proofs of Propositions~\ref{prop:Fp},
\ref{prop:Ep} and \ref{prop:Hp} to evaluate the matrices
\ben
E^{-}(u)^{-1}=(\id\otimes \pi_V)\Rc^{+}(u_+^{-1}),\qquad
F^{-}(u)^{-1}=(\id\otimes \pi_V)\Rc^{-}(u_+^{-1}),
\een
and
\ben
H^{-}(u)^{-1}=(\id\otimes \pi_V)(\Rc^{0}(u_+^{-1}))
\prod_{m=0}^{\infty}z^{-}(u\tss\xi^{2m-1})\ts z^{-}(u\tss\xi^{2m-2})^{-1}.
\een
The required expressions are then obtained by inverting the respective matrices.
\epf

Taking into account Propositions~\ref{prop:Fp}, \ref{prop:Ep}, \ref{prop:Hp}
and \ref{prop:Lm} we arrive at the following result.
\bco\label{cor:isom}
The homomorphism
\ben
RD: U(R)\rightarrow U^{\ext}_q(\wh{\spa}_{2n})
\een
defined in \eqref{inversemap}
is the inverse map to the homomorphism $DR$ defined in Theorem~\ref{thm:homom}.
Hence the algebra $U(R)$ is isomorphic to $U^{\ext}_q(\wh{\spa}_{2n})$.
\qed
\eco

Corollary~\ref{cor:isom} together with the results of Secs~\ref{subsec:eqaa} and \ref{subsec:fsz}
complete the proof of the Main Theorem.

\centerline{\bf Acknowledgments}

Jing acknowledges the National Natural Science Foundation of China grant 11531004
and Simons Foundation grant
523868. Liu acknowledges the National Natural Science Foundation of China grants 11531004 and 11701182.
Liu and Molev acknowledge the support of the Australian Research Council, grant DP180101825.

\bigskip\bigskip

\small

\noindent
N.J.:\qquad\qquad\qquad\qquad\\
Department of Mathematics\\
North Carolina State University, Raleigh, NC 27695, USA\\
jing@math.ncsu.edu

\vspace{3 mm}

\noindent
N.J. \& M.L.:\newline
School of Mathematical Sciences\\
South China University of Technology\\
Guangzhou, Guangdong 510640, China\\
mamliu@scut.edu.cn

\vspace{3 mm}

\noindent
M.L. \& A.M.:\newline
School of Mathematics and Statistics\newline
University of Sydney,
NSW 2006, Australia\newline
ming.liu2@sydney.edu.au\newline
alexander.molev@sydney.edu.au

\end{document}